\documentclass[aap,MSNbibl,seceqn,citesort,dvips]{arximspdf}
\usepackage{mathrsfs}

\doi{10.1214/11-AAP769}
\volume{22}
\issue{1}
\pubyear{2012}
\firstpage{285}
\lastpage{336}

\makeatletter

\renewcommand{\mid}{|}

\newcommand{\real}{\mathbb R}

\newcommand{\PP}{\mathbb{P}}
\newcommand{\E}{\mathbb{E}}
\newcommand{\Q}{\mathbb{Q}}
\newcommand{\de}{\delta}

\newtheorem{prop}{Proposition}[section]
\newtheorem{lemma}[prop]{Lemma}
\newtheorem{corollary}[prop]{Corollary}
\newtheorem{theorem}[prop]{Theorem}

\newproclaim{remark}[prop]{Remark}

\makeatother

\begin{document}
\begin{frontmatter}

\title{Differentiability of quadratic BSDEs generated by continuous martingales}
\runtitle{Differentiability of quadratic BSDEs}

\begin{aug}
\author[A]{\fnms{Peter} \snm{Imkeller}\ead[label=e1]{imkeller@mathematik.hu-berlin.de}},
\author[A]{\fnms{Anthony} \snm{R\'{e}veillac}\corref{}\thanksref{t1}\ead[label=e2]{areveill@mathematik.hu-berlin.de}} and
\author[A]{\fnms{Anja} \snm{Richter}\thanksref{t2}\ead[label=e3]{richtera@mathematik.hu-berlin.de}}
\runauthor{P. Imkeller, A. R\'{e}veillac and A. Richter}
\affiliation{Humboldt-Universit\"{a}t zu Berlin}
\address[A]{Institut f\"{u}r Mathematik\\
Humboldt-Universit\"{a}t zu Berlin\\
Unter den Linden 6\\
10099 Berlin\\
Germany\\
\printead{e1}\\
\hphantom{E-mail: }\printead*{e2}\\
\hphantom{E-mail: }\printead*{e3}} %adresu isvedimo komanda gale!
\end{aug}

\thankstext{t1}{Supported by DFG Research center \textsc{Matheon}, Project
E2.}

\thankstext{t2}{Supported by International
Research Training Group 1339 ``SMCP.''}

% HISTORY:
\received{\smonth{7} \syear{2009}}
\revised{\smonth{1} \syear{2011}}

% ABSTRACT
%
\begin{abstract}
In this paper we consider a class of BSDEs with drivers of quadratic
growth, on a stochastic basis generated by continuous local
martingales. We
first derive the Markov property of a forward--backward system (FBSDE) if
the generating martingale is a strong Markov process. Then we establish the
differentiability of a FBSDE with respect to the initial value of its
forward component. This enables us to obtain the main result of this
article, namely a representation formula for the control component of its
solution.
The latter is relevant in the context of securitization of random
liabilities arising from exogenous risk, which are optimally hedged by
investment in a given financial market with respect to exponential
preferences. In a purely stochastic formulation, the control process of the
backward component of the FBSDE steers the system into the random
liability and describes its optimal derivative hedge by investment in the
capital market, the dynamics of which is given by the forward component.
The representation formula of the main result describes this delta
hedge in
terms of the derivative of the BSDEs solution process on the one hand and
the correlation structure of the internal uncertainty captured by the
forward process and the external uncertainty responsible for the market
incompleteness on the other hand. The formula extends the scope of validity
of the results obtained by several authors in the Brownian setting. It is
designed to extend a~genuinely stochastic representation of the optimal
replication in cross hedging insurance derivatives from the classical
Black--Scholes model to incomplete markets on general stochastic bases. In
this setting, Malliavin's calculus which is required in the Brownian
framework, is replaced by new tools based on techniques related to a
calculus of quadratic covariations of basis martingales.\looseness=-1
\end{abstract}

% KEYWORDS
%

\end{frontmatter}

%s1 ###
\section{Introduction}
\label{introduction}

In recent years backward stochastic differential equations (BSDEs for
short) with drivers of quadratic growth have shown to be relevant in
several fields of application, for example, the study of properties of
PDEs (see,
e.g., \cite{NziOuknineSulem,BriandConfortola}). Closer to the
subject of this work, they were employed to provide a genuinely stochastic
approach to describe optimal investment strategies in a financial
market in
problems of hedging derivatives or liabilities of a small trader whose
business depends on market external risk. The latter scenario was addressed,
for instance, in
\cite
{HuImkellerMuller,AnkirchnerImkellerDosReis1,AnkirchnerImkellerDosReis2,Morlais1}.
A~small trader, such as an energy retailer, has a natural source of
income deriving
from his usual business. For instance, he may have a random position of
revenues from heating oil sales at the end of a heating season. To (cross)
hedge his risk arising from the partly market external uncertainty present
in the temperature process during the heating season, for example, via
derivatives written on temperature, he decides to invest in the capital
market, the inherent uncertainty of which is only correlated with this index
process. If the agent values his total income at terminal time by
exponential utility, or his risk by the entropic risk measure, he may be
interested in finding an optimal investment strategy that maximizes his
terminal utility, respectively, minimizes his total risk. The
description of such
strategies, even under convex constraints for the set of admissible ones,
is classical and may be achieved by convex duality methods and formulated
in terms of the analytic Hamilton--Jacobi--Bellman equation. In a genuinely
stochastic approach, \cite{HuImkellerMuller} interpreted the martingale
optimality principle by means of BSDEs with drivers of quadratic growth to
come up with a solution of this optimal investment problem even under
closed constraints that are not necessarily convex. The optimal investment
strategy is described by the control process in the solution pair of
such a
BSDE with an explicitly known driver. Using this approach, the authors of
\cite{AnkirchnerImkellerDosReis2} investigate utility indifference prices
and delta hedges for derivatives or liabilities written on nontradable
underlyings such as temperature in incomplete financial market models.
A~sensitivity analysis of the dependence of the optimal investment strategies
on the initial state of the Markovian forward process modeling the external
risk process provides an explicit delta hedging formula from the
representation of indifference prices in terms of forward--backward systems
of stochastic differential equations (FBSDEs). In the framework of a
Brownian basis, this analysis requires both the parametric as well as
variational differentiability in the sense of Malliavin calculus of the
solutions of the BSDE part (see
\cite{AnkirchnerImkellerDosReis1,AnkirchnerImkellerDosReis2,BriandConfortola}).
Related optimal investment problems have been investigated in
situations in
which the Gaussian basis is replaced by the one of a continuous martingale
(\cite{ManiaSchweizer} and \cite{Morlais1}, see also \cite{ElKarouiHuang}).

In this paper we intend to extend this utility indifference based explicit
description of a delta hedge to much more general stochastic bases. Our
main result will provide a probabilistic representation of the optimal
delta hedge of \cite{AnkirchnerImkellerDosReis2},\vadjust{\goodbreak} obtained there in the
Brownian setting, to more general scenarios in which pricing rules are
based on general continuous local martingales. We do this through a
sensitivity analysis of related systems of FBSDEs on a~stochastic basis
created by a continuous local martingale. As the backward component of our
system, we consider a BSDE of the form (\ref{eqn1}) driven by a continuous
local martingale $M$ with dynamics
%
%e1.1 ###
%
\begin{eqnarray}
\label{eqn1}
Y_t&=&B-\int_t^T Z_s \,dM_s + \int_t^T f(s,Y_s,Z_s) \,dC_s\nonumber\\[-8pt]\\[-8pt]
&&{} -
\int_t^T dL_s+ \frac{\kappa}{2} \int_t^T d\langle L, L \rangle
_s, \qquad t\in[0,T],\nonumber
\end{eqnarray}
where the \textit{generator} $f$ is assumed to be quadratic as a function
of $Z$, the \textit{terminal condition} $B$ is bounded, $C$ is an increasing
process defined as $C:=\arctan(\sum_i \langle M^{(i)}, M^{(i)}
\rangle)$, $L$ is a martingale
orthogonal to $M$ with quadratic variation $\langle L, L\rangle$ and
$\kappa$ is a positive constant.
A solution of (\ref{eqn1}) is given by a triplet $(Y,Z,L)$.
The forward component of our system is of the form
%
%e1.2 ###
%
\begin{equation}
\label{eqn2} X_s=x+\int_0^s \sigma(r,X_r,M_r) \,dM_r+ \int_0^s
b(r,X_r,M_r) \,dC_r, \qquad  s\in[0,T].\hspace*{-25pt}
\end{equation}
We first prove in Theorem \ref{theorem:Markovproperty} that the solution
processes $Y$ and $Z$ satisfy the Markov property, provided the terminal
condition $B$ is a smooth function of the terminal value of the forward
process (\ref{eqn2}) and that the local martingale~$M$ is a strong Markov
process. There is a subtlety in this setting which goes beyond causing a
purely technical complication, namely, that only the pair~$(X,M)$ is a
Markov process (as proved, e.g., in
\cite{ProtterEx,CinlarJacodProtterSharpe,Protter}). Only if $M$ has
independent increments it is
a stand-alone Markov process. We then show in Theorem \ref{DiffLast} that
the process $Y$ is differentiable with respect to the initial value of the
forward component (\ref{eqn2}) and that the derivatives of~$Y$ and $Z$
again satisfy a BSDE. The two properties then combine to allow us to state
and prove the main contribution of this paper. Thereby our delta hedge
representation (Theorem \ref{theorem:main}) generalizes the formula
obtained in the Gaussian setting (see
\cite{AnkirchnerImkellerDosReis2}, Theorem 6.7, for the quadratic case
and~\cite{ElKarouiPengQuenez}, Corollary~4.1, for the Lipschitz case).
More precisely, we show
that there exists a~deterministic function $u$ such that
%
%e1.3 ###
%
\begin{equation}
\label{eqn3}
Z_s=\partial_2 u(s,X_s,M_s) \sigma(s,X_s,M_s) + \partial_3 u(s,X_s,M_s),
\end{equation}
where $Y_s=u(s,X_s,M_s), s\in[0,T]$, and $\partial_i$ denotes the partial
derivative with respect to the $i$th variable (see Theorem
\ref{theorem:main}). In addition, we show that if $M$ has independent
increments and the coefficients of the forward process do not depend
on~$M$, then $Y_s=u(s,X_s)$ and equality (\ref{eqn3}) becomes
$Z_s=\partial_2
u(s,X_s) \sigma(s,X_s)$ which coincides with the formula known for the case
in which $M$ is a Brownian motion. To the best of our knowledge, relation
(\ref{eqn3}) is known only in the Brownian setting and the proof used in
the literature\vadjust{\goodbreak} relies on the representation of the stochastic process $Z$
as the trace of the Malliavin derivative $D$ (i.e., $Z_s=D_s Y_s,
s\in[0,T]$) relative to the underlying Brownian motion. Since Malliavin's
calculus is not available for general continuous local martingales, we
propose a new approach based on stochastic calculus techniques, in which
directional variational derivatives of Malliavin's calculus are
replaced by
absolute continuity properties of mixed variation processes of local basis
martingales. Note also reference~\cite{BallyPardouxStoica}, where a
Markovian representation of the solution $(Y,Z)$ of the solution of a BSDE
driven by a symmetric Markov process is given and whose driver is Lipschitz
in $z$ and satisfies a monotonicity condition in $y$ (see
\cite{BallyPardouxStoica}, Condition~(H2), page~35). However, the
representation of the
component $Z$ is not exactly similar to our representation (compare
\cite{BallyPardouxStoica}, Theorems 5.4 and \ref
{theorem:main}) due
to a~lack of regularity of the BSDEs driver in the setting of
\cite{BallyPardouxStoica}. Note finally that the method employed in
\cite{BallyPardouxStoica} relies on the calculus of Fukushima (see
\cite{FukushimaOshimaTakeda}) for symmetric Markov processes. We finally
emphasize that the local martingale $M$ considered in this paper is not
assumed to satisfy the martingale representation property.

The layout of this article is as follows. In Section \ref{prel} we state the main
notation and assumptions used in the paper. We discuss the Markov property
of an FBSDE in Section \ref{section:Markovproperty}. In Section
\ref{section:differentiability} we give sufficient conditions on
the FBSDEs to be differentiable in the initial values of its forward
component, while Section \ref{mainsection} is devoted to the representation formula
(\ref{eqn3}). Section \ref{sec6} is devoted to the finance and insurance application
of our main result.

%s2 ###
\section{Preliminaries}
\label{prel}
\subsection*{Notation}

Let $(M_t)_{t \in[0,T]}$ be a continuous $d$-dimensional local
martingale with $M_0=0$ which is defined on a
probability basis $(\Omega,\mathcal{F},(\mathcal{F}_t)_{t\in
[0,T]},\PP)$ where $T$ is a fixed positive real number.
We assume that the filtration $(\mathcal{F}_t)_{t\in[0,T]}$ is
continuous and complete so that every $\PP$-martingale is of the form
$Z \cdot M + L$, where $Z$ is a predictable $d$-dimensional process and
$L$ a~$\real$-valued martingale strongly orthogonal to $M$, that is,
$\langle L,M^{(i)}
\rangle=0$ for $i=1, \ldots, d$.
Here and in the following $M^{(i)}$, $i=1, \ldots,d$, denotes the entries
of the vector $M$. We assume that there exists a positive constant $Q$
such that
%
%e2.1 ###
%
\begin{equation}
\label{eq:boundedbrackets}
\bigl\langle M^{(i)}, M^{(j)} \bigr\rangle_T \leq Q\qquad  \forall1 \leq i, j
\leq d,  \qquad\PP\mbox{-a.s.}
\end{equation}
The Euclidean norm is denoted by \mbox{$|\cdot|$} and with $\mathcal E$ we refer
to the stochastic exponential.

From the Kunita--Watanabe inequality it follows that there exists a continuous,
adapted, bounded and increasing real-valued process $(C_t)_{t\in[0,T]}$
and a $\real^{d\times d}$-valued predictable process $(q_t)_{t\in[0,T]}$
such that the quadratic variation process $\langle M, M\rangle$ can be written
as
\[
\langle M, M\rangle_t=\int_0^t q_r q_r^\ast \,dC_r, \qquad  t \in[0,T],\vadjust{\goodbreak}
\]
where $^*$ denotes the transposition. We choose as in \cite{Morlais1},
$C:=\break\arctan( \sum_{i=1}^d \langle M^{(i)}$, $M^{(i)} \rangle
)$.
We write $\mathcal P$ for the predictable $\sigma$-field on $\Omega
\times[0,T]$.
Next we specify several spaces which we use in the sequel.
Given the arbitrary nonnegative and progressively measurable
real-valued process~$(\psi_t)_{t \in[0,T]}$, we define $\Psi$ by
$\Psi_t:=\int_0^t$, $\psi_s^2 \,dC_s, 0\le t\le T$.
For any $\beta>0$, $n \in\mathbb N$ and $p \in[1,\infty)$ we set:
\begin{itemize}
\item$\mathcal S^\infty:=\{ X\dvtx\Omega\times[0,T]\to\real\mid X$
adapted, bounded and continuous process$\}$,
\item$\mathcal S^p:=\{ X\dvtx\Omega\times[0,T]\to\real\mid X$
predictable process and \mbox{$\E[\sup_{t \in[0,T]} |X_t|^p]\,{<}\,\infty\}$},
\item$L^p(d\langle M,M \rangle\otimes d\PP)\vspace*{2pt}
:= \{Z\dvtx\Omega\times[0,T]\to\real^{1 \times d} \mid Z$
predictable process and
$\E[( \int_0^T |q_s Z_s^*|^2 \,dC_s )^{{p/2}}] <
\infty\}$,\vspace*{2pt}
\item$\mathcal M^2:= \{X\dvtx\Omega\times[0,T]\to\real\mid X$
square-integrable martingale$\}$,
\item$\mathbb{L}^\infty:=\{\xi\dvtx\Omega\to\real\mid \xi,
\mathcal{F}_T$-measurable bounded random variable$\}$,
\item$\mathbb{L}^p:=\{\xi\dvtx\Omega\to\real\mid \xi,
\mathcal{F}_T$-measurable random variable and $\E[|\xi|^p]
< \infty\}$,
\item$\mathbb L_\beta^2(\real^{n \times1}):=\{\xi\dvtx\Omega\to
\real\mid \xi, \mathcal{F}_T$-measurable random variable
and\break $\E[e^{\beta\Psi_T}|\xi|^2] < \infty\}$,
\item$\mathbb{H}^2_\beta:=
\{X\dvtx\Omega\times[0,T]\to\real^{d \times1} \mid
X$ predictable process and $\|X\|_\beta^2:=\break\E[\int_0^T
e^{\beta\Psi_t} | X_t |^2 \,dC_t]<\infty\}$,
\item$\mathbb{S}^2_\beta:=
\{X\dvtx\Omega\times[0,T]\to\real^{d \times1} \mid
X$ adapted continuous process and
$\|X\|_\beta^2:=\E[\sup_{t \in[0,T]} e^{\beta\Psi_t} |X_t |^2
]<\infty\}$.

\end{itemize}
Throughout\vspace*{1pt} this paper we will make use of the notation $(M^{t,m})_{s\in
[t,T]}$ ($t<T$, $m \in\real^{d \times1}$) which refers to the martingale
\[
M^{t,m}_s:=m+M_s-M_t
\]
defined with respect to the filtration $(\mathcal{F}_s^t)_{s\in
[t,T]}$ with $\mathcal{F}_s^t:=\sigma(\{M_u-M_t,  t\leq u \leq s\})$.
Obviously, all the preceding definitions can be introduced with
$M^{t,m}$ in place of $M$ and will inherit the superscript ${}^{t,m}$.
For convenience, we write $M^m:=M^{0,m}$.

Note that within this paper $c>0$ denotes a constant which can change
from line to line.

\subsection*{FBSDEs driven by continuous martingales}

In this subsection we present the main hypotheses needed in this
paper.
Let\vspace*{1pt} us fix $x \in\real^{n \times1}$ and $m \in\real^{d \times
1}$ and consider the process $X^{x,m}:=(X_t^{x,m})_{t\in
[0,T]}$ which is defined as a solution of the following stochastic
differential equation (SDE):
%
%e2.2 ###
%
\begin{eqnarray}
\label{SDE}
X_t^{x,m}&=&x+\int_0^t \sigma(s,X_s^{x,m},M_s^{m}) \,dM_s\nonumber\\[-8pt]\\[-8pt]
&&{} + \int_0^t
b(s,X_s^{x,m},M_s^{m}) \,dC_s, \qquad  t \in[0,T],\nonumber
\end{eqnarray}
where the coefficients $\sigma\dvtx[0,T] \times\real^{n\times1} \times
\real^{d\times1} \to\real^{n \times d}$ and $b\dvtx[0,T]\times
\real^{n\times1} \times\real^{d\times1} \to\real^{n \times1}$
are Borel-measurable functions.
By \cite{DoleansDade}, Theorem 1, and \cite{ProtterEx}, Theorem 3.1, this
SDE has a unique solution $X^{x,m} \in\mathcal S^p$ for all $p \geq1$
if the following hypothesis is satisfied.
\begin{longlist}[(H0)]
\item[(H0)] The functions $\sigma$ and $b$ are continuous
in $(s,x,m)$ and there exists a~$K>0$
such that for all $s \in[0,T]$, $x_1,x_2 \in\real^{n \times1}$ and
$m_1, m_2 \in\real^{d\times1}$
\begin{eqnarray*}
&&|\sigma(s,x_1,m_1) - \sigma(s,x_2,m_2)| + |b(s,x_1,m_1) -
b(s,x_2,m_2)| \\
&&\qquad\leq K (|x_1-x_2|+|m_1-m_2|).
\end{eqnarray*}
\end{longlist}
Next we give some properties of BSDEs which depend on the forward process~$X^{x,m}$.
More precisely we consider BSDEs of the form
%
%e2.3 ###
%
\begin{eqnarray}
\label{BSDE}
Y_t^{x,m}&=&F(X_T^{x,m},M_T^{m})-\int_t^T Z_r^{x,m} \,dM_r\nonumber\\
&&{} + \int_t^T
f(r,X_r^{x,m},M_r^m,Y_r^{x,m},Z_r^{x,m}q_r^*) \,dC_r\\
&&{}-\int_t^T dL_r^{x,m}
+ \frac{\kappa}{2} \int_t^T d\langle L^{x,m},L^{x,m}\rangle_r,
\qquad  t \in
[0,T],\nonumber
\end{eqnarray}
where $F\dvtx\real^{n \times1} \times\real^{d \times1} \to\real$ and
$f\dvtx\Omega\times[0,T]\times\real^{n \times1} \times\real^{d\times
1} \times\real
\times\real^{1\times d} \to\real$ are $\mathcal B(\real^{n \times
1})$, respectively, $\mathcal P \otimes\mathcal B(\real^{n \times1})
\otimes\mathcal B(\real^{d \times1}) \otimes\mathcal B(\real)
\otimes\mathcal B(\real^{1 \times d})$-measurable functions.
By $\mathcal B(\real^{d})$ we denote the Borel $\sigma$-algebra.
A solution of the BSDE with \textit{terminal condition}
$F(X_T^{x,m},M_T^{m})$, a constant $\kappa$ and \textit{generator} $f$
is defined to be a triple of processes $(Y^{x,m},Z^{x,m},L^{x,m}) \in
\mathcal S^\infty\times L^2(d \langle M,M \rangle\otimes d\PP) \times
\mathcal M^2$ satisfying (\ref{BSDE}) and such that $\langle
L^{x,m},M^i \rangle=0$, $i=1,\ldots,d$, and $\PP$-a.s. $\int_0^T
|f(r,X_r^{x,m},M_r^m,Y_r^{x,m},Z_r^{x,m}q_r^*)|\,dC_r < \infty$.

Let $\mathcal{V}:=\real^{n \times1} \times\real^{d\times1} \times
\real\times\real^{1\times d}$ and assume that (H0) holds.
Furthermore, we define the measure $\nu(A)= E[\int_0^T \mathbf
1_{A}(s) \,dC_s ] $ for all \mbox{$A \in\mathcal B([0,T]) \otimes\mathcal F$}.
Under the following conditions, existence and uniqueness of a solution
of the backward equation (\ref{BSDE}) was recently discussed in
\cite{Morlais1}, Theorem 2.5:
\begin{longlist}[(H1)]
\item[(H1)] The function $F$ is bounded.
\item[(H2)] The generator $f$ is continuous in $(y,z)$ and there
exists a nonnegative predictable process $\eta$ such that $\int_0^T
\eta_s
\,dC_s \leq a$, where $a$ is a positive constant as well as positive numbers
$b$ and $\gamma$, such that $\nu$-a.e.
\[
\vert f(s,x,m,y,z) \vert\leq\eta_s + b \eta_s \vert y \vert+\frac
{\gamma}{2}
| z |^2\qquad  \mbox{with } \gamma\geq\vert\kappa\vert,
\gamma\geq
b,  (x,m,y,z) \in\mathcal{V}.
\]
\end{longlist}
An additional assumption is needed to obtain uniqueness (see
\cite{Morlais1}, Theorem~2.6).
\begin{longlist}[(H3)]
\item[(H3)] For every $\beta\geq1$ we
have $\int_0^T \vert f(s,0,0,0,0) \vert \,dC_s \in L^\beta(\PP)$. In
addition, there exist two constants $\mu$ and $\nu$, a
nonnegative\vadjust{\goodbreak}
predictable process $\theta$ satisfying $\int_0^T |q_s \theta_s|^2 \,dC_s
\leq c_\theta$ ($c_\theta\in\real$), such that $\nu$-a.e.
\begin{eqnarray*}
&&(y_1-y_2)\bigl(f(s,x,m,y_1,z)-f(s,x,m,y_2,z)\bigr)\\
&&\qquad \leq\mu|y_1-y_2|^2, \qquad (x,m,y_i,z)
\in\mathcal{V},  i=1,2,
\end{eqnarray*}
and
\begin{eqnarray*}
&&|f(s,x,m,y,z_1)-f(s,x,m,y,z_2)| \\
&&\qquad\leq\nu(|q_s \theta_s| + |z_1| +
|z_2|)|z_1-z_2|,\qquad (x,m,y,z_i) \in\mathcal{V},  i=1,2.
\end{eqnarray*}
\end{longlist}
In this paper we will deal with martingales of bounded mean oscillation,
briefly called BMO martingales.
We recall that $ Z \cdot M$ is a BMO martingale if and only if
\[
\|Z \cdot M\|_{\mathrm{BMO}_2}= \sup_{\tau\leq T} \E\biggl[\int_{\tau}^T
|q_s Z_s^*|^2 \,dC_s \Big| \mathcal F_{\tau} \biggr]^{1/2} <
\infty,
\]
where the supremum is taken over all stopping times $\tau\leq T$.
We refer the reader to \cite{Kazamaki} for a survey. Specifically we
need the following hypothesis.
\begin{longlist}[(H4)]
\item[(H4)] There exist a $\real^{1 \times d}$-valued
predictable process $K$ and a constant $\alpha\in(0,1)$ such that $K
\cdot M$ is a BMO martingale satisfying $\nu$-a.e.
\[
(y_1-y_2)\bigl(f(s,x,m,y_1,z)-f(s,x,m,y_2,z)\bigr) \leq\vert q_s K_s^* \vert
^{2\alpha}|y_1-y_2|^2
\]
for all $(x,m,y_i,z) \in\mathcal{V}$, $i=1,2$, and
\[
|f(s,x,m,y,z_1)-f(s,x,m,y,z_2)| \leq\vert q_s K_s^* \vert|z_1-z_2|
\]
for all $(x,m,y,z_i) \in\mathcal{V}$, $i=1,2$.
\end{longlist}
Throughout this paper we also consider a second type of BSDEs associated
with the forward process $X^{x,m}$ solving (\ref{SDE}), that is,
%
%e2.4 ###
%
\begin{eqnarray}
\label{BSDEohne}
U^{x,m}_t&=& F(X^{x,m}_T,M_T^{m}) - \int_t^T V^{x,m}_s \,dM_s\nonumber\\
&&{} + \int_t^T
f(s,X_s^{x,m},M^m_s,U^{x,m}_s,V^{x,m}_sq_s^*) \,dC_s\\
&&{} + \int_t^T
dN^{x,m}_s,\nonumber
\end{eqnarray}
$t \in[0,T]$, where $N$ is a square-integrable martingale.
This type of BSDE has been studied by El Karoui and Huang in
\cite{ElKarouiHuang}. Under the following assumptions on terminal
condition $F(X_T^{x,m},M_T^{m})$
and generator $f$, there exists a~unique solution
$(U^{t,x,m},V^{t,x,m},N^{t,x,m}) \in\mathbb S^2_{\beta} \times
\mathbb H^2_{\beta} \times\mathcal M^2$ to the BSDE~(\ref{BSDEohne}):
\begin{longlist}[(L2)]
\item[(L1)] The function $F$ satisfies $F(
X_T^{x,m},M_T^{m}) \in\mathbb L_\beta^2(\real^{n \times1} \times
\real^{d \times1})$ for some~$\beta$ large enough.\vadjust{\goodbreak}

\item[(L2)] The generator $f$ satisfies $\nu$-a.e.
\begin{eqnarray*}
&&\vert f(s,x,m,y_1,z_1) -f(s,x,m,y_2,z_2)\vert\\
&&\qquad\leq r_s |y_1-y_2| +
\theta_s |z_1-z_2|,  \qquad (x,m,y_i,z_i) \in\mathcal{V}, i=1,2,
\end{eqnarray*}
where $r$ and $\theta$ are two nonnegative predictable processes.
Let $\alpha_s^2=r_s+\theta_s^2$. We assume $\nu$-a.e. that $\alpha
_s^2>0$ and $\frac{f(\cdot,0,0)}{\alpha} \in\mathbb H_\beta^2$ for
some $\beta>0$ large enough.
\end{longlist}
We conclude this section by presenting assumptions which
will be useful in Section \ref{section:differentiability}, where we
find sufficient conditions for FBSDEs to be differentiable in their
initial values $(x,m) \in\real^{n \times1} \times\real^{d \times
1}$.
Given a function $g\dvtx[0,T] \times\real^{n \times1} \times\real^{d
\times1} \to\real$ we denote the partial derivatives with respect to
the $i$th variable by $\partial_i g(s,x,m)$ and, if no confusion can
arise, we write $\partial_2 g(s,x,m):=(\partial_{1+j}g(s,x,m))_{j=1,
\ldots, n}$ and $\partial_3 g(s,x,m):=(\partial
_{1+n+j}g(s,x,m))_{j=1, \ldots, d}$.
\begin{longlist}[(D4)]
\item[(D1)] The coefficients $\sigma$ and $b$ have locally
Lipschitz partial derivatives in~$x$ and $m$ uniformly in time.
\item[(D2)] The functions $F$ and $\nabla F$ are globally Lipschitz.
\item[(D3)] The generator $f$ is differentiable in $x,m,y$
and $z$ and there exist a constant $C>0$ and a
nonnegative predictable process $\theta$ satisfying $\int_0^T |q_s
\theta_s|^2 \,dC_s \leq c_\theta$ ($c_\theta\in\real$),
such that the partial derivatives satisfy $\nu$-a.e.
\[
|\partial_i f(s,x,m,y,z)| \leq C (|q_s \theta_s| + |z|),\qquad
(x,m,y,z) \in\mathcal{V}, i=2, \ldots, 5.
\]

\item[(D4)] The generator $f$ is differentiable in $x,m,y$
and $z$ and there exist a constant $C>0$ and a
nonnegative predictable process $\theta$ satisfying $\int_0^T |q_s
\theta_s|^2 \,dC_s \leq c_\theta$ ($c_\theta\in\real$), such that
the partial derivative $\partial_5 f$ is Lipschitz in $(x,m,y,z)$
and for all $i=2, \ldots, 4$ the following inequality holds $\nu$-a.e.:
\begin{eqnarray*}
&& |\partial_i f(s,x_1,m_1,y_1,z_1)-\partial_i f(s,x_2,m_2,y_2,z_2)|\\
&&\qquad \leq C (|q_s \theta_s| + \vert z_1 \vert+\vert z_2 \vert) (\vert
x_1-x_2 \vert+\vert m_1-m_2\vert+\vert y_1-y_2 \vert+\vert z_1-z_2
\vert)
\end{eqnarray*}
for all $(x_j,m_j,y_j,z_j) \in\mathcal{S}, j=1,2$.
\end{longlist}

%s3 ###
\section{The Markov property of FBSDEs}
\label{section:Markovproperty}

For a fixed initial time $t \in[0,T)$ and initial values $x \in\real
^{n \times1}$ and $m \in\real^{d \times1}$ we consider a SDE of the form
%
%e3.1 ###
%
\begin{eqnarray}
\label{SDEbis}
X_s^{t,x,m}&=&x+\int_t^s \sigma(u,X_u^{t,x,m},M_u^{t,m}) \,dM_u \nonumber\\[-8pt]\\[-8pt]
&&{}+ \int
_t^s b(u,X_u^{t,x,m},M_u^{t,m}) \,dC_u,  \qquad s \in[t,T],\nonumber
\end{eqnarray}
where $M$ is a local martingale as in Section \ref{prel} with values
in $\real^{d \times1}$, $\sigma\dvtx[0,T] \times\real^{n\times
1} \times\real^{d\times1} \to\real^{n \times d}$ and
$b\dvtx[0,T]\times\real^{n\times1} \times\real^{d\times1} \to
\real^{n \times1}$.
Throughout this chapter the coefficients $\sigma$ and $b$ satisfy (H0)
and hence, (\ref{SDEbis}) has a
unique solution $X^{t,x,m}$.
Before stating and proving the main results of this section we recall
the following proposition which is a combination of
\cite{CinlarJacodProtterSharpe}, Theorem~(8.11)
(see also~\cite{Protter}, Theorem V.35) and \cite{Protter77}, Theorem 5.3.
\begin{prop}
\label{CinlaretalProtter}
\textup{(i)}
If $M$ is a strong Markov process then
$(M_s^{t,m},\allowbreak X_s^{t,x,m})_{s\in[t,T]}$ is a strong Markov process.

\textup{(ii)} If $M$ is a strong Markov process with independent
increments and if the coefficients $\sigma$ and $b$ do not depend on
$M$, that is to say
\[
X_s^{t,x}=x+\int_t^s \sigma(u,X_u^{t,x}) \,dM_u + \int_t^s b(u,X_u^{t,x})
\,dC_u,
\]
then the process $(X_s^{t,x})_{s\in[t,T]}$ itself is a strong Markov process.
\end{prop}

Note that in
\cite{AnkirchnerImkellerDosReis1,AnkirchnerImkellerDosReis2,ElKarouiPengQuenez}
the martingale
considered is a standard Brownian motion so that situation (ii) of Proposition
3.1 applies.
In fact, this case presents at least two major advantages; first, the process
$X$ is a Markov process itself and second, the quadratic variation of $M$
is deterministic.

This section is organized as follows. We first prove in
Proposition \ref{MarkovpropertyLip} that the solution of a Lipschitz
BSDE associated to a forward SDE of the form~(\ref{SDEbis}) is
already determined by the solution $X^{t,x,m}$ of (\ref{SDEbis}) and
the Markov process $M^{t,m}$. In Theorem
\ref{theorem:Markovproperty} we then extend this result to quadratic
BSDEs.

Consider a BSDE of the form
%
%e3.2 ###
%
\begin{eqnarray}
\label{BSDElip}
U_s^{t,x,m}&=&F(X_T^{t,x,m},M_T^{t,m})-\int_s^T V_r^{t,x,m} \,dM_r\nonumber\\
&&{} + \int
_s^T f(r,X_r^{t,x,m},M_r^{t,m},U_r^{t,x,m},V_r^{t,x,m}q^*_r)
\,dC_r \\
&&{}- \int_s^T dN^{t,x,m}_r,  \qquad s \in[t,T].\nonumber
\end{eqnarray}
We suppose that the driver does not depend on $\Omega$ and hence, is a
deterministic Borel measurable function
$f\dvtx[0,T]\times\real^{n \times1} \times\real^{d \times1} \times
\real\times\real^{1 \times d} \to\real$.
If~$F$ and $f$ satisfy hypotheses (L1) and (L2) then the BSDE (\ref
{BSDElip}) admits a~unique solution $(U^{t,x,m},V^{t,x,m},N^{t,x,m})
\in\mathbb S^2_{\beta} \times\mathbb H^2_{\beta} \times\mathcal
M^2$ (see~\cite{ElKarouiHuang}, Theorem~6.1).
By $\mathcal B_e(\real^{n \times1} \times\real^{d \times1})$ we
denote the $\sigma$-algebra generated by the family of functions
$(x,m) \mapsto\E[ \int_t^T \phi(s,X_s^{t,x,m},M_s^{t,m}) \,dC_s
]$, where $\phi\dvtx \Omega\times[0,T] \times\real^{n \times1} \times
\real^{d \times1} \to\real$ is predictable, continuous and bounded.
\begin{prop}
\label{MarkovpropertyLip} Assume that $M$ is a strong Markov
process and that~\textup{(L1)} and~\textup{(L2)} are in force. Then there exist
deterministic functions $u\dvtx[0,T]\times\real^{n \times1}\times
\real^{d \times1} \to\real$, $\mathcal{B}([0,T])\otimes
\mathcal{B}_e(\real^{n \times1}\times\real^{d \times
1})$-measurable and $v\dvtx[0,T]\times\real^{n\times1}\times\real^{d
\times1} \to\real^{1\times d}$, $\mathcal{B}([0,T])\otimes
\mathcal{B}_e(\real^{n \times1}\times\real^{d \times
1})$-measurable such that
%
%e3.3 ###
%
\begin{eqnarray}
\label{EqMarkovPropertyLip}
U_s^{t,x,m}&=&u(s,X_s^{t,x,m},M_s^{t,m}),\nonumber\\[-8pt]\\[-8pt]
V_s^{t,x,m}&=&v(s,X_s^{t,x,m},M_s^{t,m}) , \qquad  s \in[t,T].\nonumber
\end{eqnarray}
\end{prop}
\begin{remark}
\label{remarkIndInc} Before turning to the proof of Proposition
\ref{MarkovpropertyLip} we stress the following point.
Assume $M$ and $X$ are as in Proposition \ref{CinlaretalProtter}(ii)
and that the driver $f$ in (\ref{BSDElip}) does not depend on $M$,
then Proposition
\ref{MarkovpropertyLip} is equivalent to the existence of
deterministic functions
$u\dvtx[0,T]\times\real^{n \times1} \to\real$, $\mathcal
{B}([0,T])\otimes
\mathcal{B}_e(\real^{n \times1})$-measurable and $v\dvtx[0,T]\times
\real^{n\times
1} \to\real^{1\times d}$, $\mathcal{B}([0,T])\otimes\mathcal
{B}_e(\real^{n
\times1})$-measurable such that
\[
U_s^{t,x}=u(s,X_s^{t,x}), \qquad  V_s^{t,x}=v(s,X_s^{t,x}), \qquad  s \in[t,T].
\]
\end{remark}
\begin{pf*}{Proof of Proposition \ref{MarkovpropertyLip}}
Consider the following sequence
$(U^{k,t,x,m}$, $V^{k,t,x,m},N^{k,t,x,m})_{k\ge0}$ of BSDEs:
%
%e3.4 ###
%
\begin{eqnarray}
\label{approximinter}
U^{0,t,x}&=&V^{0,t,x}=0,\nonumber\\
U_s^{k+1,t,x}&=&F(X_T^{t,x,m},M_T^{t,m})\nonumber\\[-8pt]\\[-8pt]
&&{}+\int_s^T
f(r,X_r^{t,x,m},M_r^{t,m},U_r^{k,t,x,m},V_r^{k,t,x,m}q_r^*)
\,dC_r\nonumber\\
&&{} - \int_s^T V_r^{k+1,t,x,m} \,dM_r - \int_s^T
dN^{k+1,t,x,m}_r.\nonumber
\end{eqnarray}
We recall an estimate obtained in \cite{ElKarouiHuang}, page 35.
Let $\alpha>0$ and $\beta>0$ be as in Section \ref{prel}.
Then
\begin{eqnarray*}
&&\|\alpha(U^{k+1,t,x,m}-U^{k,t,x,m})\|_\beta^2\\
&&\quad{} + \|q
(V^{k+1,t,x,m}-V^{k,t,x,m})^{\ast}\|_\beta^2+ \|
(N^{k+1,t,x,m}-N^{k,t,x,m})\|_\beta^2 \\
&&\qquad\leq\varepsilon\bigl(\|\alpha(U^{k,t,x,m}-U^{k-1,t,x,m})\|_\beta
^2+ \|q (V^{k,t,x,m}-V^{k-1,t,x,m})^{\ast}\|_\beta^2\\
&&\qquad\quad\hspace*{147.3pt}{} + \|
(N^{k,t,x,m}-N^{k-1,t,x,m})\|_\beta^2\bigr),
\end{eqnarray*}
where $\varepsilon$ is a constant depending on $\beta$ which can be chosen
with $\varepsilon<1$.
Applying the result recursively we obtain
\begin{eqnarray*}
&&\|\alpha(U^{k+1,t,x,m}-U^{k,t,x,m})\|_\beta^2\\
&&\quad{}+ \|q
(V^{k+1,t,x,m}-V^{k,t,x,m})^{\ast}\|_\beta^2+\|
(N^{k+1,t,x,m}-N^{k,t,x,m})\|_\beta^2\\
&&\qquad\leq\varepsilon^k \bigl(\|\alpha(U^{1,t,x,m}-U^{0,t,x,m})\|_\beta
^2+ \|q (V^{1,t,x,m}-V^{0,t,x,m})^{\ast}\|_\beta^2\\
&&\qquad\quad\hspace*{140pt}{}+\|
(N^{1,t,x,m}-N^{0,t,x,m})\|_\beta^2\bigr).
\end{eqnarray*}
Since
\begin{eqnarray*}
&&\sum_{k=0}^\infty\|\alpha
(U^{k+1,t,x,m}-U^{k,t,x,m})\|_\beta^2+ \|q
(V^{k+1,t,x,m}-V^{k,t,x,m})^{\ast}\|_\beta^2\\
&&\qquad{}+\|
(N^{k+1,t,x,m}-N^{k,t,x,m})\|_\beta^2<\infty,
\end{eqnarray*}
the sequence
$(U^{k,t,x,m},V^{k,t,x,m},N^{k,t,x,m})_k$ converges $\nu$-a.e. to
$(U^{t,x,m},V^{t,x,m}$, $N^{t,x,m})$ as $k$ tends to
infinity.

We show by induction on $k\geq1$ the following property (\textit{Prop$_k$}):
\begin{property*}
There exist deterministic functions
$\Phi^{k}\dvtx[0,T]\times\real^{n \times1} \times\real^{d \times1}
\to\real$, $\mathcal{B}([0,T])\otimes\mathcal{B}_e(\real^{n
\times
1}\times\real^{d \times1})$-measurable and $\Psi^{k}\dvtx[0,T]\times
\real^{n\times1}\times\real^{d \times1} \to\real^{1\times d}$,
$\mathcal{B}([0,T])\otimes\mathcal{B}_e(\real^{n \times1}\times
\real^{d \times1})$-measurable such that
$U_s^{k,t,x,m}=\Phi^{k}(s$, $X_s^{t,x,m},M_s^{t,m})$ and
$V_s^{k,t,x,m}=\Psi^{k}(s,X_s^{t,x,m},M_s^{t,m}) $, for $t\le s\le
T, k\in\mathbb{N}$.
\end{property*}
\begin{pf*}{Proof of (\textit{Prop$_1$})} From the definition of
$U^{1,t,x,m}$ and since $N^{1,t,x,m}$ is a martingale, we have for
$s\in[t,T]$
%
%e3.5 ###
%
\begin{eqnarray}
\label{P1.1}\qquad
U_s^{1,t,x,m}&=&\E[U_s^{1,t,x,m} \vert\mathcal{F}_s^t]
\nonumber\\[-8pt]\\[-8pt]
&=&\E\biggl[F(X_T^{t,x,m},M_T^{t,m})-\int_s^T
f(r,X_r^{t,x,m},M_r^{t,m},0,0) \,dC_r \Big|
\mathcal{F}_s^t \biggr]. \nonumber
\end{eqnarray}
The Markov property and Doob--Dynkin's lemma give
\begin{eqnarray*}
U_s^{1,t,x,m}
&=&\E\biggl[F(X_T^{t,x,m},M_T^{t,m})-\int_s^T
f(r,X_r^{t,x,m},M_r^{t,m},0,0) \,dC_r \Big| \mathcal{F}_s^t
\biggr]
\\
&=&\E\biggl[F(X_T^{t,x,m},M_T^{t,m})\\
&&\hphantom{\E\biggl[}{}-\int_s^T
f(r,X_r^{t,x,m},M_r^{t,m},0,0) \,dC_r \Big|
(X^{t,x,m}_s,M_s^{t,m}) \biggr]
\\
&=&\Phi^1(s,X_s^{t,x,m},M_s^{t,m}),
\end{eqnarray*}
where $\Phi^{1}\dvtx[0,T]\times\real^{n \times1} \times\real^{d
\times1} \to\real$. Now let
\[
R_s^{1,t,x,m}=U_s^{1,t,x,m}+\int_t^s f(r,X_r^{t,x,m},M_r^{t,m},0,0)
\,dC_r, \qquad  s \in[t,T].
\]
Then for $s \in[t,T]$
%
%e3.6 ###
%
\begin{equation}
\label{QuadraVaria1}
R_s^{1,t,x,m}=\int_t^s V_r^{1,t,x,m} \,dM_r+N^{t,x,m}_s-N^{t,x,m}_t,
\end{equation}
hence, using the localization technique, we can assume that
$R^{1,t,x,m}$ is a~stron\-gly additive (in the sense of
\cite{CinlarJacodProtterSharpe}, page 169) square integrable
martingale. Now we apply \cite{CinlarJacod}, Theorem (2.16), to
$\mathcal{Y}^1:=M$ and $\mathcal{Y}^2:=R$. Thus, there exist two
additive locally square integrable martingales $\mathcal{M}^1$ and
$\mathcal{M}^2$, two deterministic functions $\Psi^{1},\Psi
^{2}\dvtx[0,T]\times\real^{n\times1}\times
\real^{d \times1} \to\real^{1\times d}$ such that \mbox{$\mathcal
{Y}^1=\mathcal{M}^1$} and $\mathcal{Y}^2_s=\int_t^s \Psi
^{1}(s,X_s^{t,x,m},M_s^{t,m}) \,d\mathcal{M}_s^1 + \int_t^s \Psi
^{2}(s,X_s^{t,x,m},M_s^{t,m}) \,d\mathcal{M}_s^2$. By definition\vspace*{1pt} of $R$
we deduce that $\mathcal{M}^2$ has to be equal to $N^{1,t,x,m}$
(showing that $N^{1,t,x,m}$ is additive) and that $\Psi^{2}\equiv1$.
This shows that
\[
V_s^{1,t,x,m}=\Psi^{1}(s,X_s^{t,x,m},M_s^{t,m}) ,\qquad  \nu\mbox{-a.e.}
\]
Letting $k\geq1$, we prove (\textit{Prop$_k$}) $\Longrightarrow$
(\textit{Prop$_{k+1}$}).
For $s \in[t,T]$ we have
\begin{eqnarray*}
U_s^{k+1,t,x,m}&=&
\E[U_s^{k+1,t,x,m} \vert\mathcal{F}_s^t]
\\
&=&\E\biggl[F(X_T^{t,x,m},M_T^{t,m})\\[-2pt]
&&\quad\hphantom{\E\biggl[}{}
-\int_s^T
f(r,X_r^{t,x,m},M_r^{t,m},U_r^{k,t,x,m},V_r^{k,t,x,m}q_r^*) \,dC_r
\Big| \mathcal{F}_s^t \biggr]
\\[-2pt]
&=&\E\biggl[F(X_T^{t,x,m},M_T^{t,m})
\\[-2pt]
&&\quad\hphantom{\E\biggl[}{}
-\int_s^T f(r,X_r^{t,x,m},M_r^{t,m},\Phi
^{k}(r,X_r^{t,x,m},M_r^{t,m}),\\[-2pt]
&&\hspace*{136.1pt}\Psi^{k}(r,X_r^{t,x,m},M_r^{t,m})q_r^*)
\,dC_r \Big| \mathcal{F}_s^t \biggr]
\\[-2pt]
&=&\E\biggl[F(X_T^{t,x,m},M_T^{t,m})-\int_s^T
f^{k}(r,X_r^{t,x,m},M_r^{t,m}) \,dC_r \Big| \mathcal{F}_s^t\biggr],
\end{eqnarray*}
where $f^{k}(r,y,z):=f(r,y,\Phi^{k}(r,y,z),\Psi^{k}(r,y,z)q_r^*)$.
Using the same argument as in the case $k=1$, we deduce that there
exists a
function $\Phi^{k+1}\dvtx[0,T]\times\real^{m\times1}\times\real^{d
\times1} \to\real$ such that
\[
U_s^{k+1,t,x,m}=\Phi^{k+1}(s,X_s^{t,x,m},M_s^{t,m}).
\]
For $s \in[t,T]$ let
\begin{eqnarray*}
R_s^{k+1,t,x,m}&=&U_s^{k+1,t,x,m}+\int_t^s
f^{k}(r,X_r^{t,x,m},M_r^{t,m}) \,dC_r\\
&&{} -N_s^{k+1,t,x,m} +N_t^{k+1,t,x,m}.
\end{eqnarray*}
Following the same procedure as before, we deduce that there
exists a function $\Psi^{k+1}\dvtx[0,T]\times\real^{n\times1}\times
\real^{d
\times1} \to\real^{1 \times d}$ such that
\[
V_s^{k+1,t,x,m}=\Psi^{k+1}(s,X_s^{t,x,m},M_s^{t,m}).
\]
Let
\[
u(r,y,z):=\limsup_{k \to\infty} \Phi^{k}(r,y,z),\qquad
v(r,y,z):=\limsup_{k \to\infty} \Psi^{k}(r,y,z).
\]
Since the sequence $(U^{k,t,x},V^{k,t,x},N^{t,x,m})_k$ converges
$\nu$-a.e. to $(U^{t,x,m},V^{t,x,m}$, $N^{t,x,m})$
as $k$ tends to infinity, we have for $s \in[t,T]$
\begin{eqnarray*}
u(s,X_s^{t,x,m},M_s^{t,m})
&=& \Bigl(\limsup_{k \to\infty} \Phi
^{k}\Bigr)(s,X_s^{t,x,m},M_s^{t,m})\\
&=&\limsup_{k \to\infty} (\Phi
^{k}(s,X_s^{t,x,m},M_s^{t,m}))
\\
&=&\limsup_{k \to\infty} U_s^{k,t,x,m}=U_s^{t,x,m}.
\end{eqnarray*}
Similarly we obtain
\[
v(s,X_s^{t,x,m},M_s^{t,m}) =V_s^{t,x,m}.
\]
\upqed\end{pf*}
\noqed\end{pf*}

We conclude this section by extending Proposition \ref
{MarkovpropertyLip} to a quadratic FBSDE.
More precisely, we consider the following BSDE:
%
%e3.7 ###
%
\begin{eqnarray}
\label{BSDEbis}
Y_s^{t,x,m}&=&F(X_T^{t,x,m},M_T^{t,m})-\int_s^T Z_u^{t,x,m} \,dM_u \nonumber\\
&&{}+ \int
_s^T f(u,X_u^{t,x,m},M_u^{t,m},Y_u^{t,x,m},Z_u^{t,x,m}q_u^*) \,dC_u
\\
&&{}- \int_s^T dL^{t,x,m}_u + \frac{\kappa}{2} \int_s^T d\langle
L^{t,x,m}, L^{t,x,m} \rangle_u, \qquad  s \in[t,T],\nonumber
\end{eqnarray}
where the forward process $X^{t,x,m}$ is a solution of
(\ref{SDEbis}). Again we suppose that the driver $f$ does not depend
on $\Omega$ and hence, is a deterministic Borel measurable function
$f\dvtx[0,T]\times\real^{n \times1} \times\real^{d \times1} \times
\real\times\real^{1 \times d} \to\real$. If $F$ satisfies~(H1)
and $f$ hypotheses (H2) and (H3), then the BSDE (\ref{BSDEbis})
admits a~unique solution $ (Y^{t,x,m},Z^{t,x,m},L^{t,x,m}) \in
\mathbb S^{\infty} \times L^2(d\langle M,M\rangle\otimes d\PP)
\times\mathcal M^2$ (see~\cite{Morlais1}, Theorem~2.5).
\begin{theorem}
\label{theorem:Markovproperty} We assume that $M$
is a strong Markov process and that~\textup{(H1)--(H3)} hold. Then there
exist deterministic functions $u\dvtx[0,T]\times\real^{n \times1}
\times\real^{d \times1} \to\real$, $\mathcal{B}([0,T])\otimes
\mathcal{B}_e(\real^{n \times1}\times\real^{d \times
1})$-measurable and $v\dvtx[0,T]\times\real^{n\times1} \times\real^{d
\times1} \to\real^{1\times d}$, $\mathcal{B}([0,T])\otimes
\mathcal{B}_e(\real^{n \times1}\times\real^{d \times
1})$-measurable such that
%
%e3.8 ###
%
\begin{eqnarray}
\label{EqMarkovPropertyQuad}
Y_s^{t,x,m}&=&u(s,X_s^{t,x,m},M_s^{t,m}),\nonumber\\[-8pt]\\[-8pt]
Z_s^{t,x,m}&=&v(s,X_s^{t,x,m},M_s^{t,m}) , \qquad  s \in[t,T].\nonumber
\end{eqnarray}
\end{theorem}
\begin{remark}
As mentioned in Remark \ref{remarkIndInc}, in the framework of Proposition
\ref{CinlaretalProtter}(ii), when the driver $f$ in
(\ref{BSDEbis})
does not depend on $M$, Theorem \ref{theorem:Markovproperty}
simplifies to
the existence of deterministic functions $u\dvtx[0,T]\times\real^{n
\times1}
\to\real$, $\mathcal{B}([0,T])\otimes\mathcal{B}_e(\real^{n
\times1})$-measurable
and $v\dvtx[0,T]\times\real^{n\times1} \to\real^{1\times d}$,
$\mathcal{B}([0,T])\otimes
\mathcal{B}_e(\real^{n \times1})$-measurable such that
\[
Y_s^{t,x}=u(s,X_s^{t,x}),\qquad
Z_s^{t,x}=v(s,X_s^{t,x}), \qquad  s \in[t,T].
\]
\end{remark}
\begin{pf*}{Proof of Theorem \ref{theorem:Markovproperty}}
Existence and uniqueness of the solution of~(\ref{BSDEbis}) under
the hypotheses (H1)--(H3) have been obtained in \cite{Morlais1},
Theorems~2.5 and~2.6. More precisely, it is shown in the proof of
\cite{Morlais1}, Theorem~2.5, that the solution of a quadratic BSDE
can be derived as the limit of solutions of a sequence of BSDEs with
Lipschitz generators. We follow this proof and begin by relaxing
condition~(H2).
Indeed, consider the following assumption~(H2$'$) where the generator $f$
does not need to be bounded in $y$ anymore.
\begin{longlist}[(H2$'$)]
\item[(H2$'$)] The generator $f$ is continuous in $(y,z)$ and there
exists a predictable process $\eta$ such that $\eta\geq0$ and $\int
_0^T \eta_s \,dC_s \leq a$,
where $a$ is a positive constant. Furthermore, there exists a constant
$\gamma>0$ such that $\nu$-a.e.
\[
\vert f(s,x,m,y,z) \vert\leq\eta_s +\frac{\gamma}{2} |z|^2\qquad
\mbox{with } \gamma\geq|\kappa|,  (x,m,y,z) \in\mathcal V.
\]
\end{longlist}
Assume that one can prove existence of a solution of (\ref{BSDEbis})
if $f$ satisfies (H2$'$) instead of (H2). Let $f_K$ be the generator
$f$ truncated in $Y$ at level $K$ (as in~\cite{Morlais1}, Lemma 3.1).
More precisely, set $f_K(s,x,m,y,z):=f(s,x,m,\rho(y)_K,z)$ with
\[
\rho_K(y):=
\cases{ -K, &\quad  if $y<-K$,\cr
y, &\quad  if $\vert y \vert\leq K$,\cr
K, &\quad  if $y>K$.}
\]
It is shown in \cite{Morlais1}, proof of Theorem 2.5, Step 1, that
$f_K$ satisfies (H2$'$).
Hence, by hypothesis, there exists a triple of stochastic processes
$(Y^{t,x,m}_K,\allowbreak Z^{t,x,m}_K,L^{t,x,m}_K)$ which solves (\ref{BSDEbis})
with generator~$f_K$.
With a comparison argument and since $f_K$ and $f$ coincide along the
sample paths of the solution $(Y^{t,x,m}_K,Z^{t,x,m}_K,L^{t,x,m}_K)$,
it can be shown that the bound of $Y^{t,x,m}_K$ does not depend on $K$,
if $K$ is large enough.
This is why $f_K$ can be replaced by~$f$ which satisfies (H2).
As a consequence, our proof is finished if we show that~(\ref
{EqMarkovPropertyQuad}) holds for the truncated generator $f_K$ which satisfies
(H2$'$).

The next step is to consider a BSDE which is shown in \cite{Morlais1}
to be in one to one correspondence with the BSDE (\ref{BSDEbis}) and
is obtained via an exponential coordinate
change. We only give a brief survey and refer to
\cite{Morlais1}, proof of Theorem 2.5, Step~2, for a complete
treatment. Setting
$U^{t,x,m}:=e^{\kappa Y^{t,x,m}}$ transforms (\ref{BSDEbis}) into
the following BSDE:
%
%e3.9 ###
%
\begin{eqnarray}
\label{BSDEInt}
U_s^{t,x,m}&=&e^{\kappa F(X_T^{t,x,m})}-\int_s^T V_r^{t,x,m} \,dM_r\nonumber\\
&&{} +
\int_s^T g(r,X_r^{t,x,m},M_r^{t,m},U_r^{t,x,m},V_r^{t,x,m}q_r^*) \,dC_r
\\
&&{} - \int_s^T dN_r^{t,x,m}, \qquad  s \in[t,T].\nonumber
\end{eqnarray}
We refer to a solution of this BSDE as
$(U^{t,x,m},V^{t,x,m},N^{t,x,m})$. Since $f_K$ satisfies~(H2$'$), the new
generator
\begin{eqnarray*}
&&g(s,x,m,u,v)\\
&&\qquad:=\biggl( \kappa\rho_{c^2}(u) f_K\biggl(s,x,m,\frac{\ln
(u \vee c^1)}{\kappa},\frac{v}{\kappa(u \vee c^1)} \biggr) - \frac
{1}{2(u \vee c^1)} |v|^2 \biggr),
\end{eqnarray*}
$(x,m,u,v) \in\mathcal V$, satisfies (H2$'$) (where $c^1$ and $c^2$ are
two explicit constants given in \cite{Morlais1}, pages 135--136,
depending only on $(a, \kappa, \|F\|_{\infty}, b)$ where we
recall\vadjust{\goodbreak}
that $a$ and $b$ are the constants appearing in the assumption (H2))
and the triple $(Y^{t,x,m},Z^{t,x,m},L^{t,x,m})$ with
%
%e3.10 ###
%
\begin{eqnarray}
\label{expchange}
Y^{t,x,m}&:=&\frac{\log(U^{t,x,m})}{\kappa}, \nonumber\\
Z^{t,x,m}&:=&\frac
{V^{t,x,m}}{\kappa U^{t,x,m}},\\
L^{t,x,m}&:=&\frac{1}{\kappa
U^{t,x,m}}\cdot N^{t,x,m}\nonumber
\end{eqnarray}
is well defined and is solution to (\ref{BSDEbis}) with generator
$f_K$ satisfying (H2$'$).

To derive\vspace*{1pt} the existence of a solution of (\ref{BSDEInt}), an
approximating sequence of \mbox{BSDEs} with Lipschitz generator $g^p$ and
terminal condition $e^{(\kappa F(X_T^{t,x,m}))}$ is introduced in
such a way that $g^p$ converges $d\nu$-almost everywhere to
$g$ as~$p$ tends to infinity. We do not specify the explicit
expression for $g^p$, since
we only need that the sequence is increasing in $y$, implying the same
property for the solution
component $(U^{p,t,x,m})_{p\in\mathbb{N}}$. For more details we
refer to \cite{Morlais1}, proof of Theorem 2.5, Step 3.

Let $p\geq1$. We consider the BSDE (\ref{BSDEInt}) with generator
$g^p$ and terminal condition $e^{(\kappa F(X_T^{t,x,m}))}$. Since
$g^p$ is Lipschitz continuous we know from~\cite{ElKarouiHuang},
Theorem~6.1, that a unique solution
$(U^{p,t,x,m},V^{p,t,x,m},N^{p,t,x,m})$ exists. Now we can apply
Proposition \ref{MarkovpropertyLip} which provides deterministic
functions~$a^p$ and~$b^p$ such that
\[
U_s^{p,t,x,m}=a^p(s,X_s^{t,x,m},M_s^{t,m})
\]
and
\[
V_s^{p,t,x,m}=b^p(s,X_s^{t,x,m},M_s^{t,m}), \qquad s\in[t,T].
\]
A subsequence, for convenience again denoted by
$(U^{p,t,x,m},V^{p,t,x,m},\break N^{p,t,x,m})_{p\in\mathbb{N}}$, converges
almost surely (with respect to $d\nu$) to the solution
$(U^{t,x,m},V^{t,x,m},N^{t,x,m})$ of (\ref{BSDEInt}).
Letting
\begin{eqnarray*}
a(s,y,m)&:=&\liminf_{p \to\infty} a^p(s,y,m), \\
b(s,y,m)&:=&\liminf
_{p \to\infty} b^p(s,y,m),
\end{eqnarray*}
$(s,y,m)\in[0,T]\times\mathbb{R}^{d\times1}\times
\mathbb{R}^{n\times1}$, we conclude that $U_s^{t,x,m} =
a(s,X_s^{t,x,m},M_s^{t,m})$ and $V_s^{t,x,m} =
b(s,X_s^{t,x,m},M_s^{t,m}), s\in[t,T]$. Since
$(U^{p,t,x,m})_{p\in\mathbb{N}}$ is increasing, we may set
\[
u:=\frac{\ln a}{\kappa},  \qquad v:=\frac{b}{\kappa a}.
\]
Hence, the result follows by (\ref{expchange}) and the one to one
correspondence.
\end{pf*}

%s4 ###
\section{Differentiability of FBSDEs}
\label{section:differentiability}

In this section we derive differentiability of the FBSDE of
(\ref{SDE}) and (\ref{BSDE}) with respect to the initial data $x$ and~$m$. The presence of the quantity $\langle L, L \rangle$ in the
equation,\vadjust{\goodbreak} where we recall that~$L$ is part of the solution of
(\ref{BSDE}), prevents us from extending directly the usual
techniques presented, for example, in
\cite{AnkirchnerImkellerDosReis1,AnkirchnerImkellerDosReis2,BriandConfortola}.
Under an additional assumption (MRP) defined in Section~\ref
{section:diff2}, we deduce the
differentiability of~(\ref{BSDE}) from that of the auxiliary
BSDE~(\ref{DiffFBSDE}).

%s4.1 ###
\subsection{Differentiability of an auxiliary FBSDE}
As mentioned above, we first prove the differentiability of an
auxiliary BSDE which will allow us
to deduce the result for (\ref{BSDE}) in Section \ref{section:diff2}.

For every $(x,m) \in\real^{(n+d) \times1}$, let us consider the
following forward--back\-ward system of equations:
%
%e4.1 ###
%
\begin{eqnarray}
\label{DiffFBSDE}
X_t^{x,m}&=&x+\int_0^t \sigma(r, X_r^{x,m},M^{m}_r) \,dM_r
+ \int_0^t
b(r,X_r^{x,m},M^{m}_r)\,dC_r,
\nonumber\\
Y_t^{x,m}&=&F(X_T^{x,m},M_T^{m})-\int_t^T Z_r^{x,m} \,dM_r\\
&&{} + \int_t^T
f(r,X_r^{x,m},M_r^m,Y_r^{x,m},Z_r^{x,m}q_r^*) \,dC_r,\nonumber
\end{eqnarray}
where $M$ is a continuous local martingale in $\real^{d\times1}$
satisfying the martingale representation property and $C, q, \sigma,
b, F, f$ are as described in Section~\ref{prel}. A~solution of this
system is given by the triple $(X^{x,m},Y^{x,m},Z^{x,m}) \in
\mathcal S^p \times\mathcal S^\infty\times L^2(d \langle M,M
\rangle\otimes d\PP)$ of stochastic processes. Note that the system
(\ref{DiffFBSDE}) has a~unique solution if the coefficients $\sigma$
and $b$ of the forward component satisfy~(H0) and the terminal
condition $F$ and the generator $f$ of the backward part satisfy
(H1)--(H3).

In this section\vspace*{1pt} we will give sufficient conditions for the system
(\ref{DiffFBSDE}) to be differentiable in $(x,m) \in\real^{(n+d)
\times1}$. Before turning to the backward SDE of the system, we
provide some material about the differentiability of the forward
component obtained in \cite{Protter}, Theorem V.7.39.
\begin{prop}
\label{differentiabilityofX} Assume that $\sigma$ and $b$ satisfy
\textup{(D1)}. Then for almost all $\omega\in\Omega$ there exists a
solution $X^{x,m}(\omega)$ of (\ref{DiffFBSDE}) which is
continuously differentiable in $x$ and $m$. In addition, the
derivatives $D^x_{ik}:=\frac{\partial}{\partial x_k} X^{(i)x,m}$,
$i,k=1,\ldots, n$, and $D^m_{ik}:=\frac{\partial}{\partial m_k}
X^{(i)x,m}$, $i=1,\ldots, n$, $k=1, \ldots, d$, satisfy\vspace*{1pt} the following
SDE for $t\in[0,T]$:
%
%e4.2 ###
%
\begin{eqnarray}
\label{partialx} D^x_{ik t} &=&\delta_{ik} + \sum_{\alpha=1}^d
\sum_{j=1}^n \int_0^t \partial_{1+j} \sigma_{i
\alpha}(s,X_s^{x,m},M_s^{m}) D^x_{jk s} \,dM^{(\alpha)}_{s}
\nonumber\\[-8pt]\\[-8pt]
&&{}+ \sum_{j=1}^n \int_0^t \partial_{1+j} b^{(i)}(s,X_s^{x,m},M_s^{m})
D^x_{jk s} \,d C_s,\nonumber\\[-1pt]
%
%e4.3 ###
%
\label{partialm} D^m_{ik t} &=& \sum_{\alpha=1}^d \sum_{j=1}^n
\int_0^t \partial_{1+j} \sigma_{i \alpha}(s,X_s^{x,m},M_s^{m})
D^m_{jk s} \,dM^{(\alpha)}_{s}\nonumber\\[-1pt]
&&{} + \sum_{j=1}^n \int_0^t \partial_{1+j}
b^{(i)}(s,X_s^{x,m},M_s^{m}) D^m_{jk s} \,d C_s
\nonumber\\[-8pt]\\[-8pt]
&&{}+\sum_{\alpha=1}^d \int_0^t \partial_{1+n+k} \sigma_{i
\alpha}(s,X_s^{x,m},M_s^{m}) \,dM^{(\alpha)}_{s}\nonumber\\[-1pt]
&&{} +\int_0^t
\partial_{1+n+k} b^{(i)}(s,X_s^{x,m},M_s^{m}) \,dC_s\nonumber
\end{eqnarray}
and $\frac{\partial}{\partial m_k}M^{(j)m}=\delta_{kj}$,
$k,j=1,\ldots,d$.
Furthermore, for all $p>1$ there exists a positive constant $\kappa$
such that the following estimate holds:
%
%e4.4 ###
%
\begin{equation}
\label{estX}
\E\Bigl[ {\sup_{t \in[0,T]}} |X_t^{x,m} - X_t^{x',m'}|^{2p} \Bigr]
\leq\kappa(|x - x'|^2 + |m-m'|^2)^{p}.
\end{equation}
\end{prop}
\begin{pf}
Let $(\tilde{X}_t^{x,m})_{t\in[0,T]}$ be the stochastic
process with values\break in $\real^{(1+n+d)\times1}$ defined as
\[
\tilde{X}_t^{x,m}=\pmatrix{
t\cr X_t^{x,m}\cr M_t^{m}}.
\]
This process is the solution of the SDE
\[
d\tilde{X}_t^{x,m}=\tilde{\sigma}(\tilde{X}_t^{x,m}) \,d\tilde{M}_t,
\qquad  \tilde{X}_0^{x,m}=(0,x,m)
\]
with
\[
\tilde{\sigma}(\tilde{X}_t^{x,m})=\pmatrix{
1&0&0\vspace*{2pt}\cr0&\sigma(\tilde{X}_t^{x,m})&b(\tilde
{X}_t^{x,m})\vspace*{2pt}\cr0&I_d&0},\qquad
\tilde{M}_t=\pmatrix{
t\cr M_t\cr C_t}.
\]
According\vspace*{1pt} to \cite{Protter}, Theorem V.39, the derivatives $D^x$,
$D^m$ and $\frac{\partial}{\partial m_k}M^{(j)m}$, $k,j=1,\ldots,d$,
exist and are continuous in $x$ and $m$. In addition, formula~\cite{Protter}, (D), page 312, leads to (\ref{partialx}) and
(\ref{partialm}). The estimate (\ref{estX}) follows immediately from
\cite{Protter}, $(***)$, page 309.
\end{pf}

We now focus on the backward part of system (\ref{DiffFBSDE}). Let
$\tilde
x:=(x,m) \in\real^{(n+d) \times1}$ and $e_i,  i=1,\ldots, n+d$, the
unit vectors in $\real^{(n+d) \times1}$. For all $\tilde x$, $h
\neq0$ and $i \in\{1, \ldots,n+d\}$, let $\xi^{\tilde
x,h,i}=\frac{1}{h}(F(X_T^{\tilde x + he_i},M_T^{\tilde x + h
e_i})-F(X_T^{\tilde x}, M_T^{\tilde x}))$. Here it is implicit that
$M^{\tilde x}$ only depends on the component $m$ in $\tilde x = (x,m)$. The
following lemma will be needed later in order to prove the
differentiability of the backward component. To simplify the
notation we suppress the superscript $i$.
\begin{lemma}
Suppose that \textup{(D1)} and \textup{(D2)} hold. Then for every $p>1$ there
exists a
constant $\kappa>0$, such that for all $\tilde x, \tilde x' \in
\real^{(n+d) \times1}$, $h,h' \neq0$
%
%e4.5 ###
%
\begin{equation}
\label{estF}
\E[ |\xi^{\tilde x,h} - \xi^{\tilde x', h'}|^{2p} ]
\leq
\kappa(|\tilde x- \tilde x'|^2 + |h - h'|^2)^p.
\end{equation}
\end{lemma}
\begin{pf}
Let $\tilde x, \tilde x' \in
\real^{(n+d) \times1}$ and $h,h' \neq0$. Given a real number $\theta
$ in $[0,1]$, we set
\[
G_i(\tilde x):= \partial_i F\bigl(X_T^{\tilde x} +\theta(X_T^{\tilde x +
he_i}-X_T^{\tilde x}),M_T^{\tilde x}+\theta(M_T^{\tilde x + h
e_i}-M_T^{\tilde x})\bigr), \qquad  i=1,2.
\]
For notational convenience, we also define
\begin{eqnarray*}
H&:=& \frac{X_T^{\tilde x + he_i}-X_T^{\tilde x}}{h}-\frac{X_T^{\tilde
x' + h'e_i}-X_T^{\tilde x'}}{h'}, \\
I&:=& \frac{M_T^{\tilde x +
he_i}-M_T^{\tilde x}}{h}-\frac{M_T^{\tilde x' + h'e_i}-M_T^{\tilde
x'}}{h'}.
\end{eqnarray*}
We have
\begin{eqnarray*}
&&\E[ |\xi^{\tilde x,h} - \xi^{\tilde x', h'}|^{2p} ]
\\
&&\qquad=\E\biggl[ \biggl|\frac{1}{h}\bigl(F(X_T^{\tilde x +
he_i},M_T^{\tilde x + h e_i})-F(X_T^{\tilde x}, M_T^{\tilde x})\bigr)
\\
&&\qquad\quad\hphantom{\E\biggl[ \biggl|}
{}
-\frac{1}{h'}\bigl(F(X_T^{\tilde x' + h'e_i},M_T^{\tilde x' + h'
e_i})-F(X_T^{\tilde x'}, M_T^{\tilde x'})\bigr)\biggr|^{2p} \biggr]
\\
&&\qquad=\E\biggl[ \biggl| \int_0^1 \biggl( G_1(\tilde x) \frac{X_T^{\tilde
x + he_i}-X_T^{\tilde x}}{h} + G_2(\tilde x) \frac{M_T^{\tilde x +
he_i}-M_T^{\tilde x}}{h} \\
&&\qquad\quad\hphantom{\E\biggl[ \biggl| \int_0^1 \biggl(}
{}-G_1(\tilde x') \frac
{X_T^{\tilde x' + h'e_i}-X_T^{\tilde x'}}{h'}-G_2(\tilde x') \frac
{M_T^{\tilde x' + h'e_i}-M_T^{\tilde x'}}{h'} \biggr) \,d\theta
\biggr|^{2p} \biggr]
\\
&&\qquad=\E\biggl[ \biggl|\int_0^1 G_1(\tilde x) H - \bigl(G_1(\tilde
x')-G_1(\tilde x)\bigr) \frac{X_T^{\tilde x' + h'e_i}-X_T^{\tilde
x'}}{h'} \\
&&\qquad\quad\hphantom{\E\biggl[ \biggl|}
{}+ G_2(\tilde x) I -\bigl(G_2(\tilde
x')-G_2(\tilde x)\bigr) \frac{M_T^{\tilde x' + h'e_i}-M_T^{\tilde
x'}}{h'} \,d\theta\biggr|^{2p} \biggr]
\\
&&\qquad\leq
c \E\biggl[ |H|^{2p}
+\biggl\vert\frac{X_T^{\tilde x' + h'e_i}-X_T^{\tilde x'}}{h'}
\biggr\vert^{2p} \biggl(\int_0^1 |G_1(\tilde x')-G_1(\tilde x)| \,d\theta
\biggr)^{2p} \biggr]\\
&&\qquad\quad{}+ c \E\biggl[|I|^{2p} + \biggl\vert\frac
{M_T^{\tilde x' + h'e_i}-M_T^{\tilde x'}}{h'} \biggr\vert^{2p}
\biggl(\int_0^1 |G_2(\tilde x')-G_2(\tilde x)| \,d\theta\biggr)^{2p}
\biggr]\\
&&\qquad=: T_1+T_2,
\end{eqnarray*}
where we have used the fact that $F$ is globally Lipschitz in the last
inequality.
Similarly, the Lipschitz property of $\nabla F$ entails for $i=1,2$
\begin{eqnarray*}
&&\int_0^1 G_i(\tilde x')-G_i(\tilde x) \,d\theta
\\
&&\qquad\leq C (|X_T^{\tilde x} - X_T^{\tilde x'}| + |X_T^{\tilde x +
he_i}-X_T^{\tilde x' + h'e_i}|\\
&&\qquad\quad\hphantom{C (}
{} + |M_T^{\tilde x} - M_T^{\tilde x'}| +
|M_T^{\tilde x + he_i}-M_T^{\tilde x' + h'e_i}|)=:J.
\end{eqnarray*}
Hence, using the H\"{o}lder inequality with $\gamma, q>1$ s.t. $\frac
{1}{\gamma} + \frac{1}{q}=1$ we get
\begin{eqnarray*}
T_1&\leq& c \E[|H|^{2p}]+ c \E\biggl[\biggl\vert\frac{X_T^{\tilde x'
+ h'e_i}-X_T^{\tilde x'}}{h'}\biggr\vert^{2p} J^{2p} \biggr]\\
&\leq& c \E[|H|^{2p}] + c \E\biggl[\biggl\vert\frac{X_T^{\tilde x' +
h'e_i}-X_T^{\tilde x'}}{h'} \biggr\vert^{2p\gamma}
\biggr]^{1/\gamma} \E[J^{2pq}]^{1/q}.
\end{eqnarray*}
Recall that $\E[|X_T^{\tilde x}|^r]<\infty$ for all $r\geq1$ and thus,
from inequality (\ref{estX}), we have
\[
\E\biggl[\biggl\vert\frac{X_T^{\tilde x' + h'e_i}-X_T^{\tilde
x'}}{h'} \biggr\vert^{2p\gamma} \biggr]^{1/\gamma} = \frac
{1}{(h')^{2p}} \E[\vert X_T^{\tilde x' +
h'e_i}-X_T^{\tilde x'}\vert^{2p\gamma} ]^{1/\gamma}
\leq c ,
\]
where $c$ is a constant which does not depend on $\tilde x, \tilde x',
h$ or $h'$.
Combining the previous estimates we finally obtain
\[
T_1 \leq c \E[|H|^{2p}] + c \E[J^{2pq}]^{1/q} \leq c (|\tilde x -
\tilde x'|^2 + |h-h'|^2)^p.
\]
The same method gives that
\[
T_2 \leq c \E[|H|^{2p}] + c \E[J^{2pq}]^{1/q} \leq c (|\tilde x -
\tilde x'|^2 + |h-h'|^2)^p
\]
and the proof is complete.
\end{pf}

The next lemma shows that we can choose the family $(Y^{\tilde x})$
to be continuous in $\tilde x \in\real^{(n+d) \times1}$.
\begin{lemma}
\label{continuity} Let \textup{(H1)}--\textup{(H3)} and \textup{(D1)}--\textup{(D3)} be satisfied.
Then for all $p>1$ there exists a constant $c>0$, such that for all
$\tilde x, \tilde x' \in\real^{(n+d) \times1}$
%
%e4.6 ###
%
\begin{eqnarray}
\label{estcont}
&&\E\Bigl[ {\sup_{t \in[0,T]}} |Y_t^{\tilde x} -
Y_t^{\tilde x'}|^{2p} \Bigr] + \E\biggl[ \biggl(\int_0^T
|q_t(Z_t^{\tilde x}-Z_t^{\tilde x'})^*|^{2} \,dC_t \biggr)^p \biggr]
\nonumber\\[-8pt]\\[-8pt]
&&\qquad\leq c |\tilde x-\tilde x'|^{2p}.\nonumber
\end{eqnarray}
Furthermore, for almost all $\omega\in\Omega$ there exists a
solution $Y^{\tilde x}(\omega)$ of (\ref{DiffFBSDE}) which is
continuous in $\tilde x \in\real^{(n+d) \times1}$.
\end{lemma}
\begin{pf}
Let $\de Y:=Y^{\tilde x}-Y^{\tilde x'}$, $\de Z:=Z^{\tilde
x}-Z^{\tilde x'}$, $\de M:=M^{m}-M^{m'}$ and $\de X:=X^{\tilde
x}-X^{\tilde x'}$. We also set for $s\in[0,T]$
\begin{eqnarray*}
A^Z_r&:=&\int_0^1 \partial_{5}f\bigl(r,X_r^{\tilde x},M_r^m,Y_r^{\tilde
x},Z_r^{\tilde x'} q_r^* +\zeta(Z_r^{\tilde x}-Z_r^{\tilde x'}) q_r^*\bigr)
\,d\zeta,\\
A^Y_r&:=&\int_0^1 \partial_{4}f\bigl(r,X_r^{\tilde x},M_r^m,Y_r^{\tilde
x'}+\zeta(Y_r^{\tilde x}-Y_r^{\tilde x'}),Z_r^{\tilde x'}q_r^*\bigr)
\,d\zeta,
\\
A^M_r&:=&\int_0^1 \partial_{3}f\bigl(r,X_r^{\tilde x},M_r^{m'}+\zeta
(M_r^{m}-M_r^{m'}),Y_r^{\tilde x'},Z_r^{\tilde x'} q_r^*\bigr) \,d\zeta,\\
A^X_r&:=&\int_0^1 \partial_{2}f\bigl(r,X_r^{\tilde x'}+\zeta(X_r^{\tilde
x}-X_r^{\tilde x'}),M_r^{m'},Y_r^{\tilde x'},Z_r^{\tilde x'}q_r^*\bigr)
\,d\zeta.
\end{eqnarray*}
Considering the difference $\de Y$ of the backward component in
(\ref{DiffFBSDE}) we see that for $t\in[0,T]$
\begin{eqnarray*}
\de Y_t
&=& F(X_T^{\tilde x},M_T^m)\\
&&{}-F(X_T^{\tilde x'},M_T^{m'})- \int_t^T \de
Z_r \,dM_r
\\
&&{}+ \int_t^T [f(r,X_r^{\tilde x},M_r^m,Y_r^{\tilde x},Z_r^{\tilde x}
q_r^*)-f(r,X_r^{\tilde x'},M_r^{m'},Y_r^{\tilde x'},Z_r^{\tilde
x'}q_r^*)] \,dC_r
\\
&=& F(X_T^{\tilde x},M_T^m)-F(X_T^{\tilde x'},M_T^{m'})- \int_t^T \de
Z_r \,dM_r
\\
&&{} +
\int_t^T \underbrace{( \de Z_r q_r^* A^Z_r+ \de Y_r A^Y_r + \de
M_r^* A^M_r + \de X_r^* A^X_r )}_{=:g(r,\de Y_r,\de Z_r q_r^*)} \,dC_r
\end{eqnarray*}
holds. Note that $(\de Y, \de Z )$ can be seen as a BSDE whose
generator $g$ satisfies (H4) and whose terminal condition
$F(X_T^{\tilde x},M_T^m)-F(X_T^{\tilde x'},M_T^{m'})$ is bounded [see
(H1)]. More
precisely, we derive with (D3) and \cite{Morlais1}, Lemma 3.1, the
existence of a constant $c$ such that for all $y, y_1, y_2 \in
\real$ and $z, z_1,z_2 \in\real^{1 \times d}$ $\nu$-a.e.
\begin{eqnarray*}
&&|g(r, y,z_1)-g(r,y,z_2)|\\
&&\qquad\leq |A^Z_r| |z_1-z_2|\\
&&\qquad\leq c \bigl( |q_r \theta_r| + |Z_r^{\tilde x'}q_r^*| + |(Z_r^{\tilde
x}-Z_r^{\tilde x'})q_r^*| \bigr) |z_1-z_2|
\end{eqnarray*}
and
\begin{eqnarray*}
&&|g(r, y_1,z)-g(r,y_2,z)| \\
&&\qquad\leq |A^M_r| |y_1-y_2|
\leq c ( |q_r \theta_r| + |Z_r^{\tilde x'}q_r^*|) |y_1-y_2|.
\end{eqnarray*}
Hence, we can apply the a priori estimates of Lemma \ref
{appendix:apriori} and hence, we know that for every $p>1$ there exist
constants $q>1$ and $c>0$ such that
%
%e4.7 ###
%
\begin{eqnarray}
\label{cont1}
&& \E\Bigl[ \sup_{t \in[0,T]} |\delta Y_t|^{2p} \Bigr] + \E
\biggl[ \biggl(\int_0^T |q_t \delta Z_t^*|^2 \,dC_t \biggr)^{p} \biggr] \nonumber\\
&&\qquad \leq c \E\biggl[ |F(X_T^{\tilde x},M_T^m)-F(X_T^{\tilde
x'},M_T^{m'})|^{2pq} \\
&&\qquad\quad\hphantom{c \E\biggl[}
{}+
\biggl( \int_0^T |\delta M_r^* A^M_r + \delta X_r^* A^X_r| \,dC_r
\biggr)^{2pq} \biggr]^{{1/q}}. \nonumber
\end{eqnarray}
By condition (D3) and H\"{o}lder's inequality we get
\begin{eqnarray*}
\hspace*{-3pt}&&\E\biggl[ \biggl( \int_0^T |\delta M_r^* A^M_r + \delta X_r^* A^X_r|
\,dC_r \biggr)^{2pq} \biggr]
\\
\hspace*{-3pt}&&\qquad\leq c \E\biggl[ \biggl(\int_0^T |\delta M_r|^2 \,dC_r\biggr)^{2pq}
\biggr]^{1/2} \E\biggl[ \biggl(\int_0^T (|q_r \theta_r| +
|Z_r^{\tilde x'} q_r^*|)^2 \,dC_r\biggr)^{2pq} \biggr]^{1/2}
\\
\hspace*{-3pt}&&\qquad\quad{}  + c \E\biggl[ \biggl(\int_0^T |\delta X_r|^2 \,dC_r
\biggr)^{2pq} \biggr]^{1/2} \E\biggl[ \biggl(\int_0^T (|q_r \theta
_r| + |Z_r^{\tilde x'} q_r^*|)^2 \,dC_r\biggr)^{2pq} \biggr]^{1/2}.
\end{eqnarray*}
Note that $\E[ (\int_0^T |q_r \theta_r|^2
\,dC_r)^{2pq} ]$ is bounded by (D3). Furthermore,
\[
\E
\biggl[ \biggl(\int_0^T |Z_r^{\tilde x'} q_r^*|^2 \,dC_r\biggr)^{2pq}
\biggr]
\]
is bounded, as is seen by applying Lemma
\ref{appendix:apriori}. Hence,
\begin{eqnarray*}
&&\E\biggl[ \biggl( \int_0^T |\delta M_r^* A^M_r + \delta X_r^* A^X_r|
\,dC_r \biggr)^{2pq} \biggr]\\
&&\qquad\leq c |m-m'|^{2pq} + C \E\Bigl[
\Bigl( \sup_{t \in[0,T]} |\delta X_t|^2 C_T\Bigr)^{2pq}
\Bigr]^{1/2}
\\
&&\qquad\leq c (|m-m'|^{2pq} + |\tilde x - \tilde x'|^{2pq} ),
\end{eqnarray*}
where the last inequality is due to (\ref{estX}). Combining (\ref
{cont1}), condition (D2) and the last inequality we obtain
\[
\E\Bigl[ \sup_{t \in[0,T]} |\delta Y_t|^{2p} \Bigr] + \E\biggl[
\biggl(\int_0^T |q_s \delta Z_s^*|^2\biggr)^{p} \biggr]
\leq c |\tilde x - \tilde x'|^{2p}.
\]
Now Kolmogorov's lemma (see
\cite{Protter}, Theorem 73, Chapter IV) implies that there exists a
version of $(Y^{\tilde
x})$ which is continuous in $\tilde x$ for almost all $\omega\in
\Omega$.
\end{pf}

For all $h \neq0, \tilde x \in\real^{(n+d)\times1}, t\in[0,T]$
let $U_t^{\tilde x,h}=\frac{1}{h} (Y_t^{\tilde x + he_i}-Y_t^{\tilde
x})$, $V_t^{\tilde x,h}=\frac{1}{h} (Z_t^{\tilde x +
he_i}-Z_t^{\tilde x})$, $\Delta_t^{\tilde x,h}=\frac{1}{h}
(X_t^{\tilde x + he_i}-X_t^{\tilde x})$, $\varpi_t^{\tilde
x,h}=\frac{1}{h} (M_t^{\tilde x + he_i}-M_t^{\tilde x})$ [where it
is implicit that $M^{\tilde{x}}$ depends only on the component $m$
of $\tilde{x} = (x,m)$] and $\xi^{\tilde x,h}=\frac{1}{h}
(F(X_T^{\tilde x + he_i},M_T^{\tilde x + h e_i})-F(X_T^{\tilde
x},M_T^{\tilde x}))$.
We define $\delta U$ by $\delta U=U^{\tilde x,h} - U^{\tilde x',h'}$
and the processes
$\delta V$, $\de\Delta$, $\de\varpi$ and $\delta\xi$ in an
analogous way.
We give estimates on the differences of difference quotients of the
family $(Y^{\tilde x})$.
\begin{lemma}
\label{lemma:differentiability} Let \textup{(H1)}--\textup{(H3)} and
\textup{(D1)}--\textup{(D4)} be
satisfied. Then for each $p>1$ there exists a constant $c>0$ such
that for any $\tilde x, \tilde x' \in\real^{(n+d) \times1}$ and~%
$h,\allowbreak h' \neq0$
%
%e4.8 ###
%
\begin{equation}
\label{estdiffY}
\E\Bigl[ {\sup_{t \in[0,T]}} |U^{\tilde x,h}_t - U^{\tilde
x',h'}_t|^{2p} \Bigr] \leq c (|\tilde x - \tilde x'|^2 + |h-h'|^2)^p.
\end{equation}
\end{lemma}
\begin{pf}
This proof is similar to that of Lemma \ref{continuity}. By definition
of $U^{\tilde x,h}$ and of $U^{\tilde x',h'}$ we have
%
%e4.9 ###
%
\begin{eqnarray}
\label{dif1}
U_t^{\tilde x,h} &=& \xi^{\tilde x,h} - \int_t^T V_r^{\tilde x,h}
\,dM_r\nonumber\\
&&{}+ \int_t^T \frac{1}{h} [f(r,X_r^{\tilde x + h e_i},M_r^{\tilde
x + h e_i},Y_r^{\tilde x+h ei},Z_r^{\tilde x+h e_i} q_r^\ast)\\
&&\hspace*{102pt}{}-
f(r,X_r^{\tilde x},M_r^{\tilde{x}},Y_r^{\tilde x},Z_r^{\tilde x}
q_r^\ast)] \,dC_r.\nonumber
\end{eqnarray}
As in the proof of Lemma \ref{continuity}, we decompose the integrand
in the last term of the right-hand side of the equality above by writing
\begin{eqnarray*}
&&\frac{1}{h}
\bigl(f(r,X_r^{\tilde x + h e_i},M_r^{\tilde x + h e_i},Y_r^{\tilde x+h
ei},Z_r^{\tilde x+h e_i} q_r^\ast)-f(r,X_r^{\tilde x},M_r^{\tilde
{x}},Y_r^{\tilde x},Z_r^{\tilde x} q_r^\ast)\bigr)\\
&&\qquad= V_t^{\tilde x,h} q_r^\ast(A^Z)_r^{\tilde x,h} + U_t^{\tilde{x},
h} (A^Y)_r^{\tilde{x}, h} + {\varpi_r^{\tilde{x}, h}}^\ast
(A^M)_r^{\tilde{x}, h} + \Delta^{\tilde{x}, h} (A^X)_r^{\tilde{x},
h},
\end{eqnarray*}
where $A^Z, A^Y, A^M, A^X$ are defined as in the proof of
Lemma \ref{continuity}, for instance,
\[
(A^Z)_r^{\tilde x,h}:=\int_0^1 \partial_{5}f\bigl(r,X_r^{\tilde{x} + h
e_i},M_r^{\tilde{x} +h e_i},Y_r^{\tilde{x} + h e_i},Z_r^{\tilde{x}}
q_r^* +
\theta(Z_r^{\tilde{x}+ h e_i}-Z_r^{\tilde{x}}) q_r^*\bigr)
\,d\theta.
\]
Taking the difference of two equations of the form (\ref{dif1}) we
obtain that~$(\de U,\allowbreak \de V)$ satisfies the BSDE
%
%e4.10 ###
%
\begin{eqnarray}
\label{dif2}\quad
\de U_t &=&
\de\xi-\int_0^T \de V_r \,dM_t\nonumber\\
&&{}+ \int_t^T \de V_r q_r^\ast(A^Z)_r^{\tilde x,h}+ \de U_r
(A^Y)_r^{\tilde x,h}\nonumber\\
&&{} + \bigl[q_r^\ast\bigl((A^Z)_r^{\tilde x,h}-(A^Z)_r^{\tilde
x',h'}\bigr)V_r^{\tilde x',h'}
+ U_r^{\tilde x',h'}
\bigl((A^Y)_r^{\tilde x,h}- (A^Y)_r^{\tilde
x',h'}\bigr)\\
&&\hphantom{{}+ \bigl[}
{}+{\varpi_r^{\tilde{x}, h}}^\ast(A^M)_r^{\tilde{x},
h}-{\varpi_r^{\tilde{x'}, h'}}^\ast(A^M)_r^{\tilde{x'},
h'}
\nonumber\\
&&\hspace*{126.4pt}{} +{\Delta^{\tilde{x}, h}}^\ast(A^X)_r^{\tilde{x},
h}-{\Delta^{\tilde{x}', h'}}^\ast(A^X)_r^{\tilde{x}', h'}\bigr]
\,dC_r.\nonumber
\end{eqnarray}
The generator of this BSDE satisfies condition (H4) due to
assumption (D3) (details are similar to those of the proof of Lemma
\ref{continuity} and are left to the reader). By Lemma
\ref{appendix:apriori}, for every $p>1$ there exist constants $q>1$
and $c>0$ such that
\begin{eqnarray*}
&&\E\biggl[ {\sup_{t\in[0,T]}} \vert\de U_t \vert^{2 p} + \biggl(
\int_0^T \vert q_s \de V_s^\ast\vert^2 \,dC_s \biggr)^p \biggr]
\\
&&\qquad\leq c \E\biggl[ \vert\de\xi\vert^{2 p q} + \biggl( \int_0^T \bigl\vert
q_r^\ast\bigl((A^Z)_r^{\tilde x,h}-(A^Z)_r^{\tilde x',h'}\bigr)\bigr\vert\vert
V_r^{\tilde x',h'} \vert\\
&&\qquad\quad\hspace*{67.8pt}{} + \vert U_r^{\tilde x',h'} \vert\bigl\vert
\bigl((A^Y)_r^{x,h} - (A^Y)_r^{\tilde x',h'}\bigr)\bigr\vert\\
&&\qquad\quad\hspace*{67.8pt}{} +\vert{\varpi_r^{\tilde{x}, h}}^\ast(A^M)_r^{\tilde{x},
h}- {\varpi_r^{\tilde{x}', h'}}^\ast(A^M)_r^{\tilde{x}', h'}\vert\\
&&\qquad\quad\hspace*{67.8pt}{} +\vert{\Delta^{\tilde{x}, h}}^\ast(A^X)_r^{\tilde{x}, h}-
{\Delta^{\tilde{x}', h'}}^\ast(A^X)_r^{\tilde{x}', h'}\vert \,dC_r
\biggr)^{2 p q} \biggr]^{1/q}.
\end{eqnarray*}
We estimate separately each part of the right-hand side of the
inequality presented. First, by Cauchy--Schwarz's inequality we have
\begin{eqnarray*}
&& \E\biggl[\biggl(\int_0^T \bigl\vert q_r^\ast\bigl((A^Z)_r^{\tilde
x,h}-(A^Z)_r^{\tilde x',h'}\bigr)\bigr\vert\vert V_r^{\tilde x',h'} \vert
\,dC_r\biggr)^{2 p q}\biggr]
\\
&&\qquad \leq\E\biggl[\biggl(\int_0^T \bigl\vert q_r^\ast\bigl((A^Z)_r^{\tilde
x,h}-(A^Z)_r^{\tilde x',h'}\bigr)\bigr\vert^2 \,dC_r\biggr)^{2 p q}\biggr]^{1/2}
\\
&&\qquad\quad{}\times\E\biggl[\biggl(\int_0^T \vert V_r^{\tilde x',h'} \vert^2 \,dC_r\biggr)^{2 p
q}\biggr]^{1/2}
\\
&&\qquad \leq C \E\biggl[\biggl(\int_0^T \bigl\vert q_r^\ast\bigl((A^Z)_r^{\tilde
x,h}-(A^Z)_r^{\tilde x',h'}\bigr)\bigr\vert^2 \,dC_r\biggr)^{2 p q}\biggr]^{1/2}
\end{eqnarray*}
since $\E[(\int_0^T \vert V_r^{\tilde x',h'} \vert^2 \,dC_r)^{2
p q}]$ is bounded by Lemma \ref{appendix:apriori}. Then hypothesis~%
(D4) and a combination of Lemma \ref{continuity} and (\ref{estX})
lead to the following estimate:
\begin{eqnarray*}
&& \E\biggl[\biggl(\int_0^T \bigl\vert q_r^\ast\bigl((A^Z)_r^{\tilde
x,h}-(A^Z)_r^{\tilde x',h'}\bigr)\bigr\vert\vert V_r^{\tilde x',h'} \vert \,dC_r
\biggr)^{2 p q}\biggr]\\
&&\qquad\leq c \E\biggl[\biggl(\int_0^T \vert q_r^\ast(Z_r^{\tilde
{x}}-Z_r^{\tilde{x}'}) \vert^2\\
&&\qquad\quad\hphantom{c \E\biggl[\biggl(}
{} +
\vert X_r^{\tilde{x}+h e_i}-X_r^{\tilde{x}'+h' e_i} \vert^2+ \vert
\underbrace{M_r^{\tilde{x}+h e_i}-M_r^{\tilde{x}'+h' e_i}}_{=\tilde
x +he_i-\tilde x' -h'e_i} \vert^2
\\
&&\qquad\quad\hphantom{c \E\biggl[\biggl(}
{}  +\vert Y_r^{\tilde{x}+h e_i}-Y_r^{\tilde{x}'+h' e_i} \vert
^2+\vert q_r^\ast(Z_r^{\tilde{x}+h e_i}-Z_r^{\tilde{x}'+h' e_i})
\vert^2 \,dC_r\biggr)^{2 p q}\biggr]^{1/2} \\
&&\qquad \leq c ( \vert\tilde{x} - \tilde{x}' \vert^2 +\vert h - h'
\vert^2 )^{pq}.
\end{eqnarray*}
Similarly, we derive
\begin{eqnarray*}
&&\E\biggl[\biggl(\int_0^T \vert U_r^{\tilde x',h'} \vert\vert(A^Y)_r^{x,h}
- (A^Y)_r^{\tilde x',h'}\vert \,dC_r\biggr)^{2 p q}\biggr]\\
&&\qquad\leq c ( \vert
\tilde{x} - \tilde{x}' \vert^2 +\vert h - h' \vert^2 )^{pq}.
\end{eqnarray*}
We next estimate
\begin{eqnarray*}
&& \E\biggl[\biggl(\int_0^T \vert{\varpi_r^{\tilde{x}, h}}^\ast
(A^M)_r^{\tilde{x}, h}- {\varpi_r^{\tilde{x}', h'}}^\ast
(A^M)_r^{\tilde{x}', h'}\vert \,dC_r\biggr)^{2 p q}\biggr]
\\
&&\qquad\leq c \E\biggl[\biggl(\int_0^T \vert\underbrace{
\varpi_r^{\tilde{x},h}-\varpi_r^{\tilde{x}',h'}}_{=0} \vert\vert
(A^M)_r^{\tilde{x} ,h} \vert \,dC_r\biggr)^{2 p q}\biggr]
\\
&&\qquad\quad{}  + c \E\biggl[\biggl(\int_0^T
\vert\underbrace{\varpi_r^{\tilde{x}',h'}}_{=e_i} \vert\vert
(A^M)_r^{\tilde{x} ,h}-(A^M)_r^{\tilde{x}' ,h'} \vert \,dC_r\biggr)^{2 p
q}\biggr]
\\
&&\qquad\leq c \E\biggl[\biggl(\int_0^T \vert(A^M)_r^{\tilde{x}
,h}-(A^M)_r^{\tilde{x}' ,h'} \vert \,dC_r\biggr)^{2 p q}\biggr]
\\
&&\qquad \leq c \E\biggl[\biggl(\int_0^T ( |q_r \theta_r| + \vert Z_r^{\tilde x}
q_s^\ast\vert+\vert Z_r^{\tilde x'} q_s^\ast\vert)^2 \,dC_r\biggr)^{2 p
q}\biggr]^{1/2} \\
&&\qquad\quad{}\times\E\biggl[\biggl(\int_0^T \bigl(\vert X_r^{\tilde x+h e_i} -
X_r^{\tilde x'+h' e_i} \vert\\
&&\qquad\quad\hspace*{52.2pt}{} + \vert Y_r^{\tilde x} - Y_r^{\tilde x'}\vert+ \vert
(Z_r^{\tilde x} - Z_r^{\tilde x'})q_r^\ast\vert\\
&&\qquad\quad\hspace*{52.2pt}{}+\vert M_r^{\tilde x}
- M_r^{\tilde x'} \vert+\vert M_r^{\tilde x+h e_i} - M_r^{\tilde x'+h'
e_i} \vert\bigr)^2 \,dC_r\biggr)^{2 p q}\biggr]^{1/2},
\end{eqnarray*}
where the last inequality is due to hypothesis (D4) and H\"{o}lder's
inequality. An application of the a priori estimates from Lemma
\ref{appendix:apriori} implies that $\E[(\int_0^T (|q_r
\theta_r| + \vert Z_r^{\tilde x} q_r^\ast\vert+\vert Z_r^{\tilde x'}
q_r^\ast\vert)^2 \,dC_r)^{2 p q}]$ is bounded.
Then, using (\ref{estX})\break and~(\ref{estcont}), we obtain
\begin{eqnarray*}
&& \E\biggl[\biggl(\int_0^T \vert{\varpi_r^{\tilde{x}, h}}^\ast
(A^M)_r^{\tilde{x}, h}- {\varpi_r^{\tilde{x}', h'}}^\ast
(A^M)_r^{\tilde{x}', h'}\vert \,dC_r\biggr)^{2 p q}\biggr]
\\
&&\qquad\leq c \E\biggl[\biggl(\int_0^T \vert X_r^{\tilde x+h e_i} - X_r^{\tilde
x'+h' e_i} \vert^2\\
&&\qquad\quad\hphantom{c \E\biggl[\biggl(}
{} + \vert Y_r^{\tilde x} - Y_r^{\tilde x'}\vert^2 +
\vert(Z_r^{\tilde x} - Z_r^{\tilde x'})q_r^\ast\vert^2
\\
&&\qquad\quad\hphantom{c \E\biggl[\biggl(}
{}  +\vert M_r^{\tilde x} - M_r^{\tilde x'} \vert^2 +
\vert M_r^{\tilde x+h e_i} - M_r^{\tilde x'+h' e_i} \vert^2 \,dC_r
\biggr)^{2 p q}\biggr]^{1/2}\\
\hspace*{-2pt}&&\qquad\leq c (\vert\tilde{x} -\tilde{x}' \vert^2 + \vert h -h'
\vert^2)^{pq}.
\end{eqnarray*}
We now consider the last term whose treatment is similar to that of the
term just discussed.
Therefore, we give the main computations without providing detailed arguments.
We have
\begin{eqnarray*}
&& \E\biggl[\biggl(\int_0^T
\vert{\Delta_r^{\tilde{x}, h}}^\ast(A^X)_r^{\tilde{x}, h}- {\Delta
_r^{\tilde{x}', h'}}^\ast(A^X)_r^{\tilde{x}', h'}\vert \,dC_r\biggr)^{2 p
q}\biggr]\\
&&\qquad\leq c \E\biggl[\biggl(\int_0^T \vert\Delta_r^{\tilde{x},h}-\Delta
_r^{\tilde{x}',h'} \vert^2 \,dC_r\biggr)^{2 p q}\biggr]^{1/2}\\
&&\qquad\quad{}\times \E\biggl[\biggl(\int
_0^T \vert(A^X)_r^{\tilde{x} ,h}
\vert^2 \,dC_r\biggr)^{2 p q}\biggr]^{1/2}\\
&&\qquad\quad{}  + c \E\biggl[\biggl(\int_0^T \vert\Delta_r^{\tilde{x}',h'} \vert
\vert(A^X)_r^{\tilde{x} ,h}-(A^X)_r^{\tilde{x}' ,h'} \vert
\,dC_r\biggr)^{2 p q}\biggr].
\end{eqnarray*}
Using (D3) and Lemma \ref{appendix:apriori}, we deduce that
$\E[(\int_0^T \vert(A^X)_r^{\tilde{x} ,h} \vert^2 \,dC_r
)^{2 p
q}]$ is bounded. Using hypothesis (D4) and (\ref{estX}), again we obtain
\begin{eqnarray*}
&& \E\biggl[\biggl(\int_0^T \vert{\Delta_r^{\tilde{x}, h}}^\ast
(A^X)_r^{\tilde{x}, h}- {\Delta_r^{\tilde{x}', h'}}^\ast
(A^X)_r^{\tilde{x}', h'}\vert \,dC_r\biggr)^{2 p q}\biggr]
\\
&&\qquad \leq c \E\biggl[\biggl(\int_0^T \vert\Delta_r^{\tilde{x},h}-\Delta
_r^{\tilde{x}',h'} \vert^2 \,dC_r \biggr)^{2pq} \biggr]^{1/2}
\\
&&\qquad\quad{}+ c \E\biggl[\biggl(\int_0^T ( \vert q_r \theta_r \vert+ \vert
Z_r^{\tilde x'} q_r^\ast\vert+ \vert Z_r^{\tilde x} q_r^\ast\vert
)\\
&&\qquad\quad\hphantom{+ c \E\biggl[\biggl(}
{}\times \bigl( \vert X_r^{\tilde x} -X_r^{\tilde x'} \vert+ \vert X_r^{\tilde
x + h e_i} -X_r^{\tilde x'+h' e_i} \vert+ \vert M_r^{\tilde x}
-M_r^{\tilde x'} \vert
\\
&&\qquad\quad\hspace*{126.4pt}{} + \vert Y_r^{\tilde x} -Y_r^{\tilde x'} \vert+
\vert(Z_r^{\tilde x} -Z_r^{\tilde x'})q_r^\ast\vert\bigr) \,dC_r\biggr)^{2
p q}\biggr]
\\
&&\qquad \leq
c (\vert\tilde{x} - \tilde{x}' \vert^2 + \vert h - h'
\vert^2)^{pq}.
\end{eqnarray*}
We derive
\[
\E[ |\de\xi|^{2pq}]
\leq
c (\vert\tilde{x} - \tilde{x}' \vert^2 + \vert h - h' \vert^2)^{pq}
\]
from (\ref{estF}). This completes the proof of
(\ref{estdiffY}).
\end{pf}
\begin{prop}
\label{PropDiff} Let \textup{(H1)}--\textup{(H3)} and \textup{(D1)}--\textup{(D4)} be
satisfied. Then
there exists a solution $(X^{\tilde x},Y^{\tilde x},Z^{\tilde x})$
of (\ref{DiffFBSDE}), such that $X^{\tilde x}(\omega)$ and
$Y^{\tilde x}(\omega)$ are continuously differentiable in $\tilde x
\in\real^{(n+d) \times1}$ for almost all $\omega\in\Omega$.
Furthermore, there exist processes $\frac{\partial}{\partial x}
Z^{x,m}, \frac{\partial}{\partial m} Z^{x,m} \in L^2(d\langle
M,M\rangle\otimes d\PP)$ such that the derivatives
$(U_k^x,V_{ik}^x):=(\frac{\partial}{\partial x_k} Y^{x,m},
\frac{\partial}{\partial x_k} Z^{(i),x,m})$, $i=1,\ldots,d$, $k=1,
\ldots, n$, and $(U_k^m,V_{ik}^m):=(\frac{\partial}{\partial m_k}
Y^{x,m}, \frac{\partial}{\partial m_k} Z^{(i),x,m})$, $i,k=1,
\ldots, d$, belong to $\mathcal S^2\times\break L^2(\langle M,M\rangle, \PP
)$ and in particular solve the following BSDEs for $t\in[0,T]$:
%
%e4.11 ###
%
\begin{eqnarray}
\label{diffBSDEx}
U^x_{kt} &=& \sum_{j=1}^n \partial_{j} F(X_T^{x,m},M_T^m) D^x_{jkT} -
\sum_{\alpha=1}^d \int_t^T V^x_{i\alpha s} \,dM^{(\alpha)}_{s}\nonumber
\\
&&{} + \sum_{j=1}^n \int_t^T \partial_{1+j}
f(s,X_s^{x,m},M_s^{m},Y^{x,m}_s,Z^{x,m}_{s}q_s^*) D^x_{jks} \,dC_s
\nonumber\\
&&{} + \int_t^T \partial_{1+n+d+1}
f(s,X_s^{x,m},M_s^{m},Y^{x,m}_s,Z^{x,m}_{s}q_s^*) U^x_{ks} \,dC_s
\\
&&{}  + \sum_{j=1}^n \int_t^T \partial_{1+n+d+1+j}
f(s,X_s^{x,m},M_s^{m},Y^{x,m}_s,Z^{x,m}_{s}q_s^*)
\nonumber\\
&&\qquad\quad\hspace*{12pt}{}\times \frac{\partial}{\partial x_k}q_{jks} Z^{(j),x,m} \,dC_s, \nonumber
\\
U^m_{kt}
&=& \sum_{j=1}^d \partial_{n+j} F(X_T^{x,m},M_T^m) D^m_{jkT}
- \sum_{\alpha=1}^d \int_t^T V^m_{i\alpha s} \,dM^{(\alpha)}_{s}
\nonumber\\
&&{} +  \sum_{j=1}^n \int_t^T \partial_{1+j}
f(s,X_s^{x,m},M_s^{m},Y^{x,m}_s,Z^{x,m}_{s}) D^m_{jks} \,dC_s
\nonumber\\
&&{}+  \int_t^T \partial_{1+n+k}
f(s,X_s^{x,m},M_s^{m},Y^{x,m}_s,Z^{x,m}_{s}q_s^*) \,dC_s
\nonumber\\
&&{} +  \int_t^T \partial_{1+n+d+1}
f(s,X_s^{x,m},M_s^{m},Y^{x,m}_s,Z^{x,m}_{s}q_s^*) U^m_{ks} \,dC_s
\nonumber\\
&&{} +  \sum_{j=1}^n \int_t^T \partial_{1+n+d+1+j}
f(s,X_s^{x,m},M_s^{m},Y^{x,m}_s,Z^{x,m}_{s}q_s^*) q_{jks} V^m_{jks}
\,dC_s.\nonumber
\end{eqnarray}
\end{prop}
\begin{pf}
From Lemma \ref{lemma:differentiability} and Kolmogorov's lemma (see
\cite{Protter}, Theorem~73, Chapter IV), we deduce that there
exists a family of solutions $(Y^{\tilde x})$ of~(\ref{DiffFBSDE})
which is continuously differentiable in $\tilde x$ for almost all
$\omega\in\Omega$. Finally, from (\ref{dif2}), taking $h
\to0$ the BSDEs follow.
\end{pf}

%s4.2 ###
\subsection{Differentiability of the initial FBSDE}
\label{section:diff2} Now we come back to the system
(\ref{SDE}) and (\ref{BSDE}). In order to obtain the differentiability
of this system we require the following additional assumption:
\begin{longlist}[(MRP)]
\item[(MRP)]
There exists a continuous square-integrable martingale $N:=\break(N_t)_{t\in
[0,T]}$ on $(\Omega, \mathcal{F},\PP)$ which is strongly orthogonal to
$M$ (i.e., $\langle M^i, N \rangle=0$ for $i=1,\ldots,d$) with
$\langle N, N \rangle_T \leq Q,   \PP$-a.s., such that
every $\PP$-martingale is of the form $Z \cdot M + U \cdot N$, where
$Z$ and $U$ are predictable square integrable processes [recall that
$Q$ is the same constant as in (\ref{eq:boundedbrackets})].
\end{longlist}
The presence of the additional bracket $\langle L, L \rangle$ in the
BSDE prevents us from applying the known techniques for
differentiability in the Brownian case as shown in
\cite{AnkirchnerImkellerDosReis1,AnkirchnerImkellerDosReis2,BriandConfortola}.
Nevertheless, under (MRP) we can show that the BSDE in (\ref{BSDE})
can be written as
%
%e4.12 ###
%
\begin{eqnarray}
\label{FBSDEMRP}
Y^{x,m}_t
&=&
F(X_T^{x,m},M_T^m)-\int_t^T Z_r^{x,m} \,dM_r -\int_t^T U_r^{x,m} \,dN_r
\nonumber\\[-8pt]\\[-8pt]
&&{}+ \int_t^T h(r,X_r^{x,m},M_r^m,Y_r^{x,m},Z_r^{x,m}q_r^*,U_r^{x,m})
\,d\tilde C_r,
\nonumber
\end{eqnarray}
$t\in[0,T]$, where $\tilde C$ and $h$ are defined as in Appendix
\ref{Appendix1}.
Due to hypothesis (MRP) and the orthogonality of the martingales $L$
and $M$, the representation of $L$ as $L= U \cdot N$ where $U$ is a
predictable square integrable stochastic process is obtained.
So the solution $(Y,Z,L)$ of the backward part (\ref{BSDE}) becomes
$(Y,Z,U)$ in (\ref{FBSDEMRP}).
The bracket $\langle L, L \rangle$ is then a component of the new
generator $h$, which is quadratic in $U$.
We refer to Appendix \ref{Appendix1}, where a~discussion of the
technical aspects is given.
Now we can write the system~(\ref{SDE})  and~(\ref{BSDE}) as
%
%e4.13 ###
%
\begin{eqnarray}
\label{systemdiff}\quad
X_t^{x,m}&=&x+\int_0^t \tilde\sigma(s,X_s^{t,x,m},M_s^{t,m}) \,d\tilde
M_s + \int_0^t \tilde b(s,X_s^{x,m},M^{m}_s)\,d\tilde C_s,
\nonumber\\
Y_t^{x,m}&=&F(X_T^{x,m},M_T^m)-\int_t^T \tilde Z_s^{x,m} \,d\tilde M_s\\
&&{} +
\int_t^T h(s,X_s^{x,m},M_s^m,Y_s^{x,m},\tilde Z_s^{x,m} \tilde q_s^*)
\,d\tilde C_s, \nonumber
\end{eqnarray}
$t\in[0,T]$, where $\tilde M$, $\tilde q$, $\tilde Z$, are defined as
in Appendix \ref{Appendix1} and $\tilde\sigma:= (
{\sigma} \enskip{0} )$, $\tilde b :=b\times\varphi_1$ where
$\varphi_1$ is a bounded predictable process defined in Appendix \ref
{Appendix1}.
A solution $(X^{x,m},Y^{x,m},\tilde Z^{x,m}) \in\mathcal S^p \times
\mathcal S^{\infty} \times L^2(d\langle\tilde M, \tilde
M \rangle\otimes d\PP)$ of this system exists for $\sigma$, $b$
satisfying (H0) and $F$, $h$ satisfying (H1)--(H3).
Therefore, we obtain the following result, whose proof follows from
Proposition \ref{PropDiff}.
\begin{theorem}
\label{DiffLast}
Assume that $M$ be a strong Markov process and that~$f$ and $F$ in
(\ref{BSDE}) satisfy \textup{(H1)}--\textup{(H3)} and \textup{(D1)}--\textup{(D4)}. Under
the assumption \textup{(MRP)} there exists a\vadjust{\goodbreak}
solution $(X^{\tilde x},Y^{\tilde x},\tilde Z^{\tilde x})$ of
(\ref{SDE})  and (\ref{BSDE}), such that~$X^{\tilde x}(\omega)$ and
$Y^{\tilde x}(\omega)$ are continuously differentiable in $\tilde x
\in\real^{(n+d) \times1}$ for almost all $\omega\in\Omega$ $[$we
recall that $\tilde{x}$ stands for $(x,m)]$.
\end{theorem}
\begin{pf}
Note that the processes $Y^{\tilde{x}}$ of the transformed BSDE (\ref
{FBSDEMRP}) and of the original BSDE (\ref{BSDE}) coincide. In
addition, the process $(Z^{\tilde{x}},L^{\tilde{x}})$ in (\ref
{BSDE}) and the processes $\tilde{Z}^{\tilde{x}}$ in (\ref
{FBSDEMRP}) are related as follows: $\tilde{Z}^{\tilde{x}}=(Z^{\tilde
{x}},U^{\tilde{x}})$ with $L^{\tilde{x}}=\int_0^\cdot U^{\tilde
{x}}_r \,dN_r$ and $N$ is the process coming from (MRP). The definition
of the driver $h$ of the BSDE (\ref{FBSDEMRP}) (see Appendix \ref
{Appendix1}), the fact that $f$ and $F$ in (\ref{BSDE}) satisfy and
(H1)--(H3) and (D1)--(D4), imply that $F$ and $h$ also satisfy the
assumptions\vspace*{1pt} (H1)--(H3) and (D1)--(D4). Thus, $Y^{\tilde{x}}$ and
$Z^{\tilde{x}}$ are continuously differentiable in $\tilde{x}$ by
Proposition \ref{PropDiff} which concludes the proof.
\end{pf}
\begin{prop}
\label{prop:diffuv}
Assume that $M$ is a strong Markov process and that~$f$ and $F$ in
(\ref{BSDE}) satisfy \textup{(H1)}--\textup{(H3)} and \textup{(D1)}--\textup{(D4)}.
From Theorem \ref{theorem:Markovproperty} there exists a
deterministic function $u$ such that
$Y_s^{t,x,m}=u(s,X_s^{t,x,m},M_s^{t,m}),   s\in[t,T]$. Under the
assumption \textup{(MRP)} we have that:
\begin{longlist}
\item $x \mapsto u(t,x,m) \in\mathscr{C}^1(\real^{n\times
1}),  (t,m) \in[0,T]\times\real^{d\times1}$,
\item $m \mapsto u(t,x,m) \in\mathscr{C}^1(\real^{d\times
1}),  (t,x) \in[0,T]\times\real^{n\times1}$,
\item there exist two constants $\zeta_1, \zeta_2$
depending only on $\|F\|_{\infty}$, $a$ and $b$ of assumption
\textup{(H2)}
such that
\[
\zeta_1 \leq u(t,x,m) \leq\zeta_2
\]
for all $(t,x,m) \in[0,T] \times\real^{n\times1} \times\real
^{d\times1}$,
\item the maps
\[
(t,x,m) \mapsto\partial_i u(t,x,m)
\]
are continuous for $i=2,3$.
\end{longlist}
\end{prop}
\begin{pf}
(i) Fix $(t,m)$ in $[0,T]\times\real^{d\times1}$. As already
mentioned, $Y_t^{t,x,m}$ is deterministic
and $u(t,x,m)=Y_t^{t,x,m}$. By differentiability of $Y^{t,x,m}$
with respect to $x$ (Theorem \ref{DiffLast}), we obtain that $x
\mapsto u(t,x,m)$ belongs to $\mathscr{C}^1(\real^{n\times1})$.

\mbox{}\hphantom{i}(ii) The proof is similar to (i).

(iii) Let $(t,x,m) \in[0,T] \times\real^{n\times1} \times
\real^{d\times1}$. By \cite{Morlais1}, Lemma 3.1(i), there exists~%
$\zeta_1, \zeta_2$ depending only on $|F|_{\infty}$, $a$ and $b$
such that $ \zeta_1 \leq Y_s^{t,x,m} \leq\zeta_2 $ for all $s$ in
$[t,T]$, $\PP$-a.s. Thus, $\zeta_1 \leq u(t,x,m)=Y_t^{t,x,m} \leq
\zeta_2$. Since the constants~$\zeta_1$ and~$\zeta_2$ do not depend
on $(t,x,m)$, the claim is proved.

(iv) For better readability we prove this claim for $d=n=1$.
The multidimensional case is a straightforward extension of the
following computations where we adapt \cite{MaZhang}, Theorem 3.1.
From (\ref{EqMarkovPropertyQuad}) we know that
$Y_s^{t,x,m}=u(s,X_s^{t,x,m},M_s^{t,m})$ and hence, $Y_t^{t,x,m}=u(t,x,m)$.
In the following we use the representation (\ref{systemdiff}) of the
forward--backward system, that is, we use the transformed FBSDE. Then by
definition of the driver $h$ (see Appendix~\ref{Appendix1}) the
properties\vadjust{\goodbreak} of $f$ carry over to $h$. Thus, by Proposition \ref
{PropDiff}, the processes $(\nabla_x Y^{t,x,m}, \nabla_x\tilde
Z^{t,x,m})$ satisfy the following BSDE:
\begin{eqnarray*}
\nabla_x Y_s^{t,x,m}&=& \nabla_x F(X_T^{t,x,m},M_T^{t,m}) \nabla_x
X_T^{t,x,m}-\int_s^T \nabla_x \tilde Z_r^{t,x,m} \,d\tilde{M}_r \\[-1pt]
&&{}+ \int_s^T \bigl( \partial_2 h (r,\Theta_r(t,x,m)) \nabla_x X_r^{t,x,m}
\\[-1pt]
&&\hphantom{{}+ \int_s^T \bigl(}
{} + \partial_4 h (r,\Theta_r(t,x,m)) \nabla_x Y_r^{t,x,m}\\[-1pt]
&&\hphantom{{}+ \int_s^T \bigl(}
{} + \partial
_5 h (s,\Theta_r(t,x,m)) \tilde q_r \nabla_x \tilde Z_r^{t,x,m}
\bigr) \,d\tilde C_s.
\end{eqnarray*}
Thus, putting $s=t$ in the above expression and taking the expectation
we get
\begin{eqnarray*}
&&\partial_x u(t,x,m)
\\[-1pt]
&&\qquad=\E\biggl[ \nabla_x F(X_T^{t,x,m},M_T^{t,m}) \nabla_x
X_T^{t,x,m}\\[-1pt]
&&\qquad\quad\hphantom{\E\biggl[}
{}
+ \int_t^T \bigl( \partial_2 h (s,\Theta_u(t,x,m)) \nabla_x
X_s^{t,x,m}\\[-1pt]
&&\qquad\quad\hspace*{44.5pt}{}+ \partial_4 h (s,\Theta_s(t,x,m)) \nabla_x Y_s^{t,x,m}\\[-1pt]
&&\hspace*{44.5pt}\qquad\quad{} + \partial
_5 h (s,\Theta_s(t,x,m)) \tilde q_s \nabla_x \tilde Z_s^{t,x,m} \bigr)
\,d\tilde C_s \biggr] .
\end{eqnarray*}
Here we have used $\Theta_s(t,x,\tilde m):=(X_s^{t,x,m},\tilde M_s^{t,
m},Y_s^{t,x, m},\tilde Z_s^{t,x, m} \tilde q_s)$.
Let us fix $(t_1,x_1, m_1)$ and $(t_2,x_2,m_2)$ with $t_1<t_2$ and
denote $\Theta_s^1:=\Theta_s^1(t_1,x_1, m_1)$ and $\Theta
_s^2:=\Theta_s^2(t_2,x_2,m_2)$.
We write $X^1:=X^{t_1,x_1,m_1}$ and analogously $X^2$, $Y^1$, $Y^2$, etc.
Furthermore, we define $\Delta_{1,2}\varphi(s):=\varphi(s,\Theta
_s^1)-\varphi(s,\Theta_s^2)$ for any function~$\varphi$ with values
in $\real$.
We have that
\begin{eqnarray*}
&&|\partial_x u(t_1,x_1,m_1)-\partial_x u(t_2,x_2,m_2)|
\\[-1pt]
&&\qquad\leq\E[ \nabla_x F(X_T^1,M_T^1) \nabla_x X_T^1- \nabla_x
F(X_T^2,M_T^2) \nabla_x X_T^2 ]
\\[-1pt]
&&\qquad\quad{}  + \E\biggl[\int_{t_1}^{t_2} |\partial_2 h (s,\Theta_s^1)|
|\nabla_x X_s^1|\\[-1pt]
&&\qquad\quad\hspace*{23pt}{}+ |\partial_4 h (s,\Theta_s^1)| |\nabla_x Y_s^1| +
|\partial_5 h (s,\Theta_s^1)| |\tilde q_s \nabla_x \tilde Z_s^1|
\,d\tilde C_s \biggr]
\\[-1pt]
&&\qquad\quad{}  +\E\biggl[ \int_{t_2}^T |\Delta_{1,2} \partial_2 h(s)| |\nabla
_x X_s^1| \\[-1pt]
&&\qquad\quad\hspace*{23pt}{}+ |\Delta_{1,2} \partial_4 h(s)| |\nabla_x Y_s^1| +
|\Delta_{1,2} \partial_5 h(s)| |\tilde q_s \nabla_x \tilde Z_s^1|\,d
\tilde C_s \biggr]
\\[-1pt]
&&\qquad\quad{}  +\E\biggl[ \int_{t_2}^T |\partial_2 h (s,\Theta_s^2)| |\nabla
_x X^1_s-\nabla_x X^2_s|\\[-1pt]
&&\qquad\quad\hspace*{23pt}{} + |\partial_4 h (s,\Theta_s^2)| |\nabla_x
Y^1_s-\nabla_x Y^2_s| \,d\tilde C_s \biggr]
\\
&&\qquad\quad{}  +\E\biggl[ \int_{t_2}^T |\partial_5 h (s,\Theta_s^2)| |\tilde
q_s(\nabla_x \tilde Z_s^1-\nabla_x \tilde Z_s^2)| \,d\tilde C_s \biggr]
\\
&&\qquad \leq\E[ |\nabla_x F(X_T^1,M_T^1) - \nabla_xF(X^2_T,M_T^2)|
|\nabla_x X_T^1| \\
&&\qquad\quad\hphantom{\E[}\hspace*{24pt}
{}+ |\nabla_x F(X_T^2,M_T^2)| |\nabla_x X_T^1 -
\nabla_x X_T^2| ]
\\
&&\qquad\quad{}  + c \E\biggl[\int_{t_1}^{t_2} (|\tilde q_s \theta_s|+|\tilde
q_s \tilde Z_s^{1}|)(|\nabla_x X_s^1|+ |\nabla_x Y_s^1| +
|\tilde q_s \nabla_x \tilde Z_s^1|) \,d\tilde C_s \biggr]
\\
&&\qquad\quad{}  +c \E\biggl[ \int_{t_2}^T (|\tilde q_s \theta_s|+|\tilde
q_s \tilde Z_s^1|+|\tilde q_s \tilde Z_s^2|)\\
&&\qquad\quad\hphantom{{}  + c \E\biggl[}
\hspace*{14.1pt}{}\times(|X_s^1-X_s^2|+|
M_s^1- M_s^2|+|Y_s^1-Y_s^2|)\\
&&\qquad\quad\hspace*{114.1pt}{} \times(\vert\nabla_x X_s^1 \vert+ \vert
\nabla_x Y_s^1 \vert) \,d\tilde C_s \biggr]
\\
&&\qquad\quad{}  +c \E\biggl[ \int_{t_2}^T (|\tilde q_s \theta_s|+|\tilde q_s
\tilde Z_s^1|+|\tilde q_s \tilde Z_s^2|)\\
&&\qquad\quad\hphantom{{}  +c \E\biggl[ \int_{t_2}^T}
{}\times \vert\tilde q_s(\tilde
Z_s^1-\tilde Z_s^2)\vert(\vert\nabla_x X_s^1 \vert+ \vert\nabla
_x Y_s^1 \vert) \,d\tilde C_s \biggr]
\\
&&\qquad\quad{}  +c \E\biggl[ \int_{t_2}^T |\tilde q_s \nabla_x \tilde Z_s^1|
(|X_s^1-X_s^2|+| M_s^1- M_s^2|+|Y_s^1-Y_s^2|) \,d\tilde C_s \biggr]
\\
&&\qquad\quad{}  +c \E\biggl[ \int_{t_2}^T |\tilde q_s \nabla_x \tilde Z_s^1|
|\tilde q_s(\tilde Z_s^1- \tilde Z_s^2)| \,d\tilde C_s \biggr]
\\
&&\qquad\quad{}  +c \E\biggl[ \int_{t_2}^T (|\tilde q_s \theta_s|+|\tilde
q_s \tilde Z_s^2|) \\
&&\qquad\quad\hphantom{{}  +c \E\biggl[ \int_{t_2}^T}
{}  \times\bigl(|\nabla_x X_s^1-\nabla_x X_s^2| +
|\nabla_x Y_s^1-\nabla_x Y_s^2|\\
&&\qquad\quad\hspace*{117pt}{}+|\tilde q_s(\nabla_x \tilde
Z_s^1-\nabla_x \tilde Z_s^2)|\bigr) \,d\tilde C_s\biggr]
\\
&&\qquad=:\sum_{i=1}^7 T_i,
\end{eqnarray*}
where we have used the assumptions (D3) and (D4) in the last inequality.
Recall that $(t_1,x_1, m_1)$ is fixed and $t_2 >t_1$.
With (\ref{estX}) and (\ref{estF}) we see
\[
\lim_{t_2\to t_1; x_2\to x_1; m_2\to m_1} T_1=\lim_{x_2\to x_1} T_1=0.
\]
By the monotone convergence theorem we deduce for the second term
\[
\lim_{t_2\to t_1; x_2\to x_1; m_2\to m_1} T_2 = \lim_{t_2\to t_1}
T_2 =0.
\]
We now deal with $T_3$,
\begin{eqnarray*}
T_3&\leq& c \E\biggl[ \sup_{s\in[0,T]} (\vert\nabla_x X_s^1
\vert+ \vert\nabla_x Y_s^1 \vert) \\
&&\hphantom{c \E\biggl[\sup_{s\in[0,T]}}
{}\times\int_{t_2}^T (|\tilde q_s
\theta_s|+|\tilde q_s\tilde Z_s^1|+|\tilde q_s\tilde Z_s^2|)\\
&&\qquad\quad\hphantom{\sup_{s\in[0,T]}}\hspace*{11.3pt}{}\times(|X_s^1-X_s^2|+| M_s^1- M_s^2|+|Y_s^1-Y_s^2|)
\,d\tilde C_s \biggr]
\\
&\leq& c \E\Bigl[ \sup_{s\in[0,T]} (\vert\nabla_x X_s^1 \vert+
\vert\nabla_x Y_s^1 \vert)^2\Bigr]^{1/2}
\\
&&{} \times\E\biggl[\biggl(\int_{t_2}^T (|\tilde q_s \theta_s|+|\tilde
q_s \tilde Z_s^1|+|\tilde q_s \tilde Z_s^2|)\\
&&\hspace*{46pt}{}\times (|X_s^1-X_s^2|+| M_s^1-
M_s^2|+|Y_s^1-Y_s^2|) \,d\tilde C_s\biggr)^2 \biggr]^{1/2}\\
&\leq& c \E\biggl[\int_{t_2}^T (|\tilde q_s \theta_s|+|\tilde q_s \tilde
Z_s^1|+|\tilde q_s\tilde Z_s^2|)^2 \,d\tilde C_s \\
&&\hspace*{16.1pt}{}\times\int_{t_2}^T
(|X_s^1-X_s^2|+| M_s^1- M_s^2|+|Y_s^1-Y_s^2|)^2 \,d\tilde C_s \biggr]^{1/2}\\
&\leq& c \E\biggl[\biggl(\int_{t_2}^T (|\tilde q_s \theta_s|+|\tilde
q_s\tilde Z_s^1|+|\tilde q_s \tilde Z_s^2|)^2 \,d\tilde C_s\biggr)^2
\biggr]^{{1/4}}\\
&&{} \times\E\biggl[\biggl(\int_{t_2}^T (|X_s^1-X_s^2|+|
M_s^1- M_s^2|+|Y_s^1-Y_s^2|)^2 \,d\tilde C_s\biggr)^2 \biggr]^{{1/4}}\\
&\leq& c \bigl(\tilde{c}+(|x_2|^2+| m_2|^2)^{p_1} (|x_2-x_1|^2+| m_2- m_1|^2)^{p_2}\bigr),
\end{eqnarray*}
where we have used the Cauchy--Schwarz inequality.
Here $p_1, p_2$ are two positive numbers given by the a priori
estimates Lemma \ref{appendix:apriori} and $\tilde{c}$ is a~positive
constant. Thus, we conclude
\[
\lim_{t_2\to t_1; x_2\to x_1; m_2\to m_1} T_3=\lim_{x_2\to x_1;
m_2\to m_1} T_3=0.
\]
Similarly, one shows
\[
\lim_{t_2\to t_1; x_2\to x_1; m_2\to m_1} T_5=0.
\]
We now estimate $T_4$ and $T_7$ but we give the details only for $T_4$,
since those for $T_7$ follow the same lines.
Applying the Cauchy--Schwarz inequality again we get
\begin{eqnarray*}
&&\E\biggl[ \int_{t_2}^T (|\tilde q_s \theta_s|+|\tilde q_s \tilde
Z_s^1|+|\tilde q_s \tilde Z_s^2|) \vert\tilde q_s(\tilde Z_s^1-\tilde
Z_s^2)\vert(\vert\nabla_x X_s^1 \vert+ \vert\nabla_x Y_s^1 \vert
) \,d\tilde C_s \biggr]
\\
&&\qquad\leq\E\biggl[\sup_{s\in[0,T]} (\vert\nabla_x X_s^1 \vert+ \vert
\nabla_x Y_s^1 \vert) \\
&&\qquad\quad\hphantom{\E\biggl[\sup_{s\in[0,T]}}
{}\times\int_{t_2}^T (|\tilde q_s \theta_s|+|\tilde
q_s\tilde Z_s^1|+|\tilde q_s \tilde Z_s^2|) \vert\tilde q_s(\tilde
Z_s^1-\tilde Z_s^2)\vert \,d\tilde C_s \biggr]
\\
&&\qquad\leq\E\Bigl[\sup_{s\in[0,T]} (\vert\nabla_x X_s^1 \vert+ \vert
\nabla_x Y_s^1 \vert)^2 \Bigr]^{1/2}\\
&&\qquad\quad\hphantom{\E\biggl[}
{} \times\E\biggl[\biggl(\int_{t_2}^T (|\tilde q_s \theta
_s|+|\tilde q_s \tilde Z_s^1|+|\tilde q_s \tilde Z_s^2|) \vert\tilde
q_s(\tilde Z_s^1-\tilde Z_s^2)\vert \,d\tilde C_s\biggr)^2 \biggr]^{1/2}
\\
&&\qquad\leq c \E\biggl[\int_{t_2}^T (|\tilde q_s \theta_s|+|\tilde q_s \tilde
Z_s^1|+|\tilde q_s\tilde Z_s^2|)^2 \,d\tilde C_s \int_{t_2}^T \vert
\tilde q_s(\tilde Z_s^1-\tilde Z_s^2)\vert^2 \,d\tilde C_s \biggr]^{1/2}
\\
&&\qquad\leq c \E\biggl[\biggl(\int_{t_2}^T (|\tilde q_s \theta_s|+|\tilde q_s
\tilde Z_s^1|+|\tilde q_s \tilde Z_s^2|)^2 \,dC_s\biggr)^2\biggr]^{
{1/4}} \\
&&\qquad\quad{}\times\E\biggl[\biggl(\int_{t_2}^T \vert\tilde q_s( \tilde Z_s^1- \tilde
Z_s^2)\vert^2 \,d \tilde C_s\biggr)^2 \biggr]^{{1/4}}
\\
&&\qquad\leq c \bigl(\tilde{c}+(|x_2|^2+| m_2|^2)^{p_1} (|x_2-x_1|^2+| m_2- m_1|^2)^{p_2}\bigr).
\end{eqnarray*}
Here, as before, $p_1, p_2$ are two positive numbers given by the a
priori estimates Lemma \ref{appendix:apriori} and $\tilde{c}$ is a
positive constant.
This leads to
\[
\lim_{t_2\to t_1; x_2\to x_1; m_2\to m_1} T_4=\lim_{x_2\to x_1;
m_2\to m_1} T_4=0.
\]
Finally, we consider the term $T_6$
\begin{eqnarray*}
&&\E\biggl[ \int_{t_2}^T | \tilde q_s \nabla_x \tilde Z_s^1| |\tilde
q_s(\tilde Z_s^1-\tilde Z_s^2)| \,d\tilde C_s \biggr]
\\
&&\qquad\leq\E\biggl[ \biggl(\int_{t_2}^T |\tilde q_s \nabla_x \tilde Z_s^1|^2
\,d\tilde C_s\biggr)^{1/2} \biggl( \int_{t_2}^T |\tilde q_s(\tilde
Z_s^1-\tilde Z_s^2)|^2 \,d\tilde C_s \biggr)^{1/2}\biggr]
\\
&&\qquad\leq c \E\biggl[ \int_{t_2}^T |\tilde q_s(\tilde Z_s^1-\tilde Z_s^2)|^2
\,d\tilde C_s \biggr]^{1/2}
\\
&&\qquad\leq c (\vert x_2-x_1\vert^2 + \vert m_2-m_1\vert^2)^p,
\end{eqnarray*}
where the positive constant $p$ is given by the a priori estimates
Lemma~\ref{appendix:apriori}.
Thus, we have
\[
\lim_{t_2\to t_1; x_2\to x_1; m_2\to m_1} T_6=0.
\]
The same methodology shows that for fixed $(t_2,x_2,m_2)$
\[
\lim_{t_1\to t_2; x_1\to x_2; m_1\to m_2} |\partial_x
u(t_1,x_1,m_1)-\partial_x u(t_2,x_2,m_2)|=0.
\]
Similarly, we can show that $\partial_m u$ is continuous in $(t,x,m)$.
\end{pf}

\subsection*{Example of stochastic basis where the condition \textup{(MRP)} is
satisfied}

Let $(B^1,B^2):=(B_s^1,B_s^2)_{s \in[0,T]}$ be a two-dimensional
Brownian motion defined on a probability space $(\Omega,\mathcal
{F},\PP)$ with a terminal time $0<T<\infty$ and with $B^1$ and $B^2$
being independent.
We denote by $(\mathcal{F}_t)_{t\in[0,T]}$ the filtration generated
by $(B^1,B^2)$.
Then the process $M:=(B_t^1)_{t\in[0,T]}$ is a continuous martingale
with respect to $(\mathcal{F}_t)_{t\in[0,T]}$ and it is a $(\mathcal
{F}_t)_{t\in[0,T]}$-strong Markov process.
Let $N=(B_t^2)_{t\in[0,T]}$.
The martingale representation property for $(B^1,B^2)$ and the strong
orthogonality between $B^1$ and $B^2$ entail that the pair $(M,N)$
satisfies the property (MRP) introduced in Section \ref{section:diff2}.

%s5 ###
\section{Representation formula}
\label{mainsection}

In this section we provide the representation formula (\ref{eqn3})
which generalizes the one obtained in
\cite{AnkirchnerImkellerDosReis1,AnkirchnerImkellerDosReis2}, where
$M$ is a~Brownian motion.
We recall that in the Gaussian setting the proof of this formula is
based on the representation of the stochastic process $Z$ as the trace
of the Malliavin derivative of $Y$.
In the general martingale setting of this paper, Malliavin's calculus
is not available, therefore, we propose a~new
proof based on stochastic calculus techniques.
We also stress that
the last term in formula (\ref{eqn3}) vanishes if we assume that $M$
has independent increments, $\sigma$ and $b$ do not depend on $M$ in
(\ref{SDE}) and that the driver $f$ in~(\ref{BSDE}) is independent of
$M$.

We present the main result of this paper. We stress that this
result does not rely on the assumption (MRP) made in Section \ref
{section:diff2} since only the regularity of the deterministic function
$u$ where $Y=u(\cdot,X,M)$ is needed.
\begin{theorem}
\label{theorem:main}
Assume that $M$ is a Markov process. Assume that \textup{(H0)},
\textup{(H1)}--\textup{(H3)} are in force for the FBSDE (\ref{SDE}) and (\ref
{BSDE}). Then by Theorem \ref{theorem:Markovproperty}, there exists
a deterministic function $u$ such that
$Y_s^{t,x,m}=u(s,X_s^{t,x,m},M_s^{t,m}),\break   s\in[t,T]$. Assume, in
addition, that $u$ satisfies:
\begin{longlist}
\item $x \mapsto u(t,x,m)$, $(t,m) \in[0,T]\times\real
^{d\times1}$, is continuously differentiable,
\item $m \mapsto u(t,x,m)$, $(t,x) \in[0,T]\times\real
^{n\times1}$, is continuously differentiable,
\item there exist two constants $\zeta_1, \zeta_2$
depending only on $\|F\|_{\infty}$, $a$ and $b$ of~assumption
\textup{(H2)} such that
\[
\zeta_1 \leq u(t,x,m) \leq\zeta_2\qquad  \forall(t,x,m) \in[0,T]
\times\real^{n\times1} \times\real^{d\times1},\vadjust{\goodbreak}
\]
\item the maps
\[
(t,x,m) \mapsto\partial_i u(t,x,m)  \mbox{ are continuous}\qquad\mbox{for
} i=2,3.
\]
\end{longlist}
Then for all $s \in[t,T]$ we have $\nu$-a.e.
%
%e5.1 ###
%
\begin{eqnarray}
\label{eq:main} Z_s^{t,x,m}&=&
\partial_2 u(s,X_s^{t,x,m},M_s^{t,m})
\sigma(s,X_s^{t,x,m},M_s^{t,m})\nonumber\\[-8pt]\\[-8pt]
&&{}+ \partial_3 u(s,X_s^{t,x,m},M_s^{t,m}).\nonumber
\end{eqnarray}
\end{theorem}
\begin{remark}
\label{indIncbis}
(i) An interesting particular case of Theorem \ref
{theorem:main} is given when $X$ and $M$ are as in Proposition
\ref{CinlaretalProtter}(ii) and when $f$ in (\ref{BSDEbis})
does not depend on $M$. In this situation, equation (\ref{eq:main}) becomes
\[
Z_s^{t,x}=\partial_2 u(s,X_s^{t,x}) \sigma(s,X_s^{t,x}),  \qquad\nu\mbox{-a.e.},
\]
which coincides with the representation formula derived in
\cite{AnkirchnerImkellerDosReis1,AnkirchnerImkellerDosReis2} when $M$
is a
standard Brownian motion.

(ii)  One may be interested in knowing when $u$ in Theorem \ref
{theorem:main} does not depend trivially on $M$, that is, when the
third term in (\ref{eq:main}) does not vanish. This is related to the
Markov property given for $Y$ and we provide in Appendix \ref
{Appendix2} an explicit example where $u$ depends nontrivially on $M$.
\end{remark}
\begin{pf*}{Proof of Theorem \ref{theorem:main}}
Fix $s$ in $[t,T]$. For simplicity of notation we drop the superscript
$(t,x,m)$. We briefly explain the idea of the proof. Assume that the
function $u$ introduced above is in $\mathscr{C}^{1,2,2}$ that is
continuously differentiable in time and twice continuously
differentiable in $(x,m)$. Then an application of It\^{o}'s formula
gives that
%
%e5.2 ###
%
\begin{eqnarray}
\label{bracket1}\qquad\quad
\langle Y, M_\cdot\rangle_s&=&\langle u(\cdot,X_\cdot,M_\cdot),
M_\cdot\rangle_s\nonumber\\[-8pt]\\[-8pt]
&=& \int_t^s [\partial_2u(r,X_r,M_r) \sigma(r,X_r,M_r) +\partial_3
u(r,X_r,M_r)] \,d\langle M, M \rangle_r,\nonumber
\end{eqnarray}
where we denote by $\langle u(\cdot,X_s,M_s), M_\cdot\rangle_s$ the
covariation vector
\[
\bigl(\bigl\langle u(\cdot,X_\cdot,M_\cdot),
M_\cdot^{(1)}\bigr\rangle_s,\ldots,\bigl\langle u(\cdot,X_\cdot,M_\cdot),
M_\cdot^{(d)}\bigr\rangle_s\bigr).
\]
Then, since $(Y,Z)$ is solution of
(\ref{BSDE}), we have that
%
%e5.3 ###
%
\begin{equation}
\label{bracket1bis}
\langle Y, M\rangle_s=\int_t^s Z_r \,d\langle M, M \rangle_r, \qquad  s
\in[t,T].
\end{equation}
The conclusion of the theorem then follows from the fact that
$Y_s=u(s,X_s,\allowbreak M_s)$, $s \in[t,T]$ and from relations (\ref{bracket1}) and
(\ref{bracket1bis}). However, we have assumed the function $u$ to be
much more regular than what it is and so we have to prove the relation
(\ref{bracket1}) for $u$ being only one time differentiable in
$(x,m)$. The rest of the proof is devoted to this fact. For this we
compute ``directly'' the quadratic variation between $u(\cdot,X_\cdot
,M_\cdot)$ and $M$.

Fix $i \in\{1,\ldots, d\}$. Let $r\geq1$ and $\pi
^{(r)}:=\{t_j^{(r)},   j=1, \ldots, r\}$ be a partition of $[t,T]$
whose mesh size $\vert\pi^{(r)} \vert$ tends to zero as $r$ goes to
infinity with
$t_0^{(r)}=t$ and $t_{r}^{(r)}=T$ such that
\begin{eqnarray*}
&&\lim_{r\to\infty} \sup_{t\leq s \leq T} \Biggl\vert\bigl\langle
u(\cdot,X_\cdot,M_\cdot), M_\cdot^{(i)}\bigr\rangle_s\\[-2pt]
&&\qquad\quad\hspace*{16.1pt}{}-\sum
_{j=0}^{\varphi_s-1}
\bigl(u\bigl(t_{j+1}^{(r)},X_{t_{j+1}^{(r)}},M_{t_{j+1}^{(r)}}\bigr)-u\bigl(t_j^{(r)},X_{t_j^{(r)}},M_{t_j^{(r)}}\bigr)
\bigr) \Delta_j M^{(i)} \Biggr\vert\\[-2pt]
&&\qquad=0,
\end{eqnarray*}
where the limit is understood in probability with respect to $\PP$,
$\Delta_j M$ denotes the
increments of the stochastic process $M$ on $[t_j^{(r)},t_{j+1}^{(r)}]$
and $\varphi_s^{(r)}$ is such that $\varphi_s^{(r)}=j$ with
$t_j^{(r)} \leq s <t_{j+1}^{(r)}$. For simplicity\vspace*{2pt} of notation, the
superscript~$(r)$ will be omitted. In addition, up to a subsequence we
can assume that convergence above is almost sure with respect to $\PP$.
We have that
%
%e5.4 ###
%
\begin{eqnarray}
\label{bracketmain}
&&\bigl\langle u(\cdot,X_\cdot,M_\cdot), M_\cdot^{(i)}\bigr\rangle_s
\nonumber\\[-2pt]
&&\qquad=\lim_{r \to\infty} \sum_{j=0}^{\varphi_s-1}
\bigl(u(t_{j+1},X_{t_{j+1}},M_{t_{j+1}})-u(t_j,X_{t_j},M_{t_j})\bigr)
\Delta_j M^{(i)}\nonumber\\[-2pt]
&&\qquad=\lim_{r \to\infty} \Biggl[\sum_{j=0}^{\varphi_s-1}
\bigl(u(t_{j+1},X_{t_{j}},M_{t_{j}})-u(t_j,X_{t_{j}},M_{t_{j}})\bigr)
\Delta_j M^{(i)}\\[-2pt]
&&\hspace*{26.4pt}\qquad\quad{}+\sum_{j=0}^{\varphi_s-1}
\bigl(u(t_{j+1},X_{t_{j+1}},M_{t_{j+1}})-u(t_{j+1},X_{t_j},M_{t_j})\bigr)
\Delta_j M^{(i)} \Biggr]\nonumber\\[-2pt]
&&\qquad=:\lim_{r \to\infty} \bigl[ S_{s,r,1}^{(i)} + S_{s,r,2}^{(i)}
\bigr].\nonumber
\end{eqnarray}
We treat the two parts separately. First, assume that the second term
converges, more precisely, that relation (\ref{main:conv1}) below
holds:
%
%e5.5 ###
%
\begin{eqnarray}
\label{main:conv1}
&&\lim_{r \to\infty} \sup_{t\leq s \leq T} \biggl\vert
S_{s,r,2}^{(i)}- \biggl(\int_t^s [\partial_2u(r,X_r,M_r) \sigma
(r,X_r,M_r)\nonumber\\[-2pt]
&&\qquad\quad\hspace*{123pt}{} +\partial_3 u(r,X_r,M_r)] \,d\langle M, M \rangle_r
\biggr)^{(i)} \biggr\vert\\[-2pt]
&&\qquad=0,   \qquad\PP\mbox{-a.s.}\nonumber
\end{eqnarray}
It then follows by relations (\ref{bracket1bis}) and (\ref
{bracketmain}) that
\[
\lim_{r \to\infty} \sup_{t\leq s \leq T}
\bigl|S_{s,r,1}^{(i)}-P_s\bigr|=0,  \qquad\PP\mbox{-a.s.}\vadjust{\goodbreak}
\]
with
\begin{eqnarray}
P_s:=\biggl(\int_t^s\! Z_a-\partial_2u(a,X_a,M_a) \sigma(a,X_a,M_a) -
\partial_3 u(a,X_a,M_a) \,d\langle M, M
\rangle_a\biggr)^{(i)},\nonumber\\
&&\eqntext{s \in[t,T].}
\end{eqnarray}
We will show that $P$ is $\PP$-a.s. identically equal to zero.
Since $u$ is not differentiable in time, one can a priori not say how
the sum $S_{s,r,1}^{(i)}$ behaves asymptotically. However, we know that
it converges and that its limit is absolutely continuous with respect
to $d\langle M, M \rangle$. Heuristically, this means that each term
of the form $u(t_{j+1},X_{t_{j}},M_{t_{j}})-u(t_j,X_{t_{j}},M_{t_{j}})$
behaves like a~process times an increment of $\Delta_j M^{(i)}$ which
is not possible since $u$ is a~deterministic function. We will show
that $P$ is a local martingale. Since by definition it is a finite
variation process, we will have $P=0$. We first make the following
assumption that we will relax later. Assume that
%
%e5.6 ###
%
\begin{equation}
\label{eq:finiteness}
\E[ \vert P_s \vert]<\infty\qquad  \forall
s \in[t,T].
\end{equation}
Now fix $t\leq s_1 \leq s_2 \leq T$. For a point $t_j$ in the
subdivision considered above we define $\delta_j
u:=u(t_{j+1},X_{t_{j}},M_{t_{j}})-u(t_j,X_{t_{j}},M_{t_{j}})$. We have that
%
%e5.7 ###
%
\begin{eqnarray}
\label{eq:proofmart}
\E[P_{s_2}|\mathcal{F}_{s_1}] &=& \E\Biggl[ \lim_{r\to\infty} \sum
_{j=0}^{\varphi_{s_2}-1} \delta_j u \Delta_j M^{(i)}
\Big|\mathcal{F}_{s_1}\Biggr]\nonumber\\[-8pt]\\[-8pt]
&=&\E\Biggl[ \lim_{r\to\infty} \sum_{j=0}^{\varphi_{s_2}-1} \delta
_j u \Delta_j M^{(i)}+ (M_{s_2}-M_{\varphi_{s_2}}) \Big|\mathcal
{F}_{s_1}\Biggr]\nonumber
\end{eqnarray}
since by\vspace*{1pt} continuity of the martingale $M$, $\lim_{r\to\infty}
M_{s_2}-M_{\varphi_{s_2}}=0$, $\PP$-a.s. (recall that $\varphi_{s_2}$
tends to $s_2$ when $r$ goes to infinity).
In addition,\vspace*{-1pt} since the function $u$ is bounded [by Proposition
\ref{prop:diffuv}(iii)], the sequence $(\sum_{j=0}^{\varphi
_{s_2}-1} \delta_j u \Delta_j M^{(i)}+ (M_{s_2}-M_{\varphi
_{s_2}}))_r$ is uniformly bounded. Indeed, we have that
\begin{eqnarray*}
&&\E\Biggl[\Biggl|\sum_{j=0}^{\varphi_{s_2}-1} \delta_j u \Delta_j
M^{(i)}+ (M_{s_2}-M_{t_{\varphi_{s_2}}})\Biggr|^2\Biggr]\\
&&\qquad=\sum_{j=0}^{\varphi_{s_2}-1} \E\bigl[ |\delta_j u|^2 \bigl|\Delta_j
M^{(i)}\bigr|^2 \bigr]+ \E[|M_{s_2}-M_{t_{\varphi_{s_2}}}|^2]\\
&&\qquad\leq c \Biggl( \sum_{j=0}^{\varphi_{s_2}-1} \E\bigl[
\bigl|M^{(i)}_{t_{j+1}}\bigr|^2 \bigr]-\E\bigl[ \bigl|M^{(i)}_{t_{j}}\bigr|^2 \bigr]+
\E[|M_{s_2}|^2]-\E[|M_{t_{\varphi_{s_2}}}|^2] \Biggr)\\
&&\qquad=c (\E[|M_{s_2}|^2]-m),
\end{eqnarray*}
thus, \mbox{$\sup_r \E[|\sum_{j=0}^{\varphi_{s_2}-1} \delta
_j u \Delta_j M^{(i)}+ (M_{s_2}-M_{t_{\varphi_{s_2}}})
|^2] \leq c (\E[|M_{s_2}|^2]-m) <\infty$}. Using the Lebesgue
dominated convergence theorem in (\ref{eq:proofmart}) we get
\begin{eqnarray*}
\E[P_{s_2}|\mathcal{F}_{s_1}]&=&\lim_{r\to\infty} \E\Biggl[\sum
_{j=0}^{\varphi_{s_2}-1} \delta_j u \Delta_j M^{(i)}+
\bigl(M_{s_2}^{(i)}-M_{t_{\varphi_{s_2}}}^{(i)}\bigr) \Bigr|\mathcal
{F}_{s_1}\Biggr]\\[-0.5pt]
&=&\lim_{r\to\infty} \Biggl( \sum_{j=0}^{\varphi_{s_1}-1} \delta_j
u \Delta_j M^{(i)} + \E\bigl[(\delta_{\varphi_{s_1}}u) \Delta
_{\varphi_{s_1}} M^{(i)} \vert\mathcal{F}_{s_1} \bigr]\\[-0.5pt]
&&\hphantom{\lim_{r\to\infty} \Biggl(}
{} + \E\Biggl[\sum_{j=\varphi_{s_1}+1}^{\varphi_{s_2}-1} \delta_j u
\Delta_j M^{(i)} + \bigl(M_{s_2}^{(i)}-M_{t_{\varphi_{s_2}}}^{(i)}\bigr)
\Bigr|\mathcal{F}_{s_1}\Biggr] \Biggr)\\[-0.5pt]
&=&\lim_{r\to\infty} \Biggl( \sum_{j=0}^{\varphi_{s_1}-1} \delta_j
u \Delta_j M^{(i)} + (\delta_{\varphi_{s_1}}u)
\bigl(M_{s_1}^{(i)}-M_{t_{\varphi_{s_1}}}^{(i)}\bigr)\Biggr)\\[-0.5pt]
&=&P_{s_1}.
\end{eqnarray*}
Thus, $P$ is a martingale which has (by definition) finite variation,
so it has zero quadratic variation and hence,
\[
P_s=0\qquad \forall s\in[t,T],
\]
which proves
\[
\lim_{r \to\infty} \sup_{t\leq s \leq T} \bigl\vert S_{s,r,1}^{(i)}
\bigr\vert=0, \qquad\PP\mbox{-a.s.}
\]
Now we have to relax the assumption (\ref{eq:finiteness}). Since $P$
is a continuous semimartingale by definition there exists a sequence of
stopping times $(T_m)_m$ with $\lim_{m\to\infty} T_m=T$, $\PP$-a.s.
such that $(P_{s\wedge T_m})_{s\in[t,T]}$ is integrable for all $m\geq
1$. Using this localization, the previous argument leads to
$ P_{s\wedge T_m}=0$ for all $s\in[t,T]$, $\PP$-a.s. By letting $m$
go to infinity we get
\[
\lim_{r \to\infty} \sup_{t\leq s \leq T} \bigl\vert S_{s,r,1}^{(i)}
\bigr\vert=0,\qquad\PP\mbox{-a.s.}
\]
It remains to show that relation (\ref{main:conv1}) holds.
Let $s\in[t,T]$. We have that
%
%e5.8 ###
%
\begin{eqnarray}
\label{Sr2suite}
\lim_{r \to\infty} S_{s,r,2}^{(i)}&=&\lim_{r \to\infty}\sum
_{j=0}^{\varphi_s-1}
\bigl(u(t_{j+1},X_{t_{j+1}},M_{t_{j+1}})-u(t_{j+1},X_{t_j},M_{t_j})\bigr)
\Delta_j M^{(i)}\nonumber\\
&=&\lim_{r \to\infty} \Biggl[\sum_{j=0}^{\varphi_s-1}
\bigl(u(t_{j+1},X_{t_{j+1}},M_{t_j})-u(t_{j+1},X_{t_j},M_{t_j})\bigr)
\Delta_j M^{(i)}\\
&&\hphantom{\lim_{r \to\infty} \Biggl[}
{} + \sum_{j=0}^{\varphi_s-1}
\bigl(u(t_{j+1},X_{t_{j+1}},M_{t_{j+1}})-u(t_{j+1},X_{t_{j+1}},M_{t_j})\bigr)
\Delta_j M^{(i)}\Biggr].\nonumber
\end{eqnarray}
In addition, we can write
\begin{eqnarray*}
&&u(t_{j+1},X_{t_{j+1}},M_{t_j})-u(t_{j+1},X_{t_j},M_{t_j})\\
&&\qquad=\sum_{k=1}^n
\bigl(u\bigl(t_{j+1},X_{t_{j}}^{(1)},\ldots
,X_{t_{j}}^{(k-1)},X_{t_{j+1}}^{(k)},\ldots
,X_{t_{j+1}}^{(n)},M_{t_j}\bigr)\\
&&\qquad\quad\hphantom{\sum_{k=1}^n\bigl(}
{}-u\bigl(t_{j+1},X_{t_{j}}^{(1)},\ldots
,X_{t_{j}}^{(k-1)},X_{t_{j}}^{(k)},
\ldots,X_{t_{j+1}}^{(n)},M_{t_j}\bigr)\bigr).
\end{eqnarray*}
Each term of this sum can be written as
%
%e5.9 ###
%
\begin{eqnarray}
\label{bracketmdim1}
&&u\bigl(t_{j+1},X_{t_{j}}^{(1)},\ldots
,X_{t_{j}}^{(k-1)},X_{t_{j+1}}^{(k)},\ldots
,X_{t_{j+1}}^{(n)},M_{t_j}\bigr)\nonumber\\
&&\quad{}-u\bigl(t_{j+1},X_{t_{j}}^{(1)},\ldots
,X_{t_{j}}^{(k-1)},X_{t_{j}}^{(k)},
\ldots,X_{t_{j+1}}^{(n)},M_{t_j}\bigr)\\
&&\qquad=(\overline{\partial_2 u}) \bigl(\Delta_j X^{(1)},\ldots,\Delta_j
X^{(n)}\bigr)^\ast,\nonumber
\end{eqnarray}
where $\overline{\partial_2 u}
:=
(\partial_{1+k}u(t_{j+1},X_{t_j}^{(1)},\ldots,X_{t_j}^{(k-1)},\bar{X}_{t_
j}^{(k)},X_{t_{j+1}}^{(k+1)},\ldots,X_{t_{j+1}}^{(n)},M_{t_j})
)_{1\le
k\le n}$ and $\bar{X}_{t_j}^{(k)}$ is a suitable random point in the
interval $[X_{t_j}^{(k)} \wedge X_{t_{j+1}}^{(k)},X_{t_j}^{(k)}\vee
X_{t_{j+1}}^{(k)}]$.
Similarly, we obtain
%
%e5.10 ###
%
\begin{eqnarray}
\label{bracketmdim2}
&&
u(t_{j+1},X_{t_{j+1}},M_{t_{j+1}})-u(t_{j+1},X_{t_{j+1}},M_{t_j})\nonumber\\[-8pt]\\[-8pt]
&&\qquad=(\overline
{\partial_3 u}) \bigl(\Delta_j M^{(1)},\ldots,\Delta_j
M^{(d)}\bigr)^\ast\nonumber
\end{eqnarray}
with
\[
\overline{\partial_3u}:=\bigl(\partial
_{1+n+k}u(t_{j+1},X_{t_{j+1}},M_{t_j}^{(1)},\ldots
,M_{t_j}^{(k-1)},\bar{M}_{t_j}^{(k)},M_{t_{j+1}}^{(k+1)},\ldots
, M_{t_{j+1}}^{(d)})\bigr)_{1\le k\le d}.
\]
Combining relations (\ref{Sr2suite}),
(\ref{bracketmdim1}) and (\ref{bracketmdim2}) we deduce that
%
%e5.12 ###
%e5.11 ###
%
\begin{eqnarray}
\label{bracket2}\qquad\hspace*{8pt}
\lim_{r \to\infty} S_{r,2}^{(i)} &=&\lim_{r \to\infty} \sum
_{j=0}^{\varphi_s-1} \bigl[\bigl(\overline{\partial_2 u}\bigr) \bigl(\Delta_j
X^{(1)},\ldots,\Delta_j X^{(n)}\bigr)^\ast\Delta_j M^{(i)}\\
&&\hphantom{\lim_{r \to\infty} \sum_{j=0}^{\varphi_s-1}\bigl[}
{}   +(\overline{\partial_3 u}) \bigl(\Delta_j M^{(1)},\ldots
,\Delta_j M^{(d)}\bigr)^\ast\Delta_j M^{(i)}\bigr] \nonumber\\
& = &\lim_{r \to\infty} \sum_{j=0}^{\varphi_s-1} \bigl[\partial_2
u(t_j,X_{t_j},M_{t_j}) \bigl(\Delta_j X^{(1)},\ldots,\Delta_j
X^{(n)}\bigr)^\ast\Delta_j M^{(i)} \nonumber\\
&&\hphantom{\lim_{r \to\infty} \sum_{j=0}^{\varphi_s-1}\bigl[}
{}  + \partial_3 u(t_j,X_{t_j},M_{t_j}) \bigl(\Delta_j
M^{(1)},\ldots,\Delta_j M^{(d)}\bigr)^\ast\nonumber\\[-8pt]\\[-8pt]
&&\hspace*{142.51pt}
{} \times\Delta_j M^{(i)}+
R(i,j,r)\bigr],\nonumber
\end{eqnarray}
where $R(i,j,r)$ is defined as
\begin{eqnarray*}
R(i,j,r)&=&\bigl((\overline{\partial_2 u})-\partial_2
u(t_j,X_{t_j},M_{t_j})\bigr) \bigl(\Delta_j X^{(1)},\ldots,\Delta_j
X^{(n)}\bigr)^\ast\Delta_j M^{(i)}\\
&&{} + \bigl((\overline{\partial_3 u})-\partial_3
u(t_j,X_{t_j},M_{t_j})\bigr) \bigl(\Delta_j M^{(1)},\ldots,\Delta_j
M^{(d)}\bigr)^\ast\Delta_j M^{(i)}.
\end{eqnarray*}
Since
%
%e5.14 ###
%e5.13 ###
%
\begin{eqnarray}\qquad
\label{bracket4}
&&\lim_{r\to\infty} \sum_{j=0}^{\varphi_s-1} \bigl[\partial_2
u(t_j,X_{t_j},M_{t_j}) \bigl(\Delta_j X^{(1)},\ldots,\Delta_j
X^{(n)}\bigr)^\ast\Delta_j M^{(i)}\\
&&\hspace*{78pt}{}
+(\partial_3 u) \bigl(\Delta_j M^{(1)},\ldots,\Delta
_j M^{(d)}\bigr)^\ast\Delta_j M^{(i)}\bigr]\nonumber\\
&&\qquad= \biggl(\int_t^s [\partial_2 u(r,X_r,M_r)
\sigma(r,X_r,M_r)\nonumber\\[-8pt]\\[-8pt]
&&\qquad\quad\hspace*{68.33pt}{} +
\partial_3 u(r,X_r,M_r)] \,d\langle M, M
\rangle_r\biggr)^{(i)},\nonumber
\end{eqnarray}
relation (\ref{bracket1}) follows from equations (\ref{bracket2})
and (\ref{bracket4}) provided the following equation
holds:
%
%e5.15 ###
%
\begin{equation}
\label{bracket3} \lim_{r \to\infty} \Biggl\vert\sup_{t\leq s \leq
T} \sum_{j=0}^{\varphi_s-1}
R(i,j,r) \Biggr\vert=0.
\end{equation}
We conclude the proof by showing relation (\ref{bracket3}). Let
\[
A^{(r)}\,{:=}\,\sup_{\vert s_1-s_2\vert\leq\vert
\pi^{(r)}\vert,   a,b \in[s_1,s_2],   k=1,\ldots,n } \{\vert
\partial_{1+k}
u(s_2,X_a,M_{s_1})\,{-}\,\partial_{1+k} u(s_1,X_b,M_{s_1})\vert\}
\]
and
\begin{eqnarray*}
B^{(r)} &:=&\sup_{\vert s_1-s_2\vert\leq\vert
\pi^{(r)}\vert,   a,b \in[s_1,s_2],   k=1,\ldots,d}
\{\vert\partial_{1+n+k} u(s_2,X_{s_2},M_a)\\
&&\hspace*{132.8pt}{}-\partial_{1+n+k}
u(s_1,X_{s_2},M_b)\vert\}.
\end{eqnarray*}
For $1\le i\le d$, $r\in\mathbb{N}$ we have for any $s$ in $[t,T]$ that
\begin{eqnarray*}
\sum_{j=0}^{\varphi_s-1} \vert R(i,j,r) \vert
&\leq& c A^{(r)} \sum_{j=0}^{\varphi_s-1} \sum_{k=1}^n \bigl\vert\Delta
_j X^{(k)}
\Delta_j M^{(i)}\bigr\vert\\
&&{}+ c B^{(r)} \sum_{j=0}^{\varphi_s-1} \sum
_{k=1}^d \bigl\vert
\Delta_j M^{(k)} \Delta_j M^{(i)}\bigr\vert
\\
&\leq&\frac{c}{2} A^{(r)} \sum_{j=0}^{\varphi_s-1} \sum_{k=1}^n
\bigl[\bigl\vert\Delta_j
X^{(k)}\bigr\vert^2+ \bigl\vert\Delta_j M^{(i)}\bigr\vert^2\bigr]\\
&&{}+\frac{c}{2} B^{(r)} \sum_{j=0}^{\varphi_s-1} \sum_{k=1}^d
\bigl[\bigl\vert\Delta_j M^{(k)}\bigr\vert^2+ \bigl\vert\Delta_j M^{(i)}\bigr\vert^2\bigr]\\
&\leq&\frac{c}{2} A^{(r)} \sum_{k=1}^n \Biggl[\sum_{j_1=1}^{r-1}
\bigl\vert\Delta_{j_1} X^{(k)}\bigr\vert^2+ \sum_{j_2=1}^{r-1} \bigl\vert\Delta_{j_2}
M^{(i)}\bigr\vert^2\Biggr]\\
&&{}+\frac{c}{2} B^{(r)} \sum_{k=1}^d
\Biggl[\sum_{j_1=1}^{r-1} \bigl\vert\Delta_{j_1} M^{(k)}\bigr\vert^2+ \sum_{j_2=1}^{r-1}
\bigl\vert\Delta_{j_2} M^{(i)}\bigr\vert^2\Biggr].
\end{eqnarray*}
Thus,
\begin{eqnarray*}
&&\sup_{t\leq s \leq T} \sum_{j=0}^{\varphi_s-1} \vert R(i,j,r)\vert
\\
&&\qquad \leq\frac{c}{2} A^{(r)} \sum_{k=1}^n \Biggl[\sum_{j_1=1}^{r-1}
\bigl\vert\Delta_{j_1} X^{(k)}\bigr\vert^2+ \sum_{j_2=1}^{r-1} \bigl\vert\Delta_{j_2}
M^{(i)}\bigr\vert^2\Biggr]\\
&&\qquad\quad{} +\frac{c}{2} B^{(r)} \sum_{k=1}^d
\Biggl[\sum_{j_1=1}^{r-1} \bigl\vert\Delta_{j_1} M^{(k)}\bigr\vert^2+ \sum_{j_2=1}^{r-1}
\bigl\vert\Delta_{j_2} M^{(i)}\bigr\vert^2\Biggr].
\end{eqnarray*}
According to Proposition \ref{prop:diffuv}(iv), we have that
\[
\lim_{r \to\infty} A^{(r)}=\lim_{r \to\infty} B^{(r)}=0,\qquad
\PP\mbox{-a.s.}
\]
On the other hand,
\[
\cases{
\displaystyle \lim_{r \to\infty} \sum_{j=1}^{r-1} \bigl\vert\Delta_{j} X^{(k)}\bigr\vert
^2 = \bigl\langle X^{(k)}, X^{(k)} \bigr\rangle_z,\vspace*{2pt}\cr
\displaystyle \lim_{r \to\infty} \sum_{j=1}^{r-1} \bigl\vert\Delta_j
M^{(k)}\bigr\vert^2=\bigl\langle M^{(k)}, X^{(k)} \bigr\rangle_z,
}
\]
which concludes the proof.
\end{pf*}

As an immediate consequence of Proposition \ref{prop:diffuv} and of
Theorem \ref{theorem:main}, we get the following corollary.
\begin{corollary}
Assume that $M$ is a Markov process. Assume that~\textup{(H0)},
\textup{(H1)}--\textup{(H3)} are in force for the FBSDE (\ref{SDE}) and (\ref
{BSDE}). Then by Theorem \ref{theorem:Markovproperty}, there exists
a deterministic function $u$ such that
$Y_s^{t,x,m}=u(s,X_s^{t,x,m},M_s^{t,m}),\allowbreak  s\in[t,T]$. Assume in
addition that the assumption \textup{(MRP)} (see Section
\ref{section:diff2}) is in force, then for all $s \in[t,T]$ we have
$\nu$-a.e.
\[
Z_s^{t,x,m}=
\partial_2 u(s,X_s^{t,x,m},M_s^{t,m}) \sigma(s,X_s^{t,x,m},M_s^{t,m})
+ \partial_3 u(s,X_s^{t,x,m},M_s^{t,m}).
\]
\end{corollary}

%s6 ###
\section{Application to utility based pricing and hedging in
incomplete markets}\label{sec6}

In this section we study the exponential utility based
indifference
price approach for pricing and\vadjust{\goodbreak} hedging insurance related derivatives in
incomplete markets.
Thereby we will interpret relation (\ref{eq:main}) as a delta hedging formula.
Since in the Brownian setting it is shown in
\cite{AnkirchnerImkellerDosReis2} that this relation can be expressed
as a
function of the gradient of the indifference price and correlation
coefficients, we only sketch the arguments here. Let us explain how these
quantities translate into our local martingale framework with the more
complex Markovian structure. Consider an $n$-dimensional process describing
nontradable risk
\begin{eqnarray*}
R^{t,r,m}_s
&=&
r + \int_t^s \sigma(u,R^{t,r,m}_u, M^{t,m}_u) \,dM_u\\
&&{}+ \int_t^s b(u,R^{t,r,m}_u, M^{t,m}_u) \,dC_u ,
\qquad s \in[t,T],
\end{eqnarray*}
where $\sigma\in\real^{n \times d}$ and $b \in\real^{n \times1}$
are measurable functions.
An agent aims to price and hedge a derivative of the form $F(R^{t,r,m}_T)$,
with $F$ being a bounded measurable function.
The hedging instrument is a financial market consisting of $k$ risky assets
in units of the numeraire that evolve according to the following SDE:
\[
dS_s
=
S_s \bigl( \beta(s,R^{t,r,m}_s,M^{t,m}_s) \,dM_s + \alpha
(s,R^{t,r,m}_s,M^{t,m}_s)\,dC_s \bigr), \qquad  s \in[t,T],
\]
where the measurable processes $\alpha$ and $\beta$ take their values
in $
\real^{k \times1}$, respectively, in $\real^{k \times d}$. Observe
that the price
processes of tradable assets $S$ are linked to the risk process via the
martingale $M$, its quadratic variation and the functions $\beta$ and
$\sigma$. In addition, we assume $k \leq d$ in order to exclude arbitrage
opportunities. The small agent's preferences are represented through the
exponential utility function with risk aversion coefficient $\kappa
>0$, that is,\looseness=1
\[
U(x)=-e^{-\kappa x},\qquad   x \in\real.
\]\looseness=0
The agent wants to maximize his expected utility by trading in the market.
His value function is given by
\[
V^F(x,t,r,m) = \sup_{\lambda} \E\Biggl[ U\Biggl(x + \sum_{i=1}^k \int_t^T
\lambda^{(i)}_s \,\frac{dS^{(i)}_s}{S_s^{(i)}} + F(R_T^{t,r,m})\Biggr) \Biggr],
\]
where $x$ is his initial capital and $\lambda^{(i)}$ denotes the momentary
value of his portfolio fraction invested in the $i$th asset. This
optimization problem can be reduced to solving a quadratic BSDE whose
generator has been given in~\cite{HuImkellerMuller} for the Brownian case
and then extended to our setting in \cite{Morlais1}. A~way to price and
hedge the derivative $F(R_T^{r,t,m})$ is to consider the indifference price
$p(t,r,m)$ defined via $V^F(x-p(t,r,m),t,r,m)=V^0(x,t,r,m)$. According to~\cite{AnkirchnerImkellerDosReis2}, the indifference price can\vspace*{1pt} be expressed
as $p(t,r,m)=Y^{F,t,r,m}-Y^{0,t,r,m}$, where\vadjust{\goodbreak}
$(Y^{F,t,r,m},Z^{F,t,r,m},L^{F,t,r,m})$ is the solution of the BSDE
%
%e6.1 ###
%
\begin{eqnarray}
\label{DeltaBSDE}\qquad
Y^{F,t,r,m}_s
&=&F(R^{t,r,m})-\int_s^T Z^{F,t,r,m}_u \,dM_u\nonumber\\
&&{} + \int_s^T
f(u,R^{t,r,m}_u,M^{t,m}_u,Z_u^{F,t,r,m}q_u^*) \,dC_u \\
&&{} -\int_s^T dL^{F,t,r,m}_u + \frac{\kappa}{2} \int_s^T d\langle
L^{F,t,r,m},L^{F,t,r,m}\rangle_u, \qquad t\in[0,T].\nonumber
\end{eqnarray}
Here the generator $f$ is obtained explicitly through the martingale
optimality principle; cf. \cite{HuImkellerMuller,Morlais1} and possesses
properties covered by the hypotheses of Theorem~\ref{theorem:main}. To
implement utility indifference, we have to describe the optimal strategies
$\hat\lambda^F$ and $\hat\lambda^0$. In \cite{HuImkellerMuller} it
is shown
that $\hat\lambda^F \beta(\cdot,R^{t,r,m},M^{t,m})$ [and $\hat
\lambda^0
\beta(\cdot,R^{t,r,m},M^{t,m})$] are given by the projection of a linear
function of $Z^{F,t,r}q^*$ (resp., $Z^{0,t,r}q^*$) on the constraint
set. Since $R^{t,r,m}$ is not tradable directly, $\beta$ plays the
role of
a filter for trading in the market. Due to\vspace*{1pt}
\cite{AnkirchnerImkellerDosReis2}, the optimal strategy to hedge
$F(R_T^{t,r,m})$ can be decomposed into a pure trading part $\hat
\lambda^0$
and the optimal hedge~$\Delta$, which is the part of the strategy that
replicates the derivative $F(R_T^{t,r,m})$. Using the Markov property\vspace*{1pt} given
in Theorem \ref{theorem:Markovproperty}, we see that there exists a
deterministic function $u^F$ such that
$Y^{F,t,r,m}=u^F(\cdot,R^{t,r,m},M^{t,m})$. Moreover, the projection
mentioned above can be explicitly expressed. Indeed from
\cite{AnkirchnerImkellerDosReis2}, proof of Theorems 4.2 and 4.4, we have
\[
\hat\lambda^F_s - \hat\lambda^0_s =
(Z^{F,t,r,m}_s-Z^{0,t,r,m}_s) q_s^* \beta^*(\beta\beta^*)^{-1}\beta
(s,R^{t,r,m}_s,M^{t,m}_s), \qquad  s\in[t,T].
\]
This leads to
\begin{eqnarray*}
\Delta(t,r,m) &=& (\hat\lambda^F - \hat\lambda^0 )\beta^*(\beta
\beta^*)^{-1}(t,r,m)\\
&=& (Z^{F,t,r,m}_t-Z^{0,t,r,m}_t) q^*_t \beta
^*(\beta\beta^*)^{-1}(t,r,m).
\end{eqnarray*}
Using formula (\ref{eqn3}), we derive
%
%e6.2 ###
%
\begin{equation}
\label{Delta}
\Delta(t,r,m)=[\partial_2 p(t,r,m) \sigma(t,r,m) + \partial_3
p(t,r,m) ] q^*_t \beta^*(\beta\beta^*)^{-1}(t,r,m).\hspace*{-28pt}
\end{equation}
We emphasize that, as a consequence of the particular form of the
driver $f$ in (\ref{DeltaBSDE}),
if $M$ has independent increments and the coefficients $\sigma$, $b$,
$\beta$ and~$\alpha$ do not depend on $M$ [see Remarks \ref
{indIncbis}(ii) and (iii)], then relation
(\ref{Delta}) is replaced by
%
%e6.3 ###
%
\begin{equation}
\label{Delta2}
\Delta(t,r)=[\partial_2 p(t,r) \sigma(t,r) ] q^*_t \beta
^*(\beta\beta^*)^{-1}(t,r).
\end{equation}
Finally, note that we obtain formulae (\ref{eqn3}) and (\ref
{Delta}) under condition
(MRP) (see Section~\ref{section:diff2}).
However, we believe that this condition is not necessary for
deriving~(\ref{Delta}).
Finally, we mention that in \cite{FreiSchweizer} the authors also
represent the indifference price
as the difference of two $Y$ processes solution to a BSDE when the
price process is generated by a general semimartingale.
However, the authors do not prove a representation formula for the $Z$
process of their BSDE but rather obtain some regularity property of $(Z,L)$,
that is, under some condition on the claim $F(R_T^{t,x,m})$ they prove
that $Z \cdot d$ and $L$ are BMO martingales for the minimal entropy
martingale measure.
Thus, the authors do not obtain a representation of the form (\ref
{Delta2}) for the delta hedge.

\section*{Concluding remarks}

In this paper we prove the representation formu\-la~(\ref{eqn3}) for
the control process of a quadratic growth BSDE driven by a~continuous local martingale. This can be used for giving an explicit
representation of the delta hedge in utility indifference based
hedging of insurance derivatives with exponential preferences. We also
provide the Markov property and differentiability
of the FBSDE (\ref{SDE})  and (\ref{BSDE}) in the initial state parameter
of its forward part. This last property is obtained under an
additional assumption (MRP).
However, we think that differentiability should hold without this
assumption and that
different techniques have to be developed for achieving this goal.

Additionally, as already mentioned in this paper, Malliavin's
calculus has been used by several authors to recover formula
(\ref{eq:main}) in the Brownian framework. Our alternative method is valid
in this setting and seems to present advantages in some practical
situations. Actually, Malliavin's calculus is known for its efficiency in
several topics, however, it also usually requires more regularity than the
problem needs intrinsically. In \cite{AnkirchnerHeyne}, the authors study
the quadratic hedging problem of contingent claims with basis risk when the
hedging instrument and the underlying of the contingent are related via a~random correlation process. As given in \cite{AnkirchnerHeyne}, the hedging
strategy is described via a representation formula of the form
(\ref{eq:main}) for the control process of the backward part of a FBSDE
driven by a Brownian motion. In this case, the coefficient of the forward
process depends on a correlation process~$\rho$ which is itself
solution of
a Brownian SDE. As explained in \cite{AnkirchnerHeyne}, a Section~3.4
comment, the use of Malliavin's calculus enforces that the
derivatives of the coefficients of the SDE defining $\rho$ have bounded
derivatives. This additional regularity is not necessary in our approach
and would allow one to consider more examples of correlation processes with
only locally Lipschitz bounded derivatives.

\begin{appendix}\label{section:appendix}
%s7 ###
\section*{Appendix}

In the first section of this Appendix we provide the transformation of
a~BSDE of the form
(\ref{BSDE}) which is needed in Section
\ref{section:differentiability} and give a priori estimates on the
solution of the transformed BSDE with respect to its terminal
condition and its generator. Then in Appendix \ref{Appendix2}, we
present an explicit example of the situation described in
Proposition \ref{CinlaretalProtter}(ii).

%s7.1 ###
\subsection{\texorpdfstring{Transformation of the BSDE (\protect\ref{BSDE}) under (MRP)}
{Transformation of the BSDE (2.3) under (MRP)}}
\label{Appendix1}

We start giving a justification that under (MRP) the BSDE of the form
%
%e7.1 ###
%
\setcounter{equation}{0}
\begin{eqnarray}
\label{BSDEbracket}
Y_t&=&B-\int_t^T Z_s \,dM_s +\int_t^T f(s,Y_s,Z_s q_s^*) \,dC_s\nonumber\\[-10pt]\\[-10pt]
&&{} - \int_t^T
dL_s + \frac{\kappa}{2}\int_t^T d\langle L,L \rangle_s\nonumber
\end{eqnarray}
can be transformed into a BSDE of the form
%
%e7.2 ###
%
\begin{equation}
\label{BSDEtilde} Y_t=B-\int_t^T \tilde Z_s \,d\tilde M_s + \int_t^T
h(s,Y_s,\tilde Z_s \tilde q_s^*) \,d\tilde{C}_s,
\end{equation}
where for all $s\in[0,T]$
\begin{eqnarray*}
\tilde M_s &:=& \pmatrix{ M_s \cr N_s },\qquad
 \tilde q_s:= \pmatrix{ q_s \sqrt{\varphi_1(s)} & 0 \cr0 &
\sqrt{\varphi_2(s)} },\\[-2pt]
\tilde{C}_s&:=&\arctan\Biggl(\sum_{i=1}^d \bigl\langle M^{(i)}, M^{(i)}
\bigr\rangle_s+\langle N, N \rangle_s\Biggr),
\end{eqnarray*}
$\tilde{Z}_s := ( Z_s,   U_s )$, with $\varphi_1$ and $\varphi_2$
denoting two nonnegative positive predictable processes defined below. Let
\[
d\mu_s^1:= \frac{\sum_{i=1}^d d\langle
M^{(i)}, M^{(i)} \rangle_s}{1+(\sum_{i=1}^d \langle M^{(i)},
M^{(i)} \rangle_s +
\langle N, N \rangle_s)^2}
\]
and
\[
d\mu_s^2:= \frac{d \langle N, N
\rangle_s}{1+(\sum_{i=1}^d \langle M^{(i)}, M^{(i)} \rangle_s +
\langle N,
N \rangle_s)^2}.
\]
For every $\omega$ in $\Omega$, the measure $d\mu_t^1(\omega)$
[resp., $d\mu_t^2(\omega)$] is absolutely continuous with
respect to $d(\mu_t^1+\mu_t^2)(\omega)$. Hence, since $\mu^1$ and
$\mu^2$ are predictable processes~\cite{DellacherieMeyer}, Theorem
VI.68 and its remark
imply that there exist two predictable processes $\varphi_1$ and
$\varphi_2$ such that
\begin{eqnarray*}
\mu_t^1&=&\int_0^t \varphi_1(s) \,d(\mu^1+\mu^2)(s),\\[-2pt]
\mu_t^2&=&\int_0^t \varphi_2(s) \,d(\mu^1+\mu^2)(s)\qquad \forall t \in[0,T].
\end{eqnarray*}
In addition, we have that $0 \leq\varphi_i(s) \leq1$ for all $s$ in
$[0,T]$ $\PP$-a.s. for $i=1,2$.
Indeed, because $\varphi_i$, $i=1,2$ is a density, it is nonnegative
and from $d(\mu^1+\mu^2)(s)=(\varphi_1+\varphi_2)(s)d(\mu^1 + \mu
^2)(s)$ it follows $(\varphi_1+\varphi_2)(s)=1$, $d(\mu^1 +
\mu^2)(s)$-a.e.

Recall that
\[
dC_s=\frac{\sum_{i=1}^d d\langle M^{(i)}, M^{(i)} \rangle_s}{1+
(\sum_{i=1}^d \langle M^{(i)}, M^{(i)} \rangle_s)^2}.\vadjust{\goodbreak}
\]
We have for $t\in[0,T]$
\begin{eqnarray*}
&&\int_t^T f(s,Y_s,Z_s q_s^*) \,d C_s + \frac{\kappa}{2} \int_t^T d
\langle L, L \rangle_s\\
&&\qquad=\int_t^T f(s,Y_s,Z_s q_s^*) \,d C_s + \frac{\kappa}{2} \int_t^T
U_s^2 \,d \langle N, N \rangle_s\\
&&\qquad=\int_t^T \tilde{f}(s,Y_s,Z_s (\tilde{q}_s^*)_{1,1}) \frac{\sum
_{i=1}^d d\langle
M^{(i)}, M^{(i)} \rangle_s}{1+(\sum_{i=1}^d \langle M^{(i)},
M^{(i)} \rangle_s +
\langle N, N \rangle_s)^2} \\
&&\qquad\quad{} + \int_t^T g(s,U_s) \frac{d \langle N, N
\rangle_s}{1+(\sum_{i=1}^d \langle M^{(i)}, M^{(i)} \rangle_s +
\langle N,
N \rangle_s)^2},
\end{eqnarray*}
where for $s\in[t,T]$
\[
\tilde{f}(s,y,z):=
\cases{
\displaystyle f(s,y,z   \varphi_1(s)^{-1/2}) \times\frac{1+(\sum_{i=1}^d
\langle M^{(i)}, M^{(i)} \rangle_s +
\langle N, N \rangle_s)^2}{1+(\sum_{i=1}^d \langle M^{(i)},
M^{(i)} \rangle_s )^2}, \vspace*{2pt}\cr
\qquad\mbox{if } \varphi_1(s) \neq0,\vspace*{2pt}\cr
\displaystyle f(s,y,0) \times\frac{1+(\sum_{i=1}^d \langle M^{(i)}, M^{(i)}
\rangle_s +
\langle N, N \rangle_s)^2}{1+(\sum_{i=1}^d \langle M^{(i)},
M^{(i)} \rangle_s )^2}, \vspace*{2pt}\cr
\qquad\mbox{if } \varphi_1(s) = 0,}
\]
and
\[
g(s,u):=\frac{\kappa}{2} u^2 \Biggl(1+\Biggl(\sum_{i=1}^d \bigl\langle M^{(i)},
M^{(i)} \bigr\rangle_s + \langle N, N \rangle_s\Biggr)^2\Biggr).
\]
With this definition we have that $ f(s,Y_s,Z_s q_s^*)=\tilde
{f}(s,Y_s,(\tilde{Z}_s \tilde{q}_s^*)_1)$.
Hence,
\begin{eqnarray*}
&&\int_t^T f(s,Y_s,Z_s q_s^*) \,d C_s + \frac{\kappa}{2} \int_t^T d
\langle L, L \rangle_s\\
&&\qquad=\int_t^T \bigl(\tilde{f}(s,Y_s,Z_s (\tilde{q}_s^*)_{1,1}) \varphi
_1(s) + g(s,U_s) \varphi_2(s)\bigr)\\
&&\qquad\quad\hphantom{\int_t^T}
{}\times \frac{\sum_{i=1}^d d\langle
M^{(i)}, M^{(i)} \rangle_s + d\langle N, N \rangle_s}{1+(\sum
_{i=1}^d \langle M^{(i)}, M^{(i)} \rangle_s + \langle N, N \rangle
_s)^2}\\
&&\qquad=\int_t^T \bigl(\tilde{f}(s,Y_s,Z_s (\tilde{q}_s^*)_{1,1}) \varphi
_1(s) + g(s,U_s) \varphi_2(s)\bigr) \,d \tilde{C}_s.
\end{eqnarray*}
As a consequence, letting
\begin{eqnarray*}
h(s,Y_s,\tilde{Z}_s \tilde{q}_s^*):\!&=&\tilde{f}(s,Y_s,(\tilde{Z}_s
\tilde{q}_s^*)_1)
\varphi_1(s) + g(s,(\tilde{Z}_s \tilde{q}_s^*)_2)\\
&=&\tilde{f}(s,Y_s,\tilde{Z}_s (\tilde{q}_s^*)_{1,1})
\varphi_1(s) + g(s,(\tilde{Z}_s \tilde{q}_s^*)_2),
\end{eqnarray*}
we obtain that
(\ref{BSDEbracket}) can be written as
\[
Y_t=B-\int_t^T Z_s \,dM_s - \int_t^T U_s \,dN_s + \int_t^T
h(s,Y_s,\tilde{Z}_s \tilde{q}_s^*) \,d\tilde{C}_s
\]
and if the initial generator $f$ satisfies the hypothesis (H3), so does
the generator $h$ since $\varphi_1$, $\varphi_2$, $\langle M^{(i)},
M^{(j)} \rangle_T$ and $\langle N, N\rangle_T$ are bounded processes
for all~$i, j$ in $\{1,\ldots,d\}$. In particular,
$h$ preserves the growth in the variables~$y,z,u$. Hence, we derive
at BSDE (\ref{BSDEtilde}).

%s7.2 ###
\subsection{A priori estimates}

Now we assume that $M$ itself satisfies the martingale representation
theorem and we consider the following BSDE:
%
%e7.3 ###
%
\begin{equation}
\label{BSDE:apriori}
Y_t=B-\int_t^T Z_s \,dM_s + \int_t^T f(s,Y_s,Z_s q_s^*) \,dC_s,
\end{equation}
where $M,q,C$ are defined as in Section \ref{prel}.
Suppose that the terminal condition $B$ is a bounded real-valued random
variable, the generator $f$ satisfies assumption (H4) and that $(Y,Z)$
is a solution to (\ref{BSDE:apriori}).
The following a priori inequality is crucial for our differentiability
and representation results.
\begin{lemma}
\label{appendix:apriori} We assume that for every $\beta\geq1$ we
have $\int_0^T \vert f(s,0,0) \vert \,dC_s \in L^\beta(\PP)$.
Let $p>1$, then there exist constants $q \in(1, \infty)$, $c>0$
depending only on $T$, $p$ and on the BMO-norm of $K \cdot M$ such that
\begin{eqnarray*}
&&\E\Bigl[\sup_{t \in[0,T]} \vert Y_t\vert^{2 p}\Bigr] + \E\biggl[
\biggl(\int_0^T \vert q_s Z_s^* \vert^2 \,dC_s \biggr)^p\biggr]
\\
&&\qquad\leq c \E\biggl[\vert B \vert^{2 pq} + \biggl( \int_0^T \vert
f(s,0,0) \vert \,dC_s \biggr)^{2pq}\biggr]^{{1/q}}.
\end{eqnarray*}
\end{lemma}
\begin{pf}
We follow \cite{BriandConfortola}, Lemmas 7, 8 and Corollary 9 (see
also \cite{AnkirchnerImkellerDosReis2}, Lem\-ma~6.1) which have been
designed for the Brownian setting.
However, as we will show below, their arguments can be extended to the
case of continuous local martingales.
We proceed in several steps.

In a first step we exploit properties of BMO martingales.
Let
\[
J_s=\cases{
\displaystyle \frac{f(s,Y_s,Z_s q_s^*)-f(s,0,Z_s q_s^*)}{Y_s}, &\quad if $Y_s \neq0$,
\vspace*{2pt}\cr
0, &\quad otherwise,}
\]
and
\[
H_s=\cases{
\displaystyle \frac{f(s,0,Z_s q_s^*)-f(s,0,0)}{|q_s Z_s^*|^2}Z_s, &\quad if $|q_s
Z_s^*|^2 \neq0$,
\cr
0, &\quad otherwise.}\vspace*{-2pt}
\]
Then BSDE (\ref{BSDE:apriori}) has the form
%
%e7.4 ###
%
\begin{eqnarray}
\label{finparts}\qquad
Y_t&=&B-\int_t^T Z_s \,dM_s\nonumber\\[-11pt]\\[-11pt]
&&{} +\int_t^T \bigl(J_s Y_s +(q_s H_s^*)(q_sZ_s^*)^*+
f(s,0,0)\bigr)\,dC_s,\qquad  t\in[0,T].\nonumber\vspace*{-2pt}
\end{eqnarray}
Due to (H4) we have $|q H^*| \leq|q K^*|$ and it follows that $H \cdot
M$ is a BMO($\PP$) martingale.
Furthermore, we know from \cite{Kazamaki}, Theorem 3.1, that there
exists a \mbox{$\hat q >1$} such that the reverse H\"{o}lder inequality holds,
that is, there exists a~constant $c>0$ such that
%
%e7.5 ###
%
\begin{equation}
\label{revholder}
\mathcal E(H \cdot M)_t^{- \hat q} \E[ \mathcal E(H \cdot
M)_T^{\hat q} | \mathcal F_t ] \leq c.\vspace*{-2pt}
\end{equation}
By \cite{Kazamaki}, Theorem 2.3, the measure $\Q$ defined by $d\Q
=\mathcal E(H \cdot M)_T \,d\PP$ is a~probability measure.
Girsanov's theorem implies that
\[
\Lambda=Z \cdot M - \int_0^\cdot(q_s H_s^*)(q_sZ_s^*)^* \,dC_s\vspace*{-2pt}
\]
is a local $\Q$-martingale.
This means that there exists an increasing sequence of stopping times
$(\tau^n)_{n\in\mathbb{N}}$ converging to $T$ such that $\Lambda
_{\cdot\wedge\tau^n}$ is a $\Q$-martingale for any $n\in\mathbb{N}$.
Letting $e_t=\exp(2\int_0^t |q_s K_s^*|^{2 \alpha} \,dC_s)$, $t \in
[0,T]$, with It\^{o}'s formula applied to $e_t Y_t^2$ we have
\begin{eqnarray*}
d (e_t Y_t^2) & = & 2 |q_t K_t^*|^{2 \alpha} e_{t} Y_{t}^2 \,dC_t + 2 e_t
Y_t \,dY_t + e_t |q_t Z_t^*|^2 \,dC_t \\[-2pt]
& = & 2 |q_tK_t^*|^{2 \alpha} e_t Y_t^2 \,dC_t + 2 e_t Y_t \,d \Lambda_t -
2 e_t Y_t^2 J_t \,dC_t\\[-2pt]
&&{} - 2 e_t Y_t f(t,0,0) \,dC_t+ e_t |q_t Z_t^*|^2 \,dC_t,\vspace*{-2pt}
\end{eqnarray*}
where we used (\ref{finparts}).
With the inequality $J_t \leq|q_t K_t^*|^{2 \alpha}$, $t \in[0,T]$,
which follows from assumption (H4), we know for $t\in[0,\tau^n]$
\begin{eqnarray*}
e_t Y_t^2
&\leq&
e_{\tau^n} Y_{\tau^n}^2
- \int_t^{\tau^n} 2 e_t Y_t \,d \Lambda_t
+ \int_t^{\tau^n} 2 e_t Y_t f(t,0,0) \,dC_t\\[-4pt]
&&{}- \int_t^{\tau^n} e_t |q_t Z_t^*|^2 \,dC_t.\vspace*{-2pt}
\end{eqnarray*}
Note that $e_t \geq1$ for all $t \in[0,T]$ and hence,
%
%e7.6 ###
%
\begin{eqnarray}
\label{basis}
e_t Y_t^2 + \int_t^{\tau^n} |q_s Z_s^*|^2 \,dC_s
&\leq&
e_{\tau^n} Y_{\tau^n}^2
- \int_t^{\tau^n} 2 e_s Y_s \,d \Lambda_s\nonumber\\[-11pt]\\[-11pt]
&&{} + \int_t^{\tau^n} 2 e_s Y_s f(s,0,0) \,dC_s.\nonumber\vspace*{-2pt}
\end{eqnarray}

In a second step we provide an estimate for $Y$.
We want to take the conditional expectation under the new measure $\Q$
in the previous inequality.
Therefore,\vadjust{\goodbreak} we need to check the integrability of the involved terms.
Observe that
%
%e7.7 ###
%
\begin{equation}
\label{ebound}
e_t \leq\exp\biggl( 2 \int_0^T |q_s K_s^*|^{2\alpha}\,dC_s \biggr),\qquad
t \in[0,T].
\end{equation}
Using successively the monotone convergence theorem and Jensen's
inequality, we derive for $p>1$
\[
\E\biggl[ \exp\biggl(p \int_0^T |q_s K_s^*|^{2 \alpha} \,dC_s\biggr)\biggr]
\leq C_T \sum_{n \geq0} \frac{p^n}{n!} \E\biggl[ \biggl(\int_0^T |q_s
K_s^*|^{2} \,dC_s\biggr)^{n \alpha}\biggr].
\]
The H\"{o}lder inequality again along with inequality \cite{Kazamaki},
page 26, gives
%
%e7.8 ###
%
\begin{eqnarray}\label{einsp}\quad
\E\biggl[ \exp\biggl(p \int_0^T |q_s K_s^*|^{2 \alpha} \,dC_s\biggr)\biggr]
& \leq & c \sum_{n \geq0} \frac{p^n}{n!} \E\biggl[ \biggl(\int_0^T |q_s
K_s^*|^{2} \,dC_s\biggr)^{n }\biggr]^{\alpha}
\nonumber\\[-8pt]\\[-8pt]
& \leq & c \sum_{n \geq0} \frac{(p \|K \cdot M\|_{\mathrm{BMO}_2}^{2\alpha
})^n}{n!^{1-\alpha}} < \infty.\nonumber
\end{eqnarray}
Thus, the process $e$ belongs to $\mathcal S^p(\PP)$ for all $p \geq1$
and using the H\"{o}lder inequality and formula (\ref{revholder}) we
see that $e_{\tau^n} Y_{\tau^n}^2 $, $e_T |B|^2$ and $\int_0^{T} 2
e_t |Y_t| |f(t,0,\allowbreak 0)| \,dC_t$ is in $\mathbb L^p(\Q)$ for all $p \geq1$.
In the same way we get the integrability of $\int_0^{\tau^n} 2 e_s
|Y_s| \,d \Lambda_s$.
Hence, we are allowed to take the conditional expectation in (\ref
{basis}) on both sides:
\[
e_t Y_t^2
\leq\E^\Q\biggl[ e_{\tau_n} Y_{\tau_n}^2 +\int_0^{T} 2 e_s |Y_s|
|f(s,0,0)|\,dC_s \Big| \mathcal F_t \biggr], \qquad  t \leq\tau_n.
\]
Now we let $n$ tend to infinity
\begin{eqnarray*}
e_t Y_t^{2} & \leq & \lim_{n \rightarrow\infty} \E^\Q\biggl[ e_{\tau
_n} Y_{\tau_n}^2 +\int_0^{T} 2 e_s |Y_s| |f(s,0,0)|\,dC_s \Big| \mathcal
F_t \biggr] \\
& \leq & \E^\Q\biggl[ e_{T} |B|^2 +\int_0^{T} 2 e_s
|Y_s||f(s,0,0)|\,dC_s \Big| \mathcal F_t \biggr],
\end{eqnarray*}
where we may apply the dominated convergence theorem because of (\ref{ebound}).
The Young inequality with a constant $c_1>0$ gives
\begin{eqnarray*}
Y_t^{2}
& \leq & \E^\Q\biggl[ \frac{e_{T}}{e_t} |B|^2 + 2 \int_0^{T} \frac
{e_s}{e_t} |Y_s||f(s,0,0)|\,dC_s \Big| \mathcal F_t \biggr]
\\[-1pt]
& \leq &
\E^\Q\biggl[e_{T}|B|^2 + \frac{1}{c_1} \sup_{t \in[0,T]} |Y_t|^2 +
c_1 e_T^2 \biggl(\int_0^{T} |f(s,0,0)|\,dC_s\biggr)^2 \Big| \mathcal F_t \biggr]
\\[-1pt]
& \leq & \E^\Q\biggl[ \frac{1}{c_1} \sup_{t \in[0,T]} |Y_t|^2 +
e^2_{T} \Theta_T \Big| \mathcal F_t \biggr],
\end{eqnarray*}
where we set $\Theta_T=|B|^2 +2c_1 (\int_0^{T} |f(s,0,0)|\,dC_s)^2$ and
we take into account that $e_s/e_t \leq e_T$ for all $s,t \in[0,T]$
and $e_T \leq e_T^2$.
Let $p>1$, then we have
\[
\sup_{t \in[0,T]} |Y_t|^{2p}
\leq
\sup_{t \in[0,T]} \E^{\Q} \biggl[ \frac{1}{c_1} \sup_{t \in
[0,T]} |Y_t|^2 + e^2_{T} \Theta_T \Big| \mathcal F_t \biggr]^{p}.
\]
We apply Doob's inequality to obtain
\begin{eqnarray*}
\E^{\Q} \Bigl[ \sup_{t \in[0,T]} |Y_t|^{2p} \Bigr]
& \leq & c \E^{\Q} \biggl[ \biggl( \E\biggl[ \frac{1}{c_1} \sup_{t
\in[0,T]} |Y_t|^2 + e^2_{T} \Theta_T \Big| \mathcal F_T \biggr]
\biggr)^{p} \biggr]
\\[-1pt]
& \leq & c \E^{\Q} \biggl[\frac{1}{c_1^p} \sup_{t \in[0,T]}
|Y_t|^{2p} + e_{T}^{2p} \Theta_T^p \biggr],
\end{eqnarray*}
and choosing $c_1$ such that $c/c_1^p <1$, we have
%
%e7.9 ###
%
\begin{equation}
\label{yestimate1}
\E^{\Q} \Bigl[ \sup_{t \in[0,T]} |Y_t|^{2p} \Bigr]
\leq c \E^{\Q} [e_{T}^{2p} \Theta_T^p ].
\end{equation}

In Step 3 we give an estimate on $Z$ under the measure $\Q$.
For $p > 1$ we deduce from (\ref{basis})
\begin{eqnarray*}
&&\biggl( \int_0^{\tau^n} |q_s Z_s^*|^2 \,dC_s \biggr)^p\\[-1pt]
&&\qquad\leq c \biggl(
|e_{\tau^n} Y_{\tau^n}^2|^p
+ \biggl| \int_0^{\tau^n} e_s Y_s \,d \Lambda_s \biggr|^p
+ \biggl(\int_0^{T} 2 e_s |Y_s| |f(s,0,0)| \,dC_s \biggr)^p \biggr).
\end{eqnarray*}
Then the Burkholder--Davis--Gundy and two times Young inequality (with
constants $\tilde{c}_1, \tilde{c}_2>0$) imply
\begin{eqnarray*}
\hspace*{-3pt}&& \E^{\Q} \biggl[\biggl( \int_0^{\tau^n} |q_s Z_s^*|^2 \,dC_s
\biggr)^p\biggr]
\\[-1pt]
\hspace*{-3pt}&&\qquad \leq
c \biggl(\E^{\Q} \Bigl[e_{T}^p \sup_{t\in[0,T]} |Y_{t}|^{2p}
\Bigr]
+ \E^{\Q} \biggl[ \biggl( \int_0^{\tau^n} e_s^2 Y_s^2 \vert q_s
Z_s^*\vert^2 \,dC_s \biggr)^{{p/2}} \biggr]\\
\hspace*{-3pt}&&\qquad\quad\hspace*{90pt}{} + \E^{\Q} \biggl[ \biggl(\int_0^{T} 2 e_s |Y_s|
|f(s,0,0)| \,dC_s \biggr)^p \biggr] \biggr)
\\
\hspace*{-3pt}&&\qquad \leq
c \biggl(\E^{\Q} \Bigl[e_{T}^p \sup_{t\in[0,T]} |Y_{t}|^{2p}
\Bigr]\\
\hspace*{-3pt}&&\qquad\quad\hphantom{c \biggl(}
{}
+ \E^{\Q} \Bigl[ (\tilde{c}_1+\tilde{c}_2) e_T^{2p} \sup_{t\in
[0,T]} |Y_t|^{2p}\Bigr]
+ \E^{\Q} \biggl[\frac{1}{\tilde{c}_1} \biggl( \int_0^{\tau^n}
\vert q_s Z_s^*\vert^2 \,dC_s \biggr)^{p} \biggr]\\
&&\qquad\quad\hspace*{143.2pt}{}+  \E^{\Q} \biggl[ \frac{1}{\tilde{c}_2} \biggl(\int_0^{T}
|f(s,0,0)| \,dC_s \biggr)^{2p} \biggr]\biggr),
\end{eqnarray*}
and because $e_T^p \leq e_T^{2p}$ and Fatou's lemma we have
\begin{eqnarray*}
&&\E^{\Q} \biggl[\biggl( \int_0^{T} |q_s Z_s^*|^2 \,dC_s
\biggr)^p\biggr]\\
&&\qquad
\leq
c \biggl(\E^{\Q} \Bigl[e_{T}^{2p} \sup_{t\in[0,T]} |Y_{t}|^{2p}
\Bigr] + \E^{\Q} \biggl[ \biggl(\int_0^{T} |f(s,0,0)| \,dC_s
\biggr)^{2p} \biggr]\biggr).
\end{eqnarray*}
We use the H\"{o}lder inequality with $r \geq1$, the estimate (\ref
{yestimate1}) and the H\"{o}lder inequality with $k \geq1$ again to deduce
%
%e7.10 ###
%
\begin{eqnarray}\label{zestimate}
&& \E^{\Q} \biggl[\biggl( \int_0^{T} |q_s Z_s^*|^2 \,dC_s
\biggr)^p\biggr]
\nonumber\\[1.2pt]
&&\qquad \leq
c \biggl(\E^{\Q} \bigl[e_{T}^{{2pr}/({r-1})} \bigr]^{({r-1})/{r}}
\E^{\Q} \Bigl[ \sup_{t\in[0,T]} |Y_{t}|^{2pr} \Bigr]^{
{1/r}}\nonumber\\[1.2pt]
&&\qquad\quad\hspace*{55.8pt}{} + \E^{\Q} \biggl[ \biggl(\int_0^{T} |f(s,0,0)| \,dC_s
\biggr)^{2p} \biggr] \biggr)\nonumber\\[-7.5pt]\\[-7.5pt]
&&\qquad \leq
c \biggl(\E^{\Q} \bigl[ e_{T}^{{2prk}/({k-1})} \bigr]^{
({k-1})/({rk})} \E^{\Q} [ \Theta^{prk} ]^{{1}/({rk})} \nonumber\\[1.2pt]
&&\qquad\quad\hspace*{42.1pt}{}+
\E^{\Q} \biggl[ \biggl(\int_0^{T} |f(s,0,0)| \,dC_s \biggr)^{2p}
\biggr] \biggr)\nonumber\\[1.2pt]
&&\qquad \leq
c \E^{\Q} \biggl[ |B|^{2prk} + \biggl(\int_0^{T} |f(s,0,0)| \,dC_s
\biggr)^{2prk} \biggr]^{{1}/({rk})}.\nonumber
\end{eqnarray}
Here we applied (\ref{einsp}) and in the last inequality we employ the
H\"{o}lder inequality with exponent $rk$ to the second summand in order
to obtain the last estimate.
We utilize the H\"{o}lder inequality with $rk$ to (\ref{yestimate1})
and hence, have
%
%e7.11 ###
%
\begin{equation}
\label{yestimate2}
\E^{\Q} \Bigl[ \sup_{t \in[0,T]} |Y_t|^{2p} \Bigr]
\leq
c \E^{\Q} \biggl[ |B|^{2prk} + \biggl(\int_0^{T} |f(s,0,0)| \,dC_s
\biggr)^{2prk} \biggr]^{{1}/({rk})}.\hspace*{-32pt}
\end{equation}

In Step 4, we finally want to take the expectation under the measure
$\PP$.
Let us define $\hat M_t=M_t - \int_0^t H_s \,d \langle M,M \rangle_s$
and note that since $H \cdot M$ is a BMO($\PP$) martingale, the process
$H \cdot\hat M$ and hence, $-H \cdot\hat M$ are BMO($\Q$)
martingales (see~\cite{Kazamaki}, Theorem 3.3).
Furthermore, by \cite{Kazamaki}, Theorem 3.1, there exists a $w, w'
>1$ such that $\mathcal E(H \cdot M)_T \in L^w(\PP)$ and $\mathcal E(-
H \cdot\hat M)_T \in L^{w'}(\Q)$.
As\vspace*{1pt} $\mathcal{E}(H \cdot M)^{-1}=\mathcal E (-H \cdot\hat M)$ we have
$d\PP=\mathcal E (-H \cdot\hat M)_T d\Q$.
Now, using the H\"{o}lder inequality with the conjugate exponent $v$ of
$w$ (and $v'$ of $w'$) and estimate (\ref{yestimate2}), we deduce
\begin{eqnarray*}
&&\E\Bigl[ \sup_{t \in[0,T]} |Y_t|^{2p} \Bigr]\\[1.2pt]
&&\qquad = \E^{\Q} \Bigl[ \mathcal E (-H \cdot\hat M)_T \sup_{t \in
[0,T]} |Y_t|^{2p} \Bigr]
\\[1.2pt]
&&\qquad \leq \E^{\Q} [ \mathcal E (-H \cdot\hat M)_T^{w'}
]^{{1/w'}} \E^{\Q} \Bigl[ \sup_{t \in[0,T]} |Y_t|^{2pv'}
\Bigr]^{{1/v'}}
\\[1.2pt]
&&\qquad \leq c \biggl( \E^{\Q} \biggl[ |B|^{2pv'rk} + \biggl(\int_0^{T}
|f(s,0,0)| \,dC_s \biggr)^{2pv'rk} \biggr]^{{1}/({rk})}
\biggr)^{{1/v'}}
\\
&&\qquad \leq c \E[ \mathcal E(H \cdot M)^w ]^{{1/w}} \\
&&\qquad\quad{}\times\E
\biggl[ |B|^{2pvv'rk} + \biggl(\int_0^{T} |f(s,0,0)| \,dC_s
\biggr)^{2pvv'rk} \biggr]^{{1}/({rkvv'})} .
\end{eqnarray*}
Setting $q=vv'rk$ and treating estimate (\ref{zestimate}) similarly
gives the desired result.
\end{pf}

%s7.3 ###
\subsection{Additional material on Markov processes}
\label{Appendix2}

We now provide an example where the function $u$ in Theorem
\ref{theorem:Markovproperty} does not depend trivially on~$M$.

Let $M:=(M_t)_{t\in[0,T]}$ be a continuous martingale with
nonindependent increments that is also a Markov process with respect to
a filtration $(\mathcal{F}_t)_{t\in[0,T]}$. Let
$X:=(X_t)_{t\in[0,T]}$ be the solution of the SDE
\[
dX_t=\int_0^t \sigma(a, X_a) \,dM_a, \qquad  t \in[0,T] \mbox{ and } X_0=0,
\]
with
\[
\sigma(a,x)=\cases{ 1+x, &\quad if $\displaystyle a \geq\frac{T}{2}$, \vspace*{3pt}\cr
0, &\quad if $\displaystyle a < \frac{T}{2}$.}
\]
Note that the coefficient $\sigma$ is Lipschitz in $x$ for every $a$
and that it is right continuous with left limits in $a$ for every
$x$; as a consequence $X$ admits an unique solution by \cite{Protter},
Theorem V.35. Consider a simple BSDE of the form~(\ref{BSDEbis}) with $f \equiv0$, $\kappa=0$ and $F(x):=\log(1+x)$.
Note that $F$ is not bounded but in this special case the existence
of a solution to the BSDE may be constructed directly.
Our aim is to show that $\E[F(X_T) \vert\mathcal{F}_t]$ is not a~trivial function of $M$
for $t \in[0,T]$. By It\^{o}'s formula we have
\[
F(X_T)=\log(1+X_t) + \int_t^T (1+X_s)^{-1} \,dX_s - \frac12 \bigl(\langle
M, M \rangle_T - \langle M, M \rangle_{t \vee({T}/{2})}\bigr),
\]
and hence,
\begin{eqnarray*}
\E[F(X_T)\vert\mathcal{F}_t]
&=&\log(1+X_t) -\frac12 \E\bigl[\langle M, M \rangle_T-\langle M, M
\rangle_{t \vee({T}/{2})}\vert\mathcal{F}_t\bigr] \\
&=&\log(1+X_t) - \frac12 \E\bigl[M_T^2-M_{t \vee({T}/{2})}^2\vert
\mathcal{F}_t\bigr]\\
&=&\log(1+X_t) - \frac12 \E\bigl[M_T^2-M_{t \vee({T}/{2})}^2\vert M_t\bigr]
\end{eqnarray*}
since $M$ is a Markov process.
Choose $0< t <\frac{T}{2}$ and then by definition of~$X$, the last
term on the right-hand side above cannot be
expressed as a~trivial deterministic function of $(t,X_t)$ since
$M_s$ cannot be deduced from $X_s$ for $s<\frac{T}{2}$. However,\vspace*{1pt} this
term is deterministic and only depends on $t$ if $M$ has independent
increments. This gives an example of a situation where the function
$u$ in Theorem \ref{theorem:Markovproperty} does not depend
trivially on $M$.
\end{appendix}

\section*{Acknowledgments}
The authors are grateful to
anonymous referees and to an Associate Editor for interesting comments and
remarks which have improved this paper. We also thank Stefan Ankirchner and
Gregor Heyne for inspiring
discussions and comments.

%suskaldyti doi

% imsref loaded by lrinkeviciute, 2011-03-25 08:41:04
%
% imsref loaded by lrinkeviciute, 2011-03-28 10:17:57

%
\printaddresses

\end{document}